\documentclass[12points]{amsart}

\usepackage{amsmath,amssymb,amsfonts}
\usepackage{epsfig,graphics,color}
\usepackage{amscd}
\usepackage[matrix,arrow,curve]{xy}

\newcommand{\DP}[2]{{\frac{\partial #1}{\partial #2}}} 
\newcommand{\p}{\partial}

\newcommand{\tx}[1]{\text{\rm #1}}
 
\newcommand{\CC}{\mathbb{C}}   
\newcommand{\NN}{\mathbb{N}}  
   
\newcommand{\RR}{\mathbb{R}}   
\newcommand{\Ss}{\mathbb{S}} 
 
\newcommand{\ZZ}{\mathbb{Z}}      
\newcommand{\PP}{\mathbb{P}}

\newcommand{\demo}{\noindent {\it \small Proof. }} 
 
\newtheorem{defi}{Definition}[section]    
\newtheorem{thm}[defi]{Theorem}    
  
\newtheorem{prop}[defi]{Proposition} 
\newtheorem{lem}[defi]{Lemma}

\newtheorem{ex}[defi]{Example} 
\newtheorem{rmk}[defi]{Remark}

\def\mc{\mathcal} 
\def\got{\mathfrak}

\title[Symplectic cohomology and the Leray-Serre spectral sequence]
{Fibered Symplectic cohomology\\and the Leray-Serre spectral sequence}
   
\author{Alexandru Oancea}

\date{July 24, 2007.}

\begin{document}

\maketitle

\vspace{-.5cm}

\begin{center} 
{\it Universit\'e Louis Pasteur, IRMA, F-67084 Strasbourg, France}

E-mail: {\tt oancea@math.u-strasbg.fr}
\end{center}

\begin{abstract}

We define Symplectic 
cohomology groups $FH^*_{[a,b]}(E)$, $-\infty \le a < b \le \infty$ 
for a class of symplectic
fibrations $F\hookrightarrow E \longrightarrow B$ 
with closed symplectic base and convex at infinity
fiber. The crucial geometric assumption on the fibration is a
negativity property reminiscent of negative curvature in complex
vector bundles.    
When $B$ is symplectically aspherical we construct 
a spectral sequence of Leray-Serre type converging to
$FH^*_{[a,b]}(E)$, and we use it  
to prove new cases of the Weinstein conjecture.

\end{abstract} 

\medskip 

{\footnotesize

}


\renewcommand{\thefootnote}{\arabic{footnote}}
\setcounter{footnote}{0}

\setcounter{tocdepth}{1}
\tableofcontents


\section{Introduction}

This paper investigates the behaviour of Symplectic
cohomology in a fibered context. The foundational papers on Symplectic
(co)homology are~\cite{FH94,CFH95,Viterbo99}, and recent developments
are presented in~\cite{Seidel07}. The symplectic fibrations
$F\hookrightarrow E \stackrel \pi \longrightarrow B$ that we consider
have a closed symplectic base $(B,\beta)$, a fiber with contact type
boundary and possess a coupling form $\Omega$. In this sense, they are
very close to being Hamiltonian~\cite{McDSal-Intro}. Moreover, they must
satisfy a negativity
property reminiscent of negative curvature in Hermitian vector
bundles, as well as a monodromy condition allowing to define
parallel transport along the boundary. 
We call them in the sequel {\it negative
  symplectic fibrations}~(see Definition~\ref{defi:main}). 
The total space $E$ is a symplectic manifold with symplectic form
$\omega_\epsilon=\pi^*\beta + \epsilon\Omega$, $\epsilon>0$ but the
boundary is not necessarily of contact type. It only satisfies a
fiberwise convexity condition which is nevertheless sufficient for us
to define Symplectic homology and cohomology groups $FH_*^{[a,b]}(E)$
and $FH^*_{[a,b]}(E)$, $-\infty\le a < b \le +\infty$. Our definition
is a fibered extension of the 
one of Viterbo~\cite{Viterbo99} -- which corresponds to the case
$B=\{\textrm{pt.}\}$ -- and features the same functorial properties. 

We call a symplectic manifold $(B,  \beta)$  {\it
  symplectically aspherical} if we have $\int f^*\beta =0$ for any
  smooth map
$f:\Ss^2 \longrightarrow B$. We say that $(B,  \beta)$ is {\it
  monotone} if there exists $\lambda \ge 0$ such that 
$\langle [\beta],  [f] \rangle 
= \lambda \langle c_1(TB),   [f] \rangle$ for any such $f$.
Our main result is the following. 

\medskip 

\noindent {\bf Theorem A.} {\it Assume 
$B$ is symplectically aspherical and $E$ is monotone. 
For any field of coefficients and any choice of real numbers 
$-\infty \le a < b \le \infty$ there 
exists a cohomology spectral sequence $E_r^{p,q}(a,  b)
\Longrightarrow FH^*_{[a,  b]}(E)$, $r\ge 2$ such that 
$$E_2^{p,q}(a,  b) \simeq H^{n+p}(B;  \mc{FH}^q_{[a,  b]}(F)), 
\ n=\frac 1 2 \dim \, B .$$
The spectral sequence is canonical and the notation $\mc {FH}^*_{[a, 
  b]}(F)$ stands for a local system of coefficients with
fiber $FH^*_{[a,  b]}(F)$. 
The spectral sequence is compatible with the truncation morphisms. 
}

\medskip 

The ``truncation'' morphisms in the statement are canonical morphisms
$$
FH^*_{[a,  b]}(E)\longrightarrow FH^*_{[a',  b']}(E), \quad a\ge a',
b\ge b'
$$
induced by the truncation of the range of the Hamiltonian
action. 

Let
$$\nu=\tx{min} \ \big\{ \ |\langle
f^*c_1(TE),   [\Ss^2]\rangle | \ : \ f:\Ss^2\longrightarrow E \ \big\}$$ 
be the minimal Chern number of $E$. The first Chern class of $TE$ is
computed with respect to an almost complex structure which is
compatible with the symplectic form $\omega_\epsilon$. 
If $\nu\neq 0$ the Floer cohomology groups are only $\ZZ/2\nu\ZZ$-graded and
the statement should be understood modulo $2\nu$. We shall not mention
anymore this grading issue in the course of the paper because it
is of a purely formal nature and has no bearing 
on the flow of the arguments.

The cohomology spectral sequence is constructed only with 
field coefficients because the
Symplectic cohomology groups are defined as an inverse limit. We need
in the proof of Theorem~A that the inverse limit functor be exact, which
is true if the terms involved in the limit are finite dimensional vector
spaces. On the other hand, the dual homology spectral sequence exists 
with arbitrary
coefficients because the Symplectic homology groups are defined as a
direct limit, which is an exact functor (see
Remark~\ref{rmk:homology}). 
Unless otherwise mentioned we
will use from now on field coefficients whenever 
Symplectic {\it cohomology} groups are involved.

\medskip 

Let us now introduce the following definition. Given a symplectically
aspherical manifold $(M,\omega)$ with boundary of contact type, we say
that $\p M$ is of {\it positive contact type} if every positively
oriented closed contractible characteristic $\gamma$ has positive
action $A_\omega(\gamma):= \int_{D^2} \bar \gamma^* \omega$ bounded away
from zero. Here $\bar \gamma:D^2\to M$ is any smooth extension of
$\gamma:\Ss^1\to M$ (see Section~\ref{sec:properties of
  Floer homology} for a discussion of this notion). Our main class of
examples is that of convex exact symplectic manifolds.  

  Let us recall that the homological structure of the fibration
  $(E,\pi,B,F)$ is captured by the Leray-Serre spectral sequence
  ${_{_{LS}}}E _{r}^{p,q}\Rightarrow H^{p+q}(E,\p E)$ with
  ${_{_{LS}}}E _2^{p,q} \simeq H^p(B;\mc{H}^q(F,\p F))$, where
  $\mc{H}^q(F,\p F)$ is the local system of coefficients on $B$ given
  by the locally constant presheaf $U\mapsto H^q(\pi^{-1}(U), \pi^{-1}(U)
  \cap \p E)$. 

\medskip 

\noindent {\bf Theorem~B.} {\it If $B$, $E$ are symplectically aspherical
  and $\partial F$ is of positive contact type in $F$ then, 
 for $\mu > 0$ small enough and $a<0$ arbitrary, the 
  spectral sequence $E_r^{p,q}(a,   \mu)$, $r\ge 2$ 
is canonically isomorphic to the Leray-Serre spectral sequence ${_{_{LS}}}E
 _{r}^{n+p,k+q}$, $n=\frac 1 2 \dim\, B$, $k=\frac 1 2 \dim\, F$. In
 particular, the local system $\mc {FH} 
 ^*_{[a, \mu]} (F) $ is canonically isomorphic to the
cohomological local system  $\mc{H}^{k+*}(F,   \partial F)$. 
}

\medskip 

Theorem~B can be read as a description of the Leray-Serre spectral sequence in
Morse homological terms. A related construction is that of
Hutchings~\cite{Hutchings03}, with the notable difference that he
views the fibration as being a \emph{family} of fibers, while we view
it as being an \emph{additional structure} on the total space. The
resulting generalizations and applications to Floer homology are very
different.  

The above comparison result yields applications to the Weinstein
conjecture. Under the hypothesis of Theorem~B we have  $FH^*_{]-\infty, 
  \mu]}(E)\simeq H^*(E,  \partial E)$ and there is a canonical
truncation morphism 
$$
c^*:FH^*(E)\longrightarrow H^*(E,  \partial E),
$$ 
where
$FH^*(E):=FH^*_{]-\infty,  +\infty[}(E)$. 
The functorial properties of the Symplectic cohomology groups $FH^*_{[a,  b]}(E)$ 
which are summarized in Section~\ref{sec:properties of Floer homology}
force dynamical consequences from algebraic assumptions. The resulting
principle concerning the Weinstein conjecture~\cite{Weinstein79} 
is the following. 

\medskip 

\noindent {\bf Main principle.} 
{\it If the morphism $c^*$ is not surjective in maximal degree then
  any contact type hypersurface bounding a compact domain in $E$
  carries a closed characteristic, i.e. the Weinstein conjecture holds
  in $E$.   
}

\medskip 

In the case $B=\{\tx{pt.}\}$ this was proved by
Viterbo~\cite{Viterbo99}. His proof carries over verbatim to our
situation once the groups $FH^*(E)$ have been defined and their
functorial properties have been established.
Following~\cite{Viterbo99} 
we call the above condition on $c^*$ the 
{\it Strong Algebraic Weinstein Conjecture (SAWC)} property. 

\medskip 

\noindent {\bf Remark.} The morphism $c^*$ is defined only if the 
manifold is aspherical and has positive contact type boundary. 
Thus the hypotheses of Theorem B are the most general ones under which 
one can apply the Main Principle stated above. 

\medskip 

\noindent {\bf Theorem C.} {\it Assume $E$, $B$ are symplectically
aspherical and $F$ has positive contact type boundary. 
The SAWC property is inherited from the
fiber by the total space. In particular, if $F$ satisfies the
SAWC then the Weinstein conjecture holds in $E$. 
}

\medskip 

\demo Under the assumption that the morphism 
$$FH^*(F) \longrightarrow
H^*(F,  \partial F)$$ 
is not surjective in degree $2k = \dim \, F$, we
have to prove that the morphism $FH^*(E)\longrightarrow H^*(E,  
\partial E)$ is not surjective in degree $2n+2k = \dim \, E$. Because
$H^{2n+2k}(E,   \partial E) \simeq H^{2n}(B;   \mc H ^{2k}(F,  
\partial F))$ and the morphism of spectral sequences $E_*^{*,*}
\longrightarrow {_{_{LS}}}E _*^{*,*}$ respects the bigrading, it is
  enough to show that the map 
$$H^{2n}(B;   \mc {FH}^{k}(F))
  \longrightarrow H^{2n}(B;   \mc H^{2k}(F,   \partial F))$$ 
is not
  surjective. We apply Poincar\'e duality on $B$ and we are left to
  show that the map $H_0(B;   \mc {FH}^{k}(F))
  \longrightarrow H_0(B;   \mc H^{2k}(F,   \partial F))$ is not
  surjective. We remark now that parallel transport in $E$ 
  is symplectic and therefore preserves
  the orientation of the fibers. This implies that the local system
  $\mc H^{2k}(F,   \partial F)$ is trivial and therefore $H_0(B;  
  \mc H^{2k}(F,   \partial F) \simeq H^{2k}(F,   \partial F)$. On
  the other hand $H_0(B;   \mc {FH}^{k}(F))$ is isomorphic to a
  quotient of $FH^{k}(F)$ (more precisely the quotient by the
  submodule $F'$ generated by elements of the form $\Phi_\alpha(u) - u$,
  where 
  $u\in FH^{k}(F)$, $\alpha \in \pi_1(B)$ and $\Phi_\alpha$ is the
  monodromy transformation along $\alpha$ -- see
  Section~\ref{sec:local systems} for details on homology with
  values in a local system). Now the hypothesis implies
  that the map $FH^{k}(F)/F' \longrightarrow H^{2k}(F,   \partial
  F)$, induced by $FH^{k}(F) \longrightarrow H^{2k}(F,   \partial
  F)$, is not surjective. 
\hfill{$\square$}

\medskip 

The most important class of symplectic manifolds which satisfy the
SAWC condition are subcritical Stein domains. In this case we actually
know by work of Cieliebak~\cite{Ci02} 
that Floer (co)homology is zero. We
obtain in particular the following vanishing theorem. 

\medskip 

\noindent {\bf Theorem D.} {\it Let $E$ be a negative symplectic
  fibration with symplectically aspherical base 
and subcritical Stein fiber. We have 
$$FH^*(E) =  0$$
and the Weinstein conjecture holds in $E$. 

In particular, the Weinstein conjecture holds if $E$ is the unit disc bundle of a 
Hermitian vector bundle with negative curvature over a symplectically
aspherical base. 
}

\medskip 

\demo By~\cite{Ci02} we know that, if the fiber $F$ is subcritical Stein,
we have $FH^*(F)=0$. The local system
giving the $E_2$-term of the spectral sequence in Theorem~A has
trivial fiber and therefore $E_2=0$, $E_\infty = 0$ and
$FH^*(E)=0$. 
\hfill{$\square$} 

\medskip 

\noindent {\bf Remark.} The homological analogue of the conclusion in
Theorem~D is that $FH_*(E)=0$. If $E\to B$ is a negative
  Hermitian disc bundle as in Definition~\ref{defi:neg line bundles}
  we deduce $FH_*^{[\mu,\infty[}(E)\simeq H_{*+n}(E,\partial
  E)\simeq H_{*+n-2}(B)$ for $n=\frac 1 2 \dim\, B$ and $\mu>0$ small
  enough. We have used here the tautological long exact sequence
  $FH_*(E)\to FH_*^{[\mu,\infty[}(E)\to H_{*+(n+1)-1}(E,\partial E)\to
  FH_{*-1}(E)$ from~\cite{Viterbo99}. This
  fits perfectly into the long exact sequence from~\cite{BOcont} which
  relates $FH_*^{[\mu,\infty[}(E)$ and linearized contact homology
  $HC_*(\partial E)\simeq \bigoplus _{k\ge 0} H_{*-2k}(B)$.  

\medskip 

In the case of a trivial fibration the spectral sequence degenerates
at $E_2$ by construction, since the Floer complex on $B\times
F$ can be identified with the tensor product of a Morse complex on $B$
with a Floer complex on $F$. The local system $\mc{FH}^*(F)$ is trivial 
and we obtain the K\"unneth formula (see also~\cite{Kunneth} for a
related statement). 

\medskip 

\noindent {\bf Theorem E.} {\it Let $B$ and $F$ be symplectically
aspherical. We have 
\begin{center}
$FH^*(B\times F)\simeq H^*(B;   FH^*(F)).$ 
\end{center} 
\hfill{$\square$}
}

\medskip 

\noindent {\bf Remark.} This result generalizes the one by Hofer and
Viterbo~\cite{FHV90}, who prove the Weinstein conjecture for
$E=B\times \CC^\ell$, $\ell \ge 1$. We have $FH^*(B\times
F)=0$ for any subcritical Stein manifold $F$, and the Weinstein
conjecture holds in $B\times F$ by the Main Principle stated above.

\medskip 

The construction of the spectral sequence is geometric. We choose on
$B$ a $C^2$-small Morse function $f$ and a generic 
almost complex structure $J_B$. We work on $E$ with Hamiltonians $K$
whose $1$-periodic orbits are nondegenerate and concentrated in the
fibers lying over the critical points of $f$. 
We consider on $E$ a modification of the standard Floer
equation having the form 
\begin{equation} \label{eq:main eq in the introduction} 
u_s+Ju_t = Y\circ u,
\end{equation} 
where $J$ is such that the projection $\pi$ is $(J, 
J_B)$-holomorphic and $\pi_*Y=\nabla f$. This 
ensures that the map $v=\pi\circ u$ will satisfy the equation 
\begin{equation} \label{eq:main eq on the base in the introduction}
v_s+J_Bv_t=\nabla f \circ v .
\end{equation} 
The main point is that we can choose many such pairs $(K,  Y)$ such
that the vector field 
\begin{equation} \label{eq:cYintro}
\mc Y(x) = J\dot x - Y\circ x
\end{equation} 
defined on the space of contractible
loops in $E$ is a strong pseudo gradient for the action functional $A_K$, i.e. 
$$dA_K
\cdot \mc Y \ge \alpha \parallel \mc Y \parallel _{L^2}^2$$
for some $\alpha>0$. This ensures that Floer cohomology can be
computed by studying moduli spaces of solutions of~(\ref{eq:main eq in
  the introduction}) (see Section~\ref{sec:defiSH}). 
Since the solutions of our pseudo-gradient Floer
equation~(\ref{eq:main eq in the introduction}) on the total space $E$
project to
solutions of the Floer equation~(\ref{eq:main eq on the base in the
  introduction}) on the base $B$, the Floer complex on $E$ can be
filtered by the index of the critical points of $f$. 
The resulting spectral sequence is the one in Theorem~A. 

The proof of Theorem~B uses that Morse homology is equal to
cellular homology as defined in~\cite[Appendix~A.4]{MS}, provided that
the unstable manifolds of the pseudo-gradient vector field give rise
to a CW-decomposition of the underlying manifold (see
Section~\ref{sec:MHlocsyst}). This last fact was shown to be true by
Laudenbach~\cite{Laudenbach92} under the mild assumption that the
pseudo-gradient vector field is equal near its zeroes to the 
gradient of a quadratic form with respect to the Euclidean metric. We
use Laudenbach's result also in order to construct the
\emph{Floer local system} $\mc{FH}^*_{[a,b]}(F)$.  

\medskip 

\noindent {\bf Remark.} 
In the case of a monotone basis our method of construction of the
spectral sequence runs into two kinds of difficulties. 
The first one is technical and concerns the proof of the
pseudo-gradient property, which involves a time-independent metric on
the base (see Remark~\ref{rmk:why symp asph is important}). The second
one is conceptual and concerns the expression of the $E_2$-term, which
has to encode quantum homological contributions from the base.

\medskip

The paper is structured as follows. We give in Section~\ref{sec:defi of
  sympl fib} the definition and first properties of negative symplectic
fibrations. Section~\ref{sec:examples} contains examples: trivial
fibrations, fibrations associated to loops of
compactly supported Hamiltonian diffeomorphisms of the fiber, 
negative vector bundles, 
convex fibrations as defined in~\cite[\S2.10]{CGMS02}.

We give in Section~\ref{sec:construction of FH} the
construction of Symplectic cohomology groups in a fibered setting (the
dual homological construction is sketched
in~Remark~\ref{rmk:homology}). 
The main difficulty of the construction is
the proof of a priori $C^0$-bounds on Floer trajectories. 
We emphasize the   following two distinctive features of our approach. 
\begin{enumerate} 
\item We allow admissible Hamiltonians to be ``asymptotically
  linear''. This is a much larger class than the ones previously
  considered in Symplectic homology
  constructions for arbitrary convex manifolds~\cite{CFH95,Viterbo99},
  and generalizes the class of asymptotically quadratic Hamiltonians
  in~\cite{FH94}.  
\item We consider a generalization of Floer's equation in which the
  zero-order term is modified so that the resulting vector
  field~\eqref{eq:cYintro} is a pseudo-gradient for the action
  functional.   
\end{enumerate} 
We need both these degrees of freedom in order to ensure that solutions
of the Floer equation on $E$ project on solutions of the Floer
equation on $B$. 

Sections~\ref{sec:main section pseudo
  gr vector fields} to~\ref{sec:construction of the sp seq} deal with
  the construction of the spectral sequence. We construct in
  Section~\ref{sec:main section pseudo gr vector fields}
  pseudo-gradient vector fields of a special form on the loop space of
  $E$, and we establish in Section~\ref{sec:main section transv}
  transversality within a geometrically meaningful class of almost
  complex structures. These technical ingredients are put together in
  Section~\ref{sec:construction of the sp seq}. The proof of Theorems~A
  and~B is given in subsection~\ref{sec:proofs}. Subsection~\ref{sec:local
  systems}, in which we explain how 
  local systems of coefficients can be encoded in the Morse complex,
  may be of independent interest. 

Appendix~\ref{sec:appendix symplectic forms} contains a proof of the
(purely linear) fact that a symplectic form $\omega$ is determined, in
its conformal class, by the set of $\omega$-compatible almost complex
structures. This is referred to in Sections~\ref{sec:defi of sympl
  fib} and~\ref{sec:neg line bdles}.


\section{Negative symplectic fibrations} \label{sec:defi of sympl fib}

\begin{defi}\label{defi:main}
  A locally trivial fibration $F \hookrightarrow E \stackrel \pi 
  \longrightarrow B$ is
  called a 
  {\rm   negative symplectic fibration with contact type boundary
  fibers} 
  (or, for short, {\rm negative symplectic fibration}) 
  if the following conditions are satisfied. 
  \begin{enumerate} 
    \item the base $B$ is closed and the fiber $F$ has a non-empty boundary. 
    \item there exists a $2$-form $\Omega \in \Omega^2(E,
      \RR)$ and a vertical vector field $Z$ defined in a
    neighbourhood of $\partial E$ such that: 
     \begin{itemize} 
      \item {\sc (symplectic fibration)} 
      $\Omega$ is nondegenerate along the fibers and globally closed;
      \item {\sc (contact type boundary)} 
       $Z$ is outward pointing and transverse to $\partial E$,
     and satisfies $L_Z\Omega=\Omega$;
      \item {\sc (monodromy)} the
      horizontal distribution $H=(\ker \, \pi_*)^{\perp_\Omega}$ is
      tangent to $\partial E$.
     \end{itemize}  
    \item {\sc (negativity)} there is a symplectic form $\beta$ on $B$
    and a nonempty open subset $\mc J' \subset \mc J(B,  
    \beta)$ such that, in a neighbourhood of $\partial E$, we have  
      $$\Omega(v,   \widetilde J_B v) \ge 0$$
    for any $v\in H$ 
    and any almost complex structure $J_B \in \mc J'$. 
    Here $\widetilde J_B$ denotes the lift of $J_B$ to $H$ and $\mc
    J(B,  \beta)$ is the set of almost complex structures $J_B$ which are
    compatible with $\beta$ in the sense that $\beta(\cdot, 
    J_B\cdot)$ defines a Riemannian metric. 
  \end{enumerate} 
\end{defi}

\noindent {\bf Notations and terminology.} 
The fiber $\pi^{-1}(b)$ at
a point $b\in B$ will be denoted either by $E_b$ or by $F_b$. We shall
refer to $H = 
(\ker \, \pi_*)^{\perp _\Omega}$ as the {\it horizontal
  distribution} or the {\it horizontal connection}, while $V=\ker \, 
\pi_*$ will be called the {\it vertical distribution}. 
  As $H$ is tangent to $\partial E$ we have a well-defined {\it parallel
  transport} 
$$\tau_\gamma : E_{\gamma(0)} \stackrel \sim
\longrightarrow E_{\gamma(1)}$$ 
associated to any continuous 
 path $\gamma :[0,   1] \longrightarrow  B$. The form
 $\Omega$ will be called the {\it connection form}, while the vector field $Z$
 will be called the {\it Liouville vector field} or the {\it Liouville
   vector field in the fibers}. We shall refer to a negative
 symplectic fibration as being a tuple $(E,   \pi,   B,   
  F,   \Omega,    Z,   \beta)$ or, in order to emphasize the role
of $\Omega$ and $Z$, we shall simply refer to it as a  triple
 $(E,   \Omega,   Z)$.

\medskip 

\noindent {\bf Remarks.} 
1. The hypothesis $d\Omega =0$
 implies in particular $d\Omega (v_1,   v_2,   \cdot) =0$ for any
 $v_1,   v_2 \in V$. This in turn is equivalent to the
 fact that $\tau_\gamma$ is a symplectic diffeomorphism, where the
 symplectic form on $E_b$ is
 $\Omega_b=\Omega|_{E_b}$~\cite[Lemma~6.11]{McDSal-Intro}. 

\medskip 

2. The conditions $L_Z\Omega = \Omega$ and $d\Omega =0$ imply that
   $\Omega $ is exact in a neighbourhood of $\partial E$. The
   primitive is 
   $$\Theta =  \iota_Z\Omega .$$
   Moreover, the assumption that $Z$ is vertical implies
   $$\Theta |_H \equiv 0 .$$
We denote by
$\tx{Spec}(\partial E)$ the set of periods of closed characteristics
on $\partial E$ normalized by the $1$-form $\Theta=\iota_Z\Omega$. 

\medskip    

3. The total space of a negative symplectic fibration is itself a
   symplectic manifold, with symplectic form 
   $$\omega_\epsilon = \pi^*\beta + \epsilon \Omega, \quad \epsilon >
   0 \tx{ small enough.}$$
   
   \medskip 

4. The boundary $\partial E$ may fail to be of contact
   type, as we {\it do not} suppose that $\pi^*\beta$ is exact in one
   of its neighbourhoods. This phenomenon happens for example in
   trivial fibrations $E=B\times F$. Nevertheless, as we shall see in
   \S\ref{sec:proof of C0 bounds}, 
a version of weak pseudoconvexity still holds and that
   will be enough in order to make use of the maximum principle.

\medskip 

5. Let $\varphi_t$ be the flow of $Z$. 
We can trivialize a neighbourhood $\mc U$ of $\partial E$
by the diffeomorphism 
$$\Psi :  \partial E \times [1-\delta,    1]
  \longrightarrow \mc U  ,$$
$$(p,   S) \longmapsto \varphi_{\ln S}(p) .$$
The condition $L_Z\Omega=\Omega$ translates into $\varphi_t^*\Omega =
e^t \Omega$. If $\Theta|$ denotes the restriction of $\Theta$ to $\partial E$ 
we  have $\Psi^* \Theta = S\Theta|$ and 
$$\Psi^*\Omega = d(S\Theta|) .$$
We can therefore complete $E$ to a fibration 
$$\widehat E = E   \bigcup_{\Psi}   \partial E \times [1,   \infty[
$$
and define the connection form $\widehat \Omega$ on $\widehat E$ by 
$$ \widehat \Omega = \left\{ \begin{array}{ll} \Omega & \tx{ on } E \
      ,  \\
d(S\Theta|) & \tx{ on } \partial E \times [1,   \infty[ \
      . \end{array} \right.$$
The Liouville vector field $Z$ is transformed by $\Psi$ into
$S\DP{}{S}$ on $\partial E \times [1-\delta,
  1]$. We extend it to $\partial E \times [1,   \infty[$ 
as $S\DP{}{S}$ and we denote the extended vector field by $\widehat Z$. 

\medskip 

The {\sc (monodromy)} condition implies
that the horizontal distribution on $\widehat E$ is tangent to every
level set $S=\tx{ct}$, $S\ge 1-\delta$. This follows from the fact that 
$\varphi_t$
preserves $H$, hence the latter is invariant under the flow of
$S\DP{}{S}$ (or, equivalently, of $\DP{}{S}$). We therefore have 
$$
dS|_H \equiv 0.
$$

\medskip 

6. One must note that the construction of the manifold 
$\widehat E$ only makes use
   of the {\sc (symplectic fibration)} and {\sc (contact type
   boundary)} conditions, while the {\sc (negativity)} condition
   ensures that the $2$-form $\widehat \omega_\epsilon =
   \pi^*\beta + \epsilon\widehat \Omega$ on $\widehat E$ is symplectic
   for $\epsilon>0$ small enough. 

   Conversely, let us start with a fibration $\widehat E$ endowed
   with a $2$-form $\widehat \Omega$ and a vertical 
   vector field $\widehat Z$ which is
   complete at infinity, satisfying the {\sc (symplectic
   fibration), (contact type boundary)} and {\sc (negativity)}
   conditions. A compact hypersurface $\Sigma \subset \widehat E$ 
   that trivializes through the flow of $\widehat Z$ a
   neighbourhood of infinity as $\Sigma \times [1,\infty[$ will be
   called a {\it trivializing hypersurface}. The choice of any such
   $\Sigma$ gives rise to a fibration $E=\overline{\tx{int}  
   \Sigma}$ which satisfies the same three conditions above. 

\medskip 

7. We have imposed the {\sc (monodromy)} condition on $E$ in order for the
   monodromy to be well defined as a symplectic diffeomorphism
   of the fiber. Note however that, if one starts directly with 
   $\widehat E$ as above, the natural condition under which monodromy
   is well-defined is some uniform non-verticality assumption on $H$, 
   strictly weaker than the requirement $dS|_H \equiv
   0$. Our choice is motivated by the fact that 
   trivializing hypersurfaces $\Sigma$ 
   such that $dS|_\Sigma \equiv 0$ are a 
   crucial ingredient in the proof of a 
   priori $C^0$-bounds for the admissible 
   Hamiltonians on $\widehat E$ that we define in
   \S\ref{sec:proof of C0 bounds}.  

\medskip 

   The condition $dS|_H \equiv 0$ along $\partial E$ is equivalent to
   the fact that the characteristic distribution on $\partial E$ is
   contained in the fibers. In particular it coincides with the
   characteristic distribution of the restriction of $\Omega$ to the
   fibers. The latter is 
   preserved by parallel transport and therefore the Reeb dynamics on 
   the boundary of the fibrations that we consider in this paper is of
   Morse-Bott type, with the meaning that one closed characteristic on
   $\partial E_z$, $z\in B$ 
   gives rise locally to a family of closed characteristics on 
   $\partial E$ parametrized by an open subset of the base. The
   reader has to keep in mind this 
   geometric picture as a motivation for the 
   construction of the geometric Hamiltonians in \S\ref{sec:the
  geometric Hamiltonians}.

The following result gives a geometric criterion for the {\sc
  (monodromy)} condition. 

\begin{prop} 
  Let $(\widehat E,   \widehat \Omega,   \widehat Z)$ be a fibration
  satisfying conditions {\sc (symplectic fibration)} and {\sc (contact type
  boundary)}, with $\widehat Z$ complete at 
  infinity. There is a choice of a trivializing hypersurface $\Sigma$
  in $\widehat E$ such that $H$ is tangent to $\Sigma$ 
  if and only if the monodromy of $\widehat E$ admits an
  invariant trivializing hypersurface in the fiber. 
 
\end{prop}

\demo Assume first that $H$ is tangent to $\Sigma$. For any loop 
$\gamma$ on $B$ based at $b$ the monodromy will send $\Sigma_b=\Sigma
\cap E_b$ diffeomorphically to itself and $\Sigma_b$ is obviously 
a trivializing hypersurface in $E_b$.

Conversely, let $\Sigma_b \subset E_b$ be an invariant trivializing
hypersurface. Define $\Sigma = \bigcup _{b'\in B}
\tau_{\gamma_{bb'}}(\Sigma_b)$, where $\gamma_{bb'}$ is an arbitrary
path from $b$ to $b'$. It is enough to prove that
$\tau_{\gamma_{bb'}}(\Sigma_b)$ is independent of $\gamma_{bb'}$ in
order to infer that $\Sigma$ is smooth and $H \subset T\Sigma$. Let
therefore $\widetilde \gamma_{bb'}$ be another path running from $b$
to $b'$. We have 
$$\tau_{\widetilde \gamma_{bb'}} (\Sigma_b) =  \tau_{\gamma_{bb'}}
\circ \tau_{\gamma_{bb'}^{-1}\cdot \widetilde \gamma_{bb'}} (\Sigma_b)
= \tau_{\gamma_{bb'}}(\Sigma_b) .$$
The last equality makes use of the fact that $\Sigma_b$ is monodromy
invariant. 
\hfill{$\square$}

\medskip 

8. One may strengthen the {\sc (negativity)} condition by
requiring $\Omega$ to be nonnegative on $H$ for 
{\it all} $\widetilde J_B$
where $J_B\in \mc J(B,   \beta)$. Appendix~\ref{sec:appendix
  symplectic forms} shows that this is
a very strong assumption: either $\Omega$ is nondegenerate on $H$ and
then it is proportional to $\pi^*\beta$, either $H$ splits at
each point into a direct sum of two subspaces which are symplectic for
$\pi^*\beta$ and such that $\Omega$ vanishes on one of them and is
proportional to $\pi^*\beta$ on the other. 

Nevertheless, requiring $\Omega$ to tame only the almost complex
structures belonging to some nonempty open subset of $\mc J (B,
  \beta)$ is enough for the subsequent transversality issues.
If the restriction of $\Omega$ to $H$ is nondegenerate on $\partial
E$, this follows by imposing the {\sc (negativity)} condition to
be true only for {\it one} almost complex structure $J_B$.


\section{Examples} \label{sec:examples}


\subsection{Products} Trivial fibrations 
$E=B\times F$ with $B$ closed symplectic and $F$
    symplectic with contact type boundary are negative symplectic
    fibrations in the sense of Definition~\ref{defi:main}. Let
    $\pi_B$, $\pi_F$ denote the projections on the two factors, let
    $\Omega_F$ be the symplectic form on $F$ and $Z_F$ be the Liouville
    vector field on $F$, defined in a neighbourhood of $\partial
    F$. Then $\Omega =\pi_F^*\Omega_F$ and $Z= (0,   Z_F) \in TB \times TF
    \simeq T(B\times F)$ make $E$ into a negative symplectic fibration. We
    have $H = TB \times \{ 0\} \subset T(B\times F)$ and $\Omega|_H
    \equiv 0$, so the negativity condition is trivially satisfied.


\subsection{Hamiltonian diffeomorphisms} Let $F\hookrightarrow E
\stackrel \pi \longrightarrow B$ be a Hamiltonian fibration with
contact type boundary fibers and 
structure group $\tx{Ham}(F,  \partial F)$, the group of Hamiltonian
diffeomorphisms which fix a neighbourhood of $\partial F$. A
Hamiltonian fibration admits a canonical coupling form $\Omega$, which
in our situation vanishes on the horizontal distribution near the
boundary (see~\cite{GLS96}) and makes $E$ into a negative symplectic
fibration. As a special case we mention fibrations $F\hookrightarrow E
\longrightarrow \Ss^2$ defined by elements of $\pi_1(\tx{Ham}(F, 
\partial F))$.


\subsection{Negative line bundles} \label{sec:neg line bdles}

\begin{defi} \label{defi:neg line bundles}  
A complex line bundle $\mc L \longrightarrow B$  over a closed 
symplectic manifold $(B,   \beta)$ is called {\rm negative} if it admits a
Hermitian metric $h$ and a Hermitian connection $\nabla$ such that the
curvature $\frac 1 {2i\pi} F^\nabla \in \Omega^2(B,   \RR)$ 
is negative: 
$$\frac 1 {2i\pi} F^\nabla (v,   J_Bv ) < 0$$
for any $J_B\in \mc J(B,  \beta)$ and any
nonzero vector $v\in TB$.  
\end{defi}  

\medskip 

\noindent {\bf Remark.} The $2$-form $- \frac 1
{2i\pi}F^\nabla$  is a symplectic form on $B$ representing $-c_1(\mc
L)$. Moreover, it tames all almost complex structures that are tamed by
$\beta$ and this ensures that $- \frac 1
{2i\pi}F^\nabla$ and $\beta$ are proportional by a positive constant
if $\dim \, B\ge 4$ or by a positive function if $\dim \, B=2$
(cf. Appendix~\ref{sec:appendix symplectic forms}). 
In particular we have $c_1(\mc L) = -\lambda [\beta]$,
$\lambda >0$ (if $\dim \, B =2$ this is true because 
$\dim \, H^2(B,   \RR) =1$). 
Conversely, assume $c_1(\mc L) = -\lambda
[\beta]$, $\lambda >0$. Then $-\lambda \beta$ represents $c_1(\mc L)$
and, for any Hermitian metric $h$ on $\mc L$, 
one can find a Hermitian connexion $\nabla$ such that $\frac 1 {2i\pi}
F^\nabla = -\lambda \beta$. In particular $-\frac 1 {2i\pi} F^\nabla$
tames the same almost complex structures as $\beta$. 
We have just proved that Definition~\ref{defi:neg line bundles} is
equivalent to

\begin{defi} \label{defi:neg line bundles 2}
A complex line bundle $\mc L \longrightarrow B$  over a closed 
symplectic manifold $(B,   \beta)$ is {\rm negative} if there exists
$\lambda >0$ such that 
$$c_1(\mc L) = -\lambda [\beta] .$$
\end{defi} 

We note here that the topological type of $\mc L$ is uniquely determined by
the choice of an integral lift of $-\lambda [\beta]$.  
The preceding discussion shows in particular that, 
up to a change of connection, we can assume that $-\frac 1
{2i\pi} F^\nabla = \lambda \beta$, $\lambda >0$.
 
\medskip 

Any linear connection $\nabla$ determines a  
  {\it transgression $1$-form} $\theta^\nabla \in 
  \Omega^1(\mathcal{L} \setminus 0_{\mathcal{L}}, \ \RR)$. Its
  definition is the following~\cite{Gauduchon94}:
 \begin{equation*}
        \left\{ \begin{array}{l} 
   \theta ^\nabla _u(u) = 0, \quad \theta ^\nabla _u (iu) =  
     1/2\pi,\quad u \in \mathcal{L}\setminus 0_\mathcal{L}\ ; \\   
    {} \\
   \theta ^\nabla | _{H ^\nabla} \equiv 0  , \tx{ with } H 
     ^\nabla \tx{ the horizontal distribution defining }
      \nabla . 
        \end{array} \right.   
  \end{equation*}   
  The transgression form is a primitive for  
 $-\pi^*(\frac 1 {2 i \pi} F^\nabla)$. In our case, this means
  $$ d\theta ^\nabla = \lambda \pi ^*\beta .$$
On the other hand the restriction of $\theta^\nabla$ to
the fibers equals, up to the factor $\frac 1 {2\pi}$, the angular
form. If $r(u)=|u|$ is the radial coordinate in the fibers we infer
that 
$$\Omega = d(r^2\theta ^\nabla) $$
equals $\frac 1 \pi   d\tx{Area}$ along the fibers. Moreover, 
$\Omega$ extends to a smooth form on $\mc L$ by $\Omega_z(\xi,
\cdot)=0$, $\xi \in T_z 0_{\mc L}$ and $\Omega_z|_{\mc L_z}=\frac 1
\pi   d\tx{Area}$, $z\in 0_{\mc L}$. This follows from the expansion
$\Omega = dr^2 \wedge \theta^\nabla+r^2\lambda \pi^*\beta$, with 
$dr^2\wedge \theta^\nabla_u(\xi, \cdot)=0$, $\xi \in H^\nabla_u$ and 
$dr^2\wedge \theta^\nabla_u|_{\mc L_{\pi(u)}} = \frac 1 \pi   d\tx{Area}$,
$u\in \mc L\setminus 0_{\mc L}$. The vertical vector field
$$Z(u)=\frac u 2$$ 
satisfies $\iota _Z\Omega = r^2\theta^\nabla$, hence $L_Z\Omega = \Omega$.
We define 
$$E = \{ u \in \mc L \ : \ |u| \le 1\} .$$
We claim that $E$ together with $\Omega$ and $Z$ as
above is a negative symplectic fibration in the sense of
Definition~\ref{defi:main}. The {\sc (symplectic fibration)} and {\sc
  (contact type boundary)} conditions are clear by construction. The
connection being Hermitian, parallel transport preserves the length of
vectors in $\mc L$ hence the {\sc (monodromy)} condition is also
satisfied. The {\sc (negativity)} condition follows from the expansion of
$\Omega$, which implies $\Omega|_{H^\nabla}  = r^2\lambda\pi^*\beta$
with $r^2\lambda \ge 0$. 

\medskip 

\noindent {\bf Remark.} The dual $\mc L^*$ of an ample line bundle
$\mc L$ over a complex manifold $B$ is a negative line bundle in the
sense of Definition~\ref{defi:neg line bundles 2}, with the meaning
that there is a symplectic form $\beta$ on $B$ such that $c_1(\mc L^*) =
-\lambda [\beta]$, $\lambda>0$. Indeed, Kodaira's embedding theorem
ensures the existence of an embedding $\phi_m:B \hookrightarrow
\PP^N$ given by the sections of $\mc L ^{\otimes m}$ such that
$\phi_m^*\mc O(1)= \mc L^{\otimes m}$. Let $\omega_{\tx{FS}}$ be the
Fubini-Study form on $\PP^N$ representing $c_1(\mc O(1))$, and define
$\beta = \phi_m^*\omega_{\tx{FS}}$. Then 
$$
c_1(\mc L^*)  = - \frac 1 m c_1(\mc L^{\otimes m})
= -\frac 1 m \phi_m^* c_1(\mc O(1)) 
= -\frac 1 m [\phi_m^*\omega_{\tx{FS}}] = -\frac 1 m [\beta] .
$$

\medskip 


\subsection{Negative vector bundles} 

Our discussion in this section follows Griffiths~\cite{Griffiths69} and
Kobayashi~\cite{Kobayashi87}. 

\begin{defi} \label{defi:neg vector bundles}
  A complex vector bundle $E \stackrel \pi \longrightarrow B$ over a
  closed 
  symplectic manifold $(B,   \beta)$ is called 
{\rm negative} if it admits a
  Hermitian metric $h$ and a Hermitian connection $\nabla$ such that
  the curvature $\frac 1 i F^\nabla \in \Omega^2(B,   \tx{End}   E)$
  is negative definite as a Hermitian matrix: 
  $$ \frac 1 i F^\nabla(v,   J_B v) < 0$$
  for any $J_B\in \mc J(B,  \beta)$ and any
  nonzero vector $v \in TB$. 
\end{defi} 

\medskip 

\noindent {\bf Remark.} {\it Negative} 
projectively flat bundles (with curvature $-\beta\tx{Id}$) are a
particular case of negative vector bundles.

\medskip 

The projectivized bundle associated to $E$ is 
$$\PP(E) = \big\{ (b,   [v]) \ : \ b\in B,   v \in E_b \setminus \{ 0 \}
\big\} .  $$
Let $p: \PP(E) \longrightarrow B$ and $\bar
p:p^*E \longrightarrow E$ be the induced projection 
and the corresponding bundle map. The 
tautological bundle $\mc L _E \longrightarrow \PP(E)$ is the
subbundle of $p^*E$ defined as 
$$ \mc L_E = \{ (b,   [v],   \lambda v) \ : \ (b,   [v])\in
\PP(E),   \lambda \in \CC \} .$$
Let $i: \mc L_E \longrightarrow p^*E$ be the canonical inclusion and
$\Phi  = \bar p \circ i$. Then
$$\Phi : \mc L_E \setminus 0_{\mc L_E} \stackrel \sim 
\longrightarrow E \setminus 0_E$$ 
is a
diffeomorphism (respectively a biholomorphism if the bundle $E
\longrightarrow B$ is holomorphic). Its inverse is 
$$\Psi : E \setminus 0_E \stackrel \sim \longrightarrow \mc L_E
\setminus 0_{\mc L_E}  ,$$
$$(b,   v) \longmapsto (b,   [v],   v) .$$

We show now that, under the negativity condition on $E$, the line
bundle $\mc L_E$ is negative in the sense that there is a canonical
connection $\nabla$ which preserves the induced
Hermitian metric on $\mc L_E$, as well as a canonical symplectic form
on $\PP(E)$, such that Definition~\ref{defi:neg line bundles} is
satisfied. Moreover, if $\Omega$ and $Z$ are
constructed on $\mc L_E$ as in the previous section, their pull-backs 
$\Omega_E=\Psi^*\Omega$ and $Z_E=\Psi^*Z$, defined a priori 
only on $E\setminus 0_E$, extend smoothly to $E$ and make 
$$ DE = \{ v \in E \ : \ |v | \le 1 \} $$
into a negative symplectic fibration. 
$$
\xymatrix
@R=13pt
{
E \ar[dd]_\pi \ar@/^1.3pc/ @{.>}[rr]^\Psi
\ar@/_1.9pc/ @{<-}[rr]_\Phi
 & p^*E \ar[l]_{\bar p} \ar@{-}[d]  & \mc L_E
\ar[l]_i  \ar[dd]_\pi \\
& \ar[d]_\pi & \\
B & \PP(E) \ar[l]_p & \PP(E) \ar@{=}[l] 
}
$$

\subsubsection{Connection on $\mc L_E$} We write 
$\PP(E)_b$ for the fiber $p^{-1}(b)$, $b\in
B$ and, for a connection $D$, we denote by $H^D$ 
the associated horizontal distribution, or simply $H$ if there is no
danger of confusion.  

The connection $\nabla$ canonically defines a parallel transport in
$\PP(E)$ and hence a horizontal
distribution $H_{\PP(E)}$ with monodromy in $PGL(r)$, $r=\tx{rk}(E)$. 
The horizontal distribution $H$ associated to any connection on $p^*E$ 
canonically decomposes as
$$ H = H_{\tx{fiber}} \oplus
H_{\tx{base}}  ,$$ 
where the subspaces $H_{\tx{fiber}}$ and 
$H_{\tx{base}}$ of $H$ are uniquely determined by
the conditions 
$$ p_*\circ \pi_* H_{\tx{fiber}} = 0, \qquad \pi_*H_{\tx{base}} = 
H_{\PP(E)} .$$
Note that $H_{\tx{fiber}} 
\subset T\big(p^*E|_{\PP(E)_b}\big)$ and $H_{\tx{base}} \pitchfork
T\big(p^*E|_{\PP(E)_b}\big)$, $b\in B$. We call them respectively 
{\it the components of
$H$ along the fibers of $\PP(E)$ and along the base $B$}. The
distribution $H_{\tx{fiber}}$ defines a linear connection on every
$p^*E|_{\PP(E)_b}$, $b\in B$.

We can give further details on 
the preceding decomposition for the induced connection
$p^*\nabla$. The associated horizontal distribution is 
$H^{p^*\nabla} = (\bar p_*)^{-1} H^\nabla$. 
Any choice of frame $(e_1, \ldots,   e_r)$ in $E_b$
gives rise to a trivialization 
$$ \PP(E)_b \times \CC^r \stackrel \sim \longrightarrow
p^*E|_{\PP(E)_b}  ,$$
$$(b,   [v],   (\lambda_1, \ldots,   \lambda_r)) \longmapsto (b,  
[v],   \lambda_1 e_1 +\ldots + \lambda_r e_r) .$$
The associated flat connection does not depend on the choice of frame and
its horizontal distribution is precisely
$H^{p^*\nabla}_{\tx{fiber}}$. 

\medskip 

Let us now go to $\mc L_E$. Its restriction to any $\PP(E)_b$ 
is clearly not preserved by parallel transport along
$H^{p^*\nabla}_{\tx{fiber}}$, otherwise $\mc L_E|_{\PP(E)_b}$ would be
trivial. On the other hand, $\mc L_E$ {\it is} preserved,
together with the induced Hermitian metric, by
parallel transport along $H^{p^*\nabla}_{\tx{base}}$. This follows
from the fact that $\pi_*H_{\tx{base}}=H_{\PP(E)}$. 

We want to associate to $h$ and $\nabla$ in a canonical way 
a connection $\widetilde \nabla$
on $\mc L_E$ which preserves the induced Hermitian
metric. Its horizontal distribution decomposes as $\widetilde H = \widetilde
H_{\tx{fiber}} \oplus \widetilde H_{\tx{base}}$ and we define 
$$ \widetilde H_{\tx{base}} = H^{p^*\nabla}_{\tx{base}} .$$
Let us define $\widetilde H_{\tx{fiber}}$. 
Each restriction $\mc L_E|_{\PP(E)_b}$, $b\in B$ is
isomorphic (as a Hermitian bundle) 
to the canonical bundle $\mc O(-1) \longrightarrow
\PP^{r-1}$ endowed with the canonical Hermitian metric. 
The isomorphism is given by the choice of some frame in
$E_b$ which is orthonormal with respect to $h$. The Chern connection
on $\mc O(-1)$ is invariant under the action of $PSU(r)$ and this
implies that the induced connection on $\mc L_E|_{\PP(E)_b}$ is
independent of the choice of orthonormal frame. We define its
horizontal distribution to be $\widetilde H_{\tx{fiber}}$. With a
slight abuse of notation, we can write the decomposition $\widetilde H
= \widetilde H_{\tx{fiber}} \oplus \widetilde H_{\tx{base}}$ as
$$ \widetilde H = H_{\mc O(-1)} \oplus H^{p^*\nabla}_{\tx{base}} .$$

\subsubsection{Symplectic form on $\PP(E)$} The curvature of the 
Chern connection on $\mc O(-1)$ is $-\omega_{\tx{FS}}$, with
$\omega_{\tx{FS}}$ the Fubini-Study form normalized by
$\langle [\omega_{\tx{FS}}],    [\CC P^1] \rangle = 1$. We infer that
$$\omega = -\frac 1 {2i\pi} F^{\widetilde \nabla}$$ 
is a $2$-form which
restricts to $\omega_{\tx{FS}}$ on every fiber $\PP(E)_b$, $b\in B$. 
We claim that $\omega$ is actually nondegenerate on $\PP(E)$. This will define
our preferred symplectic form on $\PP(E)$. 

Firstly we show that $T(\PP(E)_b)$ and $H_{\PP(E)}$ are orthogonal with
respect to $\omega$. This amounts to proving that $H_{\mc O(-1)}$ and
$H^{p^*\nabla}_{\tx{base}}$ are in involution, as the value of the 
curvature at two vectors is given by the vertical projection of the
Lie bracket of their horizontal lifts (see e.g. Gauduchon~\cite{Gauduchon94}). 

Let $u(s,  t)$, $s,   t \in [0,  1]$ be a parametrized surface on 
$\PP(E)$ such that 
\begin{itemize} 
 \item $u(\cdot,   0)$ is tangent to some $\PP(E)_b$, $b\in B$; 
 \item $u(\cdot,   t)$ is the parallel transport of $u(\cdot,   0)$
 along some curve $\gamma$ on $B$ with $\gamma(0)=b$.  
\end{itemize} 

Let us fix a point $q \in \mc L_{E,   u(0,   0)}$ and a
horizontal lift $\widetilde u$ of $u(\cdot,   0)$ at $q$. This allows to lift
horizontally every curve $u(s, \cdot)$ with initial point
$\widetilde u(s,   0)$. We still denote by $\widetilde u$ the
resulting lift of $u(\cdot, \cdot)$ and we have to show that every
$\widetilde u(\cdot,   t)$ is horizontal.

This amounts to show that $H_{\mc O(-1)}$ is preserved
by parallel transport along $H^{p^*\nabla}_{\tx{base}}$. But $H_{\mc
  O(-1)}$ corresponds via the isomorphism $\CC^n \setminus \{ 0 \}
\simeq \mc O(-1) \setminus 0_{\mc O(-1)}$ to the distribution of
hyperplanes $(\CC\cdot v)^\perp$, $v\in \CC^n \setminus \{ 0 \}$,
which is clearly preserved by Hermitian parallel transport in
$E$. The latter in turn corresponds to parallel transport along
$H^{p^*\nabla}_{\tx{base}}$ in $\mc L_E$.

Secondly we show that the negativity condition in 
Definition~\ref{defi:neg vector bundles} is equivalent to the fact 
that $\omega$ is 
positive on $H_{\PP(E)}$, with the meaning that 
$\omega(X,   \widetilde J_B X) >0$ for any nonzero vector $X\in
H_{\PP(E)}$, where $\widetilde J_B$ 
is the lift to $H_{\PP(E)}$ of an almost complex structure
$J_B$ compatible with $\beta$. We denote $X'=p_*X$ and we have
\begin{eqnarray*} 
- \frac 1 i F^{\widetilde \nabla}_{(b,   [v])} (X,  
\widetilde J_B X) & = & 
- \frac 1 {i|v|^2} 
{^t \bar v} \cdot F^{p^*\nabla}(X,   \widetilde J_B X) \cdot v  \\
& = & - \frac 1 {i|v|^2} 
{^t \bar v} \cdot F^\nabla(X',    J_B X') \cdot v  > 0 .
\end{eqnarray*}

This shows that $\mc L_E$ is a negative line bundle in the sense of
Definition~\ref{defi:neg line bundles}. If $\theta ^{\widetilde
  \nabla}$ is the transgression $1$-form associated to
$\widetilde \nabla$, then the connection $2$-form and the Liouville vector
field on $\mc L_E$ are 
$$ \Omega = d(r^2 \theta^{\widetilde \nabla}), \qquad Z(u) = \frac u 2
.$$

\subsubsection{Connection form and Liouville vector field on $E$}
We define 
$$\Omega_E = \Psi^*\Omega, \qquad Z_E=\Psi^*Z .$$
We claim that $\Omega_E$ and $Z_E$ extend smoothly to the whole of $E$
and they verify Definition~\ref{defi:neg vector bundles}. The key
step is to consider 
$$\Theta= \Psi^*\big( r^2 \theta^{\widetilde \nabla}\big) .$$ 
We clearly have $\Theta|_{E_b\setminus \{ 0 \}} = \Psi^*\big(
r^2\theta^{\widetilde \nabla} |_{\PP(E)_b}\big)$ and we claim that
this is the positive $U(r)$-invariant Liouville form on $E_b$. By
choosing a unitary frame on $E_b$ we can work within the explicit
model of the biholomorphism $\CC^r\setminus \{ 0 \} \simeq \mc O(-1)
\setminus 0_{\mc O(-1)}$, $v\longmapsto ([v],   v)$. We have already
mentioned that the horizontal distribution of the Chern connection on
$\mc O(-1)$ corresponds to the distribution of hyperplanes $(\CC\cdot
v)^\perp$, $v\in \CC^r \setminus \{ 0 \}$ and we therefore have 
$\Theta _v |_{(\CC \cdot v)^\perp} \equiv 0$, $\Theta_v(v) =0$,
$\Theta_v(iv) = \frac {|v|^2} {2\pi}$, or else stated 
\begin{equation*} 
\Theta_v  =  \frac 1 {2\pi} \langle iv, \cdot \rangle  
=  \frac 1
{2\pi} \sum_{j=1}^r x_jdy_j - y_j dx_j.
\end{equation*}
As a consequence, $\Theta$ extends smoothly over the origin in every
fiber. But it is clear that this argument can be performed
in families and the extension is smooth on $E$ so that $d\Theta$ is
a smooth extension of $\Omega_E$ which is closed and equal to $\frac
1 \pi d\tx{Area}$ in the fibers. It is also clear that $Z_E$ 
extends smoothly by $0$ over $0_E$, with $\iota_{Z_E}\Omega_E=\Theta$. 
This accounts for the
\textsc{(symplectic fibration)} and \textsc{(contact type boundary)}
conditions. The \textsc{(monodromy)} condition is automatic as the
connection $\nabla$ was supposed from the very beginning to be
Hermitian. In order to verify the 
\textsc{(negativity)} condition let us recall that 
$\Phi_*\widetilde H_{\tx{base}}= H^\nabla$. We then have 
$$\big(\Psi^*\Omega \big)|_{H^\nabla} = \Psi^*\big( \Omega|_{\widetilde
  H_{\tx{base}}} \big) = \Psi^* \big( \lambda r^2 (\pi^* \omega)
  |_{\widetilde H_{\tx{base}}} \big)  = \lambda r^2 \Psi^* \pi^*
  (\omega|_{H_{\PP(E)}})  .$$ 

Let $J_B\in \mc J(B,  \beta)$ and denote
$\widetilde J_B$ the lift to $H^\nabla$ and $J'_B$ the lift to
$H_{\PP(E)}$. Let $X$ be a  
vector field on $B$ and denote $\widetilde
X$ the lift to $H^\nabla$ and $X'$ the lift to $H_{\PP(E)}$. We
have $\pi_* \Psi_* \widetilde X = X'$ and therefore 

$$\Omega_E (\widetilde X,   \widetilde J_B \widetilde X) =
\Psi^*\Omega (\widetilde X,   \widetilde J_B \widetilde X) =
\lambda r^2 \omega(X',   J'_B X') ,$$
or 
$$\Omega_E\big|_{(b,  v)}(\widetilde X,   \widetilde J_B \widetilde
X) = -\frac 1 {2i\pi}   {^t\bar v} \cdot F^\nabla (X,  J_B X) \cdot v .
$$
The last expression is positive for $v\neq 0$. 

\medskip 

\noindent {\bf Remark.} The case of projectively flat negative vector
bundles (with curvature equal to $-i\beta \tx{Id}$) corresponds precisely to a
connection form $\Omega_E$ which depends only on $|v|$ (and, of course,
on $b$). 

\medskip


\subsection{Convex fibrations} We
explain now a variation on an example from~\cite[~\S2.10]{CGMS02}. Let $G$ be
a compact Lie group with Lie algebra $\got g$. Let $X\stackrel \pi
\longrightarrow B$ be a principal
$G$-bundle with connection $\theta_A\in \Omega^1(X,  \got
g)$. We denote by $F_A\in \Omega^2(B,  X\times
_{\tx{ad}} \got g)$ its curvature and by $\tx{Hor}_A$ its horizontal
distribution. We assume that $B$ is symplectic
with symplectic form $\beta$. 

Let $(F,  \omega_F)$ be a symplectic manifold endowed with a
Hamiltonian action of $G$ with moment map 
$\phi_F: F\longrightarrow \got g^*$. We impose the following
conditions. 

\begin{itemize}
\item {\sc ($G$ - contact type boundary)} the boundary $\partial F$ is
  $G$-invariant and admits a $G$-invariant Liouville vector
  field $Z$ which is also conformal for the moment map:
  $$d\phi_F\cdot Z = \phi_F .$$
\item {\sc ($G$ - negativity)} there exists $J_B\in \mc J(B,  \beta)$ such
  that 
  $$\langle   F_A(X,  J_B X),  \phi_F(f)   \rangle \le 0$$
for all $X\in \tx{Hor}_A$ and $f$ in a neighbourhood of $\partial F$.
\end{itemize}

\medskip 

\noindent {\bf Remark.} The {\sc ($G$ - negativity)} condition is related to
Weinstein's notion of {\it fat} bundles i.e. $G$-principal bundles admitting a
connection $\theta_A$ such that the two form $\langle  
F_A(\cdot,\cdot),  \eta   \rangle$ is nondegenerate for all nonzero
$\eta\in \got g^*$. The {\sc ($G$ - negativity)} 
condition is also the crucial ingredient of the
construction in~\cite{CGMS02}.

\medskip 

We claim that the associated bundle 
$$X_F= X\times _G F$$
is a negative symplectic fibration. In order to see this we recall
Weinstein's construction of symplectic fibrations through symplectic
reduction as explained in~\cite[~\S2]{GLS96}. The connection
$\tx{Hor}_A$ defines the subbundle
$$M=\{\eta\in T^*X \ : \ \eta|_{\tx{Hor}_A} =0\} \subset T^*X .$$
One can show that $M\simeq X\times \got g ^*$. In any case, $M$
inherits from $T^*X$ a $2$-form $\omega_\Gamma$ which restricts to the
canonical 
symplectic form on the fibers $X_b\times \got g^* \simeq T^*G$. Then
$M\times F$ is a Hamiltonian $G$-space with moment map 
$$\phi(x,  \eta,  f) = \phi_F(f) + \eta .$$
The zero set $\phi^{-1}(0)$ is naturally identified with $X\times F$
and the symplectic reduction $\phi^{-1}(0)/G$ is isomorphic to
$X_F$. The (pre)symplectic form 
$\omega_\Gamma+\omega_F$ on $M\times F$ is $G$-invariant, and the same
is true for its restriction to $\phi^{-1}(0)$. Moreover, for any
$\xi\in \got g$ we have $\iota(\xi_M, 
\xi_F)(\omega_\Gamma+\omega_F)=-d\langle \phi_F+ \tx{pr}_2, 
\xi\rangle =0$ on $\phi^{-1}(0)$, hence
$(\omega_\Gamma+\omega_F)|_{\phi^{-1}(0)}$ is the pull-back of a
(pre)symplectic form on $M\times _G F=X_F$ which we denote by
$\omega_{\Gamma,  F}$. We have denoted by $\xi_M$, $\xi_F$ the
infinitesimal generators of the action of $G$ on $M$ and $F$
respectively. The form $\omega_{\Gamma,  F}$ restricts to the
symplectic form $\omega_F$ in the fibers of $X_F$. 
We define 
$$\Omega = \omega _{\Gamma,  F} .$$
The {\sc (negativity)} condition for $\omega_{\Gamma,  F}$ is now
equivalent to the {\sc ($G$ - negativity)} condition above because 
$\omega_{\Gamma,  F}$ acts at a point $[x,  f]\in X_F$ as $-\langle
  \phi_F(f),  F_A(x)(\cdot,\cdot)   \rangle$
(cf.~\cite{GLS96}). One can also prove that the horizontal
distribution of $\omega_{\Gamma,  F}$ is the distribution induced by
the connection $\theta_A$. On the other hand, parallel transport
$\tau_\gamma$ along a curve $\gamma$ in $B$ with respect to 
the latter horizontal distribution acts as $\tau_\gamma([x, 
f])=[\tau_\gamma(x),  f]$. Because $\partial F$ is invariant under
$G$ we infer that parallel transport preserves $\partial X_F=X\times
_G \partial F$ and the  {\sc (monodromy)} condition is 
satisfied. The {\sc (symplectic fibration)} condition is satisfied by
construction of $X_F$ and we are left to verify the {\sc (contact type
  boundary)} condition. The natural Liouville vector field on $X\times
F$ is 
$$\bar Z_{(x,  \eta,  f)} = (0,  \beta,  Z_f) .$$
By $G$-invariance $\bar Z$ 
descends to a Liouville vector field on $X_F$ provided it is
tangent to $\phi^{-1}(0)$, and this is equivalent to
$\phi_F(f)=d\phi_F(f)\cdot Z_f$ in a neighbourhood of $\partial
F$. Moreover, if the last condition holds then verticality in $X_F$ is
automatic.  

\medskip 

\noindent {\bf Remark.} Negative vector bundles, seen as associated
bundles of the corresponding frame bundles, are a special instance
of the above construction.


\section{Fibered symplectic cohomology groups} \label{sec:construction of FH}

We define now the Floer or Symplectic 
cohomology groups $FH^*(E)$ for negative symplectic
fibrations.  The key concept is that of an {\it asymptotically linear
Hamiltonian}, and from this point of view our definition can also be
thought of as a bridge between the one of Viterbo~\cite{Viterbo99}, 
who uses Hamiltonians that are linear at infinity, and
the one of Floer and Hofer~\cite{FH94}, who use Hamiltonians
that are asymptotically quadratic on $\CC^n$ (with the somewhat
surprising remark that ``quadratic'' is the same as ``linear'' after
the change of variables $S=r^2$). 

The main feature of the Floer cohomology groups that we define in this
paper is that they have the same functorial properties as those
of Viterbo in~\cite{Viterbo99} (see~\S\ref{sec:properties of
  Floer homology}). 

\medskip 

 \noindent {\bf Convention.} We shall assume in this section that the form
      $\omega_\epsilon$ is symplectic for $0<\epsilon \le 1$. For
       clarity we drop the subscript $\epsilon$ and work with
       $\omega=\omega_1$. 


\subsection[$C^0$-bounds for asymptotically linear
Hamiltonians]{\hspace{-.1cm}Admissible Hamiltonians and almost complex
  structures. 
  $C^0$-bounds}  \label{sec:proof of C0 bounds}

The crucial ingredient of the construction is the proof of a priori
$C^0$-bounds for solutions $u:\RR \times \Ss^1 
\longrightarrow \widehat{E}$ of the equation 
\begin{eqnarray}  \label{eq Floer parametre} 
u_s + \widehat{J}(s,   t,   u) (u_t - X_H(s,   t ,   u)) =0  
 , \\ 
-\infty < \inf_{s \in \RR} A_{H(s)}(u(s)), \ \sup_{s \in \RR} 
A_{H(s)}(u(s)) < + \infty .  \label{energie bornee} 
\end{eqnarray}  
 
Here  $H(s,   t ,   u), \  \widehat{J}(s,   t,   u)$ is a homotopy
of Hamiltonians and almost complex structures on which we impose 
additional constraints as described below.
The constraints on $H$ and $\widehat J$, 
as well as the proofs of the $C^0$-estimates, are adapted
from the papers of Cieliebak, Floer and Hofer~\cite{CFH95,FH94}. 

A point $u \in \widehat{E}$ which belongs
  to $\partial E \times [1,   \infty[$ will be denoted $u=(\bar{u},  
  S)$.

\begin{defi} \label{defi:standard almost complex structure on
    symplectic cone}
 
 Let $(\Sigma,   \lambda)$ be a contact manifold. The manifold 
 $$\big(\Sigma \times ]0,   \infty[,  
  d(S\lambda) \big), \ S\in ]0,   \infty[$$ 
  is called \emph{the symplectic cone} over $\Sigma$.
  The \emph{Reeb vector field}
  $X_{\tx{Reeb}}$ on $\Sigma$ is defined by $\iota_{X_{\tx{Reeb}}}
 d\lambda \equiv 0$, $\lambda(X_{\tx{Reeb}}) = 1$. The \emph{contact
 distribution} 
 $\ker \,  \lambda$ is denoted by $\xi$.  
  An almost complex structure $J$ on $\Sigma \times ]0,   \infty[$ is
  called {\rm standard} if 
\begin{eqnarray}  J_{(\bar u,   S)}(\DP{}{S}) & = &   
      \frac 1 {CS} X_{\tx{Reeb}}(\bar u)  , \nonumber \\  
    J_{(\bar u,   S)}(X_{\tx{Reeb}}(\bar u) ) & = &  -CS \DP{}{S}  ,
    \label{def:standard almost complex structure}\\ 
    J_{(\bar u,   S)}|_\xi & = & J_0  , \nonumber
\end{eqnarray}
where $J_0$ is an almost complex structure compatible with $d\lambda$
on $\xi$ and $C>0$ is a positive constant.    
\end{defi}

Standard almost complex structures are precisely the ones that
are preserved by homotheties in the $S$ variable. As an example, for
$\Sigma = \Ss^{2n-1}(1) \subset \CC^n$ and $\lambda = \frac 1 2 \sum
x_i dy_i - y_i dx_i$, the manifold $\big(\Sigma \times ]0,   \infty[,
  d(S\lambda|) \big)$ is
    symplectomorphic to $\big(\CC^n \setminus \{0\},   \sum
    dx_i \wedge dy_i \big)$ through the
    map $(\bar u,   S) \longmapsto \varphi^{\ln S}_X(\bar u)$. Here
    $\varphi ^t_X$ stands for the flow of $X(z) = \frac 1 2 z$. 
    The inverse map satisfies  
    $S(z) = \sqrt{|z|}$ and the canonical complex structure on
    $\CC^n \setminus \{0\}$ translates into a standard
        almost complex structure on $\Sigma \times ]0,   \infty[$
        which satisfies (\ref{def:standard almost complex structure}) 
            with $C=4$.

The metric $d(S\lambda)(\cdot,   J\cdot)$ 
associated to a standard almost complex structure will be
called {\it conical}. The following 
homogeneity property is straightforward:

\begin{equation}\label{homogeneite}   
   \big| v+a\DP{}{S} \big| ^2 _{(\bar u,S)}= S \big| v+\frac{a}{S}\DP{}{S}  
   \big| ^2 _{(\bar u,1)},  
   \qquad v\in T_{\bar u} \Sigma, \ a \in \RR  . 
  \end{equation}

\begin{defi} \label{defi:split almost complex structures}
  Let $F\hookrightarrow E \stackrel \pi \longrightarrow B$ be a
  fibration satisfying the assumptions \textsc{(symplectic fibration)} and
  \textsc{(contact type boundary)}. Let $H$ be the horizontal
  distribution on $\widehat E$. An almost complex
  structure $J$ on $\partial E \times [1,   \infty[$ is called 
  \emph{(standard) split} if 
  $$J = J_V \oplus \widetilde J_B  ,$$
  where $J_V$ is a (standard) almost complex structure in the fibers and
  $\widetilde J_B$ is the lift to $H$ of an almost complex structure
  $J_B$ on $B$ which is $\beta$-tame. 
\end{defi} 

\medskip 

The thrust of the present section is that the a priori $C^0$-bounds on
Floer trajectories (and ultimately a variant of the maximum principle) hold in
$\widehat E$ with respect to  
almost complex structures that are standard split at infinity.  

\begin{defi} \label{defi:admissible deformations a c structures}

An \emph{admissible homotopy} of almost complex structures on
$\widehat E$ is a smooth
family $\widehat{J}(s,   t)$, $s\in \RR$,
$t\in \Ss^1$   
of almost complex structures tamed by  $\omega = \pi^*\beta + \Omega$ 
such that the following conditions hold.
 
\renewcommand{\theenumi}{\roman{enumi}}

\begin{enumerate} 
  \item $\widehat J$ is standard split for
    $S$ large enough, i.e. there exists $R\ge 1$ such that 
   \begin{equation} \label{standard S grand} 
      \widehat{J}(s,   t,   \bar{u},   S) =  
J_V(s,   t,   \bar{u},   S) \oplus \widetilde J_B(s,   t,   \bar
   u), \quad S \ge R .
   \end{equation} 
  \item $\widehat J$ is constant for $|s|$ large enough, i.e. there
  exists $s_0 >0$ such that   
    \begin{equation} \label{homotopie J}
     \begin{array}{rcl}  
   \widehat{J}(s,  t,   u) & = & J_-(t,   u), \qquad s \le -s_0  ,
    \\  
   \widehat{J}(s,  t,   u) & = & J_+(t,   u), \qquad s  
     \ge s_0  .
\end{array} 
    \end{equation} 
\end{enumerate} 

Here $J_V$ is a standard almost complex structure in the fiber and
$\widetilde J_B$ is the horizontal lift of an almost complex structure
$J_B$ on $B$ which is $\beta$-tame.   
 
\end{defi} 

\medskip

 We now define admissible homotopies of Hamiltonians. In the usual
 setting of Floer homology, these are functions  
 $H(s,   t,   \bar  u,   S)$, $s\in \RR$, $t \in \Ss^1$ 
 with a special asymptotic behaviour that ensures compactness for the
 moduli spaces of finite energy solutions of Floer's equation 
 $u_s = -\nabla  A_{H(s,\cdot)}(u(s, \cdot))$. Here 
$$
A_{H(s,\cdot)}(x) = -\int_{D^2}\bar x^*\omega - \int_{\Ss^1}H(s)
\circ x
$$
is the symplectic action defined on the space
$\Lambda_0 \widehat E$ of $1$-periodic contractible loops
in $\widehat E$, $\bar x$ denotes an extension of the loop $x$ over
a disc and the $L^2$-gradient of $A_H$ is 
$$
\nabla A_H(x) = \widehat J   \dot x - \nabla  H (x)  .
$$   
The construction of the spectral sequence will crucially require
the use of negative \emph{pseudo-gradient} trajectories for the
action functional. We shall
be interested in solutions $u(s,   t)$ of  
\begin{equation} \label{eq:pseudo grad Floer} 
u_s = - \mc Y(u(s, \cdot))  ,
\end{equation} 
\begin{equation} \label{eq:energie bornee} 
-\infty < \inf_{s \in \RR} A_{H(s)}(u(s)), \ \sup_{s \in \RR} 
A_{H(s)}(u(s)) < + \infty ,  
\end{equation} 
where $\mc Y$ is a pseudo-gradient for some action functional 
$A_H$ on $\Lambda_0 \widehat E$, i.e. 
$$dA_H(x) \cdot  \mc Y(x) \ge 0  , $$
with equality if and only if the loop $x$ is a critical point of
$A_H$ (hence a periodic orbit of $X_H$). We shall actually need the
stronger pseudo-gradient condition 
$$ dA_H(x) \cdot  \mc Y(x) \ge a ^2 \parallel \mc Y(x) \parallel
^2, \quad a > 0, $$
with $\mc Y$ and $\nabla A_H$ having the same zeroes. We shall use
vector fields  $\mc Y$ of the type 
\begin{equation} \label{eq:pseudo grad is CR} 
\mc Y(x) = \widehat J   \dot x - Y(x)  ,
\end{equation} 
where $Y$ is a vector field on $\widehat E$. The vector field $Y$ will
therefore be, along with the Hamiltonian $H$, part of the data
defining an \emph{admissible deformation}. In the following we let 
$X=X^{\tx{v}} + X^{\tx{h}}$ be the decomposition of a vector $X\in
T\widehat E$ in its vertical and horizontal parts.

\begin{defi} \label{defi:admissible Hamiltonian deformation}
 
Let $\widehat J$ be an admissible homotopy of almost complex
structures. An \emph{admissible pseudo-gradient deformation} consists of a one
parameter family $H(s,  t,   u)$, $s\in \RR$, $t\in \Ss^1$ of Hamiltonians 
and of a one parameter family of vector fields $Y(s, 
t,   u)$ on $\widehat E$, which satisfy the following properties. 
 
\renewcommand{\theenumi}{\roman{enumi}}
\begin{enumerate} 
 
 \item (strong pseudo-gradient) Let 
   $$\mc Y(s,   x) = \widehat J(s)   \dot x - Y(s)\circ x .$$
   We require the existence of a
   function $a:\RR \longrightarrow [0,   \infty[$ 
   with nowhere dense vanishing locus such that, for every loop $x \in
   \Lambda_0\widehat E$, we have 
  \begin{equation}\label{Y is pseudo grad} 
   dA_{H(s)}(x) \cdot \mc Y(s,   x) \ge a(s)^2 \parallel \mc Y
   (s,   x) \parallel ^2 _{_{\widehat J(s)}}.  
  \end{equation}
  Moreover, $\mc Y(s,   x)$ is required to have the same critical
  points as $\nabla A_{H(s)}$. 
 
\item (monotonicity) $H$ is increasing  
   \begin{equation} \label{croissance de H}
     \DP{H}{s}(s,  t,   u) \ge 0   
   \end{equation} 
   and there exists $s_0>0$ such that 
   \begin{eqnarray} \label{homotopie H} & & \\
   H(s,  t,   u) & = & H_-(t,   u), \quad Y(s,   t,   u) \ = \ Y_-(t,
     u), \quad a(s)=a_->0, \quad s \le -s_0  , \nonumber \\  
   H(s,  t,   u) & = & H_+(t,   u), \quad  Y(s,   t,   u) \ = \ Y_+(t,
     u),\quad a(s) = a_+>0, \quad s \ge s_0  .   \nonumber 
   \end{eqnarray} 
 
\item (asymptotes) There exists $f(s,   \bar u): \RR   \times   \partial
  E  \longrightarrow ]1-\delta,   \infty[$  so that the function $F(s,   
\bar{u},   S) = S f(s,   \bar u) $ satisfies
   \begin{eqnarray} \label{homotopie asymptotiquement lineaire}  
     |Y^{\tx{v}}- (\nabla{F})^{\tx{v}}|_{J_V}
     / \sqrt{S}  
\longrightarrow 0, & &  \\ 
  \label{cond asy sur H}
     |\nabla H- \nabla{F}|_{J_V\oplus\widetilde J_B} / \sqrt{S} 
\longrightarrow 0, & &  \\
  \label{second derivatives s and S}
     \left| \DP{^2H}{s\partial S} - \DP{^2F}{s\partial S} \right|
     \longrightarrow 0, & & \quad S \longrightarrow \infty,     
   \end{eqnarray} 
  uniformly in $s$, $t$ and $\bar u$. Moreover, the function $f$
  is required to satisfy the following conditions: 
   \begin{itemize} 
     \item For every $s\in \RR$ and every large constant $c$, 
       the horizontal distribution $H$ is tangent to 
       \begin{equation} \label{eq:cond on graph of f}
         \big( \tx{graph} \  c/f(s) \big) \subset  \partial E \times
       ]1-\delta,   \infty[ .
       \end{equation} 
     \item There exists $s_0>0$ such that 
    \begin{eqnarray} \label{homotopie f} 
       f(s,   \bar u) & = & f^-(\bar u), \qquad s \le -s_0  , \\
       f(s,   \bar u) & = & f^+(\bar u), \qquad s \ge s_0  \
       . \nonumber  
    \end{eqnarray} 
     \item The vector fields $X_{F^\pm}$ have no $1$-periodic orbits
       at infinity, where  
     \begin{equation} \label{pas d_orbites periodiques pour F et G}  
       F^\pm (\bar{u},   S) = S f^\pm (\bar u).
     \end{equation} 
   \item \begin{equation} \label{croissance de f}
       \partial_s f \ge 0 . 
     \end{equation} 
   \item If the $1$-periodic orbits of $X_{F(\widehat s)}$ are not
   contained in a compact set, then 
    \begin{equation} \label{croissance stricte si orbite periodique}   
    \partial_s f|_{s=\widehat s} \ge \epsilon(\widehat s) >0 .  
    \end{equation}  
   \end{itemize} 

\item (boundedness) There exists a constant $c>0$ such that 
  \begin{equation} 
    \label{eq:horizontal bound on Y}
    |Y^{\tx{h}}|_{\beta,   \widetilde J_B} \le c  , 
  \end{equation}
   \begin{equation} \label{controle dt}
     | \partial_t  Y (s,   t,   \bar{u},   S) |_{J_V\oplus
         \widetilde J_B} \le c(1 + \sqrt{S})  , 
   \end{equation} 
   \begin{equation} \label{controle derivee seconde}  
| \nabla _X Y(s,   t,   \bar{u},   S)|_{J_V \oplus \widetilde J_B} \le 
c|X|_{J_V\oplus \widetilde J_B}, 
\quad X \in T_{(\bar{u},   S)}\big( \partial E \times [1,  \infty[
    \big) .  
   \end{equation} 

\end{enumerate} 
\end{defi} 

\medskip 

 \noindent {\bf Remarks.} 

1. Condition~(\ref{eq:cond on graph of f}) is equivalent to saying
that the graph of $f(s)$ restricted to some fiber is a monodromy
invariant hypersurface, and determines the values of $f(s)$ on the
whole of $\partial E$ through parallel transport along $H$. This
condition is clearly void if the base $B$ is a point. In case the base
$B$ is not a point it ensures that $X_{F(s)}$, $s\in \RR$ is preserved
by parallel transport along $H$, being colinear with the Reeb vector field
on the level sets $c/f(s)$, $c >> 1$. We shall crucially use this
fact in the proof of Proposition~\ref{borne H 1 2}. 

2. Condition~(\ref{pas d_orbites periodiques pour F et G}) means that,
for any choice of a large enough constant $c$, the graphs of $c/f^\pm$ 
have no closed characteristics of period $1$. This
can be achieved for example if $f^\pm$ are equal to constants not
belonging to the period spectrum of $\partial E$ (this is a generic
condition). 
 
3. The pseudo-gradient condition~(\ref{Y is pseudo grad})
obviously holds with $a(s) \equiv 1$ in case $\mc Y(s,   x) =
\nabla^{\widehat J(s)} A_{H(s)}(x)$, which corresponds to the usual Floer
equation. 


4. The function $F$ defined above satisfies
$$|X_F|_\omega = O(\sqrt S) .$$

5. The above conditions are satisfied in the $s$-independent case 
   by vector fields $Y = \nabla h + \widetilde {\nabla f}$, where
   $h=h(S)$ is linear for $S$ big enough and $f: B \longrightarrow
   \RR$ is a smooth function. We must take in this case $F= h$ and $H
   = h+ \widetilde f$. The pseudo-gradient property is the only 
   nontrivial one, and we refer to~\S\ref{sec:main section
   pseudo gr vector fields} for a proof. There are two other nonempty
   properties, namely~\eqref{eq:horizontal
   bound on Y} which holds since $Y^{\tx{h}}=\widetilde J_B
   \widetilde X_f$, and~\eqref{cond asy sur H} which holds because
   $|X_f|_{\widetilde J_B}$ is bounded. This in turn is implied by the
   fact that $\Omega|_{H_0}=Sd\Theta|_{H_0}$, hence the component of
   $X_f$ on some nondegeneracy subspace of $\Omega|_{H_0}$ goes to
   zero as $S \to\infty$, whereas the component on the degeneracy
   subspace stays bounded.

\medskip  

The proof of the $C^0$-bounds follows the arguments in~\cite{FH94}. We
recall the notations and a crucial technical result therein. 

\medskip 
 
Let $\alpha: \RR \times \Ss^1 \longrightarrow \RR$ be a smooth
function. Let $\delta >0$. We denote by  $\Gamma_\delta$ the set of
all sequences  $(s_k)_{k \in \ZZ}$ such that 
\begin{equation} \label{Gamma delta}  
0 < s_{k+1} - s_k \le \delta, \quad  k \in \ZZ ,
\end{equation}  
$$ s_k \longrightarrow \pm \infty, \ k \longrightarrow \pm \infty .$$ 
For $s=(s_k) \in \Gamma_\delta$ we define
$$[\alpha ] ^s \ = \ \sup \{ \alpha(s_k,   t): \ k\in \ZZ, \ t \in \Ss^1  
\} .$$ 
Let 
$$[\alpha ] _\delta  \  = \ \inf\{ [\alpha]^s: \ s \in \Gamma_\delta 
\} .$$ 
 
\begin{prop}[\cite{FH94}, Prop. 8] \label{inequation elliptique}  
Let $A, \ B, \ \lambda \ge 0$  be nonnegative real numbers and let  
$\delta >0$ satisfy 
$$\delta^2 \lambda < \pi ^2 .$$ 
There exists a positive constant $C= C(A,   B,   \lambda,   \delta) >0$ 
such that, for any function $\alpha : \RR \times \Ss^1 \longrightarrow [0,
    \infty[$ satisfying 
\begin{eqnarray*}
-\Delta \alpha - \lambda \alpha \le A & \tx{ on } \RR \times 
\Ss^1 , \\ 
 \lbrack \alpha \rbrack  _\delta < B ,  & 
\end{eqnarray*} 
we have 
$$\sup \{ \alpha(s,   t) : \ (s,   t) \in \RR \times \Ss^1 \} \le 
C .$$ 
\end{prop} 
 
This is a key result involving the maximum principle, which is
unavoidable in the proof of the $C^0$-bounds for all versions
of Floer homology defined on open
manifolds. Proposition~\ref{inequation elliptique} will be applied to
the function $\alpha = S\circ 
u$, with $u$ an arbitrary solution of 
(\ref{eq:pseudo grad Floer}  - \ref{eq:energie 
  bornee}), and will ultimately yield the following result. 
 
\medskip 

\begin{thm}[a priori $C^0$-bound; compare \cite{FH94}, Thm. 12] 
\label{estimation C 0} 
Let $\widehat{J}$, $H$ and $Y$ satisfy conditions  
(\ref{standard S grand} -  \ref{homotopie J}) and 
(\ref{Y is pseudo grad} - \ref{controle derivee seconde}). 
There is a constant $d = d(\widehat{J},   H,   Y) > 0$ such that
any solution of 
\begin{eqnarray*}   
u_s + \widehat{J}(s,   t,   u) u_t - Y(s,   t ,   u) =0 
 , \\ 
-\infty < \inf_{s \in \RR} A_{H(s)}(u(s)), \ \sup_{s \in \RR} 
A_{H(s)}(u(s)) < + \infty   
\end{eqnarray*}  
 satisfies
\begin{equation} \sup_{(s,   t) \in \RR \times \Ss^1} S \circ u (s, 
  t) \ \le \ d . 
\end{equation}  
\end{thm}  
 
\medskip 

The rest of this section is devoted to the proof of
Theorem~\ref{estimation C 0}. 

\medskip 

\noindent {\bf Notation.} 

\medskip 

1. a. We {\it extend the function $\sqrt S$}, 
which is canonically defined on $\partial
E \times [1-\delta,   \infty[$, to a smooth function on the whole of
    $\widehat E$ as follows. Consider a strictly increasing 
    smooth bijective map $\rho :
    [1-\delta,   1] \longrightarrow [0,   1]$ whose derivatives
    vanish at infinite order at $0$ and $1$. 
Our smooth extension of $\sqrt S$ is 
defined to be equal to $\rho(S) \sqrt S$ on $\partial E \times [1-\delta,  
  \infty[$, and identically equal to zero 
on $\widehat E \setminus \partial E \times [1-\delta,   \infty[$. We
still denote this smooth extension by $\sqrt S$, and its square is 
a well defined function on $\widehat E$ which we denote by $S$. 

\medskip 

1. b. Every point $u\in \partial E   \times   [1-\delta,   \infty[$ 
can be uniquely written as 
$u=(\bar u,   S)$, $\bar u \in \partial E   \times   \{
    1-\delta \}$. The map {\it $u \longmapsto \bar u$ continuously 
    extends} as 
    the identity over $E \setminus \partial E \times [1-\delta,   1]$
    and we  denote the extension again by $u
    \longmapsto \bar u$. 

\medskip 

1. c. Every point $u\in \widehat E$ is now uniquely characterized by
the pair $(\bar u,   S)$ and we shall identify the two in the sequel.

\medskip 

2. We define the following three norms. 
   \begin{eqnarray*}
   | X |_\beta ^2  & = & \Omega(X^{\tx{v}},   J_V X^{\tx{v}}) +
   \pi^*\beta (X^{\tx{h}},   
   \widetilde J_B X^{\tx{h}})  , \\ 
   | X |_\omega ^2  & = & \Omega(X^{\tx{v}},   J_V X^{\tx{v}}) +
   \omega (X^{\tx{h}},   \widetilde J_B X^{\tx{h}})  , \\  
   | X |_\Omega^2  & = & \Omega(X^{\tx{v}},   J_V X^{\tx{v}}) +
   \Omega (X^{\tx{h}},   
   \widetilde J_B X^{\tx{h}}) . 
   \end{eqnarray*}

  The last expression only defines $|\cdot |_\Omega$ as a semi-norm on
  $H$. Its utility will nevertheless become appearant in the sequel. 
  The above norms all satisfy the homogeneity property
  (\ref{homogeneite}) along the vertical
  distribution, whereas on the horizontal
  distribution, which is preserved by the Liouville flow, they satisfy
  the inequality    
  \begin{equation} \label{horizontal homogeneity}
     |X^{\tx{h}}|_{(\bar u,   S)}^2 \le S |X^{\tx{h}}|_{(\bar u,  
      1)}^2, \ S\ge 1 
      .
  \end{equation} 

\medskip

\noindent {\bf Convention.}  {\it The default norm used in the sequel is
   $|\cdot|_\omega$.} 

\medskip 

3. We define 
\begin{eqnarray*} 
L^2(\Ss^1,   \widehat E)  & = & 
\{ x : \Ss^1 \longrightarrow \widehat E
    \tx{ measurable} \ : \ \sqrt{S\circ x} \in  L^2(\Ss^1,   \RR)
    \}, \\ 
H^1(\Ss^1,   \widehat{E}) & = & \{ x\in L^2(\Ss^1,   
\widehat{E}) \ : \ \dot{\bar x} \in 
L^2(x^*T E), \quad \sqrt{S\circ x} ^{ \ \prime} \in L^2(\Ss^1,   \RR)
\} .
\end{eqnarray*}  
Here $\widehat{E}$ is endowed with 
the metric $\langle X,   Y \rangle _\omega = \frac 1 2 \big(\omega (
X,   JY) + \omega (Y,   JX)\big)$, $J=J_V\oplus \widetilde {J_B}$ 
and with its associated Lebesgue measure. The meaning of the
condition $\dot {\bar x} \in L^2(x^* TE)$ is the following:  $\dot
{\bar x}$ is well defined as a distribution once we choose
an embedding of $E$ into some Euclidean space. We require it to be
an $L^2$ function and this does not depend on the choice of the 
embedding.

\begin{lem} \label{lem:compact embedding} 
  There is a compact embedding 
  $$H^1(\Ss^1,   \widehat E) \hookrightarrow C^0(\Ss^1,  
  \widehat E) .$$
\end{lem}

\demo Let $(y_k)$ be a sequence in $H^1(\Ss^1,   \widehat E)$. We
denote $S_k = S\circ y_k$. The embedding $H^1(\Ss^1,  
\RR)\hookrightarrow C^0(\Ss^1,   \RR)$ is compact and therefore 
a subsequence of $\sqrt{S_k}$ converges
uniformly to a continuous function $\sqrt{S_0}$. 
As a consequence, the 
corresponding subsequence of $y_k$ takes values in a compact 
set $\widehat E_c = \{ S\le c \}$. 
The embedding $H^1(\Ss^1,  
\widehat E_c)\hookrightarrow C^0(\Ss^1,   \widehat 
E_c)$ is again compact 
 and we get a subsequence converging
uniformly to a continuous limit $y$, with $S\circ y = S_0$. 
\hfill{$\square$}


%
%

\begin{prop}[compare~\cite{FH94}, Lemma 10] \label{borne H 1 2}  
Let $\widehat{J}$, $H$ and $Y$ satisfy  conditions  
(\ref{standard S grand} -  \ref{homotopie J}) and 
(\ref{Y is pseudo grad} - \ref{controle derivee seconde}).
Let $s\in \RR$ be such that $a(s)>0$. 
For any choice of $c > 0$ there is a constant $d=d(c,   s) > 0$ 
such that, if we have 
\begin{equation}  \label{eq:temporary label} 
 dA_{H(s)}(x)\cdot \mc Y(s,   x) + 
\int_{0}^1\DP{H}{s}(s,   t,   x(t))   dt \ \le \ c 
\end{equation}  
for some $x \in H^1(\Ss^1,   \widehat{E})$, then  
\begin{equation}  
||x||_{H^1} \le d . 
\end{equation}
We have denoted $||x||^2_{H^1} = || \sqrt{S\circ x}||^2_{H^1} + 
|| \dot{\bar x}||^2_{L^2}$ for some fixed embeding of $E$ in a
Euclidean space. 
\end{prop}

\noindent {\bf Remark.} The statement is false if $a(s)$ is allowed to
vanish. As a counterexample 
one can consider a weak pseudo-gradient vector field of the form 
$Y=J_VX_h+ \widetilde J_B \widetilde{X_f}$ as in \S\ref{sec:main 
section pseudo gr vector fields}, where $h=h(S)$ has critical slope at
infinity. The expression
(\ref{eq:temporary label}) is bounded if $x$ is a periodic orbit, but
$x$ can nevertheless go to infinity.

Note also that, if
$a(s)\ge a_0 >0$ and $c\le c_0<\infty$, then $d(c,s)\le d(a_0,c_0)$. 

\medskip

\demo  
The strong pseudo-gradient condition (\ref{Y is pseudo grad})
implies that $\parallel \mc Y(s,   x) \parallel ^2 _{_{\widehat J(s)}}$
is bounded by $c'=c/a(s)^2$. 
The asymptotic behaviour of $\widehat J$ ensures that the
norms $|\cdot|_{\beta,   \omega,   \Omega}$ 
defined with respect to $\widehat J(s,   t)$, $s\in
\RR$, $t\in \Ss^1$ are equivalent to the corresponding norms defined
with respect to some fixed almost complex structure $J$ which is
standard split  at infinity. Moreover, the operator norm of $\widehat
J$ is bounded. We can therefore assume in the sequel,
without loss of generality, that $\widehat
J = J$. 

\medskip 

{\it Claim~1.} If $\parallel \mc Y(x) \parallel_{L^2}$ and $\parallel \sqrt
{S \circ x} \parallel_{L_2}$ are bounded, then $\parallel x \parallel
_{H^1}$ is bounded. 

{\it Proof of Claim~1.} We remind that $\mc Y(x) = J\dot x - Y(x)$ and
we have $|Y|_\omega \le \bar c_1 (1+\sqrt S)$. The hypothesis implies
therefore that $\parallel \dot x \parallel _{L^2}$ is bounded. 

Let $S(t) = S\circ x(t)$. At a point $t$ where $S(t)\ge 1$ we have
$x(t) = (\bar x (t),   S(t))$ and $|\dot x (t) | ^2 = |\dot{\bar x}(t)|_{x(t)}^2
+|S'(t)\DP{}{S}|^2 = |\dot{\bar x}(t)|_{x(t)}^2 + S'(t)^2/S(t)$. This ensures  
  $ { \sqrt{S}^{ \ \prime} } ^{ \ 2 } \le \frac 1 4 | \dot{x}(t) |^2$. 
At a point $t$ where $S(t)\le 1$ we have   
 ${ \sqrt{S}^{ \ \prime} } ^{ \ 2 } = \big[ d(\sqrt{S})\cdot \dot{x}(t) 
 \big] ^2 \le c |\dot{x}(t) |^2 $, the norm of $d\sqrt{S}$ being bounded
 on $E$. 

On the other hand we clearly have 
 $|\dot{\bar x}(t)|_{\dot{\bar x}(t)} \le |\dot{\bar
    x}(t)|_{x(t)} \le |\dot x (t) |$ at a point $t$ where $S(t) >0$,
  while for $S(t)=0$ we have $\dot{\bar x}(t) = \dot x(t)$ hence
  $|\dot{\bar x}(t)| = |\dot x(t)|$. 

This shows that $\parallel \dot x \parallel_{L^2}$ bounds $\parallel
\sqrt{S}^{ \ \prime} \parallel _{L^2}$ and $\parallel \dot{\bar x}
\parallel_{L^2} $. Claim~1 is proved.

\medskip 

 We are left to prove that $\parallel \sqrt
 {S\circ x} \parallel _{L^2}$ is bounded. Arguing by contradiction, let
 us suppose the existence of a sequence $(s_k,   x_k)$ such that   
 \begin{equation} \label{eq:hypothesis for contradiction} 
 \parallel \sqrt{S\circ x_k} \parallel \longrightarrow \infty  ,
 \end{equation}
 $$\parallel \mc Y(s_k,   x_k) \parallel^2 \le c', \qquad \int_0^1 
 \DP{H}{s} (s_k,   t,   x_k(t))   dt \le c .$$ 
 By (\ref{homotopie H}) we can suppose $s_k \longrightarrow
 \widehat{s}$. Let $S_k = S\circ x_k$ and   
 $\lambda_k = \parallel \sqrt{S_k} \parallel / \Lambda$, with $\Lambda 
 > 1$ a constant to be chosen below. 
 
\medskip 
 
{\it Claim~2.} The sequence $v_k(t) = (\bar{x}_k(t),  
S_k(t)/\lambda_k^2)$ has a uniformly convergent subsequence. 

{\it Proof of Claim~2.} Let $\widetilde S_k = S\circ v_k$. 
We have $\parallel \sqrt {\widetilde
  S_k} \parallel = \Lambda < \infty$ and, by Lemma~\ref{lem:compact
  embedding}, it is enough to prove that $\parallel \dot v_k
\parallel$ is bounded. We shall drop the subscript $k$ in the next
paragraph. For $S>0$ we denote by $S^o$ (resp. $\widetilde S^o$) the
``true'' $S$-coordinate corresponding to $S$ (resp. $\widetilde
S$), with values in $]1-\delta,   \infty[$ and with respect to which the
homogeneity property~(\ref{homogeneite}) is verified. 

We first prove the following inequality, with the
convention $\widetilde S^o/S^o=1$ if $S=0$:  
\begin{equation}
  \label{eq:rescaling}
  |\dot v(t)|^2 \le \bar c \frac {\widetilde S^o}{S^o} |\dot x|^2_\Omega +
  |\dot{\bar x} ^{\tx{h}} |^2_\beta  ,
\end{equation}
where $\bar c \ge 1$ is some constant. 
At a point $t$ where $S(t)=0$ we have $\dot v = \dot
x$ and (\ref{eq:rescaling}) is clear. Let us examine the situation at
a point $t$ where $S(t) > 0$. After expansion (\ref{eq:rescaling})
reduces to
\begin{equation*}
  \frac {(\widetilde S^{o   \prime})^2} {\widetilde S^o} \le
  \bar c \frac{\widetilde S^o}{S^o} \frac{(S^{o   \prime})^2}{S^o},
\end{equation*}
which is equivalent to 
\begin{equation}
  \label{eq:rescaling variant}
  \big| \frac{\widetilde S^{o   \prime}}{\widetilde S^o} \big| \le
  \bar c_1
  \big| \frac{S^{o   \prime}}{S^o} \big| .
\end{equation}
 We prove (\ref{eq:rescaling variant}) under the
 assumption that $\widetilde S^{o   \prime}$ and $S^{o   \prime}$
 are both strictly positive (note that they necessarily have the same
 sign). The proof applies as such to the negative case as well. 

 We have $S^o=g(S)$ where $g:]0,   \infty[ \longrightarrow ]1-\delta,
   \infty[$ is a strictly increasing diffeomorphism, with $g(y)= y$
 for $y\ge 1$ and $g'(0^+)=\infty$. The inverse $f=g^{-1}$ has to
 vanish at infinite order at $1-\delta$ and has as typical profile
 $f(x)=\exp \big(-1/(x-(1-\delta)) \big)$. Therefore $g$ can be chosen
 to be equal to $g(y)= (1-\delta) - 1/\ln y$ near $0$, say on $]0,  
 \epsilon]$, with $g'(y)=1/y(\ln y)^2$, 
 whereas on $[\epsilon,   1]$ we have $0<c \le g' \le C$.

 a. Assume $0<S\le1$. This means that $\widetilde S^o$ and $S^o$ belong
 to $]1-\delta,   1]$. It is enough to prove that, for any
 $\lambda\ge 1$, we have  $g'(S/\lambda)\le
 \bar c_2 \lambda g'(S)$ for some constant $\bar c_2>0$. Now: 
 \begin{itemize}
 \item If $S\in ]0,   \epsilon]$ we have $g'(S)=1/S(\ln S)^2$ and
   $g'(S/\lambda) 
\le \lambda g'(S)$. 
 \item If $S\in ]\epsilon,   \epsilon \lambda]$ we have $g'(S/\lambda)\le
\epsilon\lambda/S \cdot g'(\epsilon) \le \lambda C \le \lambda C/c
\cdot g'(S)$. 
 \item If $S\in ]\epsilon \lambda,   1]$ we have $g'(S/\lambda) \le
   C/c \cdot
g'(S) \le \lambda C/c \cdot g'(S)$. 
\end{itemize}

   b. Assume $S\ge 1$. We have to prove that $S/\lambda \cdot g'(S/\lambda) \le
   \bar c_2 g(S/\lambda)$, $\lambda \ge 1$. 
   \begin{itemize}
   \item If $1\le S \le \epsilon \lambda$ then $S/\lambda \cdot
     g'(S/\lambda) \le 1/ 
(\ln \epsilon)^2 \le 1/(1-\delta)(\ln \epsilon)^2 \cdot g(S/\lambda)$. 
   \item If $1 \le \epsilon \lambda \le S \le \lambda$ we have
     $S/\lambda \cdot g'(S/\lambda) \le 1\cdot C \le C/(1-\delta)
     \cdot g(S/\lambda)$. 
   \item If $S\ge \lambda$ we have $g'(S/\lambda)=1$ and $S/\lambda =
g(S/\lambda)$. 
\end{itemize}
Inequality  (\ref{eq:rescaling}) is now proved. 

\medskip 

We can now proceed to the proof of Claim~2. We have 

\begin{eqnarray*}
  \lefteqn{\parallel \dot v \parallel^2 \ = \ \int_{_{\widetilde S(t) \le 1}}
  |\dot v (t)|^2dt    + \int_{_{\widetilde S(t) > 1}} |\dot v (t)|^2dt }  \\
  & \le & 2 \int_{_{\widetilde S(t) \le 1}} \bar c \frac {\widetilde S^o}
  {S^o} |\dot x - JY\circ x|^2_\Omega + \bar c \frac {\widetilde S^o}
  {S^o} |JY\circ x|^2_\Omega + |\dot{\bar x}^{\tx{h}} - JY\circ
  x^{\tx{h}} |^2_\beta + |JY\circ x^{\tx{h}} |^2_\beta \\ 
  & & + 2 \int_{_{\widetilde S(t) > 1}} \frac 1 {\lambda^2} |\dot x -
  JY\circ x|^2_\Omega + \frac 1 {\lambda^2} |JY\circ x|^2_\Omega + 
  |\dot{\bar x}^{\tx{h}} - JY\circ
  x^{\tx{h}} |^2_\beta + |JY\circ x^{\tx{h}} |^2_\beta \\ 
  & \le & 2 \Big( 
  \bar c \parallel \dot x - JY\circ x\parallel ^2 + \parallel
  JY\circ x^{\tx{h}} \parallel ^2_\beta  +  
  \displaystyle {\int_{_{\widetilde S(t) \le 1}}} 
  \hspace{-.4cm} \bar c \ \frac {\widetilde S^o |JY\circ
    x|^2_\Omega} {S^o} +  \int _{_{\widetilde S(t) >
    1}} \hspace{-.4cm} \frac {|JY\circ x|^2_\Omega} {\lambda^2} 
\Big) .
\end{eqnarray*}

 We claim that all four terms in the last sum are bounded. The first
 equals $\parallel \mc Y\circ x \parallel$. The second is bounded by
 assumption (\ref{eq:horizontal bound on Y}) on $Y^{\tx{h}}$. In
 the third term we have $\widetilde S^o\le 1$ whereas $|JY\circ
 x|^2_\Omega /S^o$ is bounded, as seen in the proof of Claim~1.
 Finally, the fourth term is bounded by $\frac C {\lambda ^2} 
 \parallel \sqrt {S\circ x} \parallel_{L^2} ^2 = C\Lambda$ and Claim~2
 is proved.  

\medskip 

 After going to a subsequence we can now assume that $v_k\stackrel {C^0}
 \longrightarrow v$ with $v$ continuous. 

\medskip 

 {\it Claim~3.} We have $S\circ v >1$ if $\Lambda$ is big enough. 

\medskip 

 We prove Claim~3 after having proved Claim~4 below. 

\medskip 

  Let now $U\subset B$ be a ball 
  such that $v \subset \widehat E|_U$, hence $v_k \subset
  \widehat E|_U$ for $k$ large enough. The existence of $U$ follows
  from the contractibility of the $v_k$'s. Choose a radial 
  contraction of $U$
  onto $b\in U$ and consider the associated parallel transport
  $\tau:\widehat E|_U \longrightarrow \widehat E_b$. We define
  $\widetilde v = \tau \circ v$, $\widetilde v_k=\tau\circ v_k$,
  $\widetilde x_k = \tau \circ x_k$, $\widetilde {\bar x}_k = \tau
  \circ \bar x_k$. We clearly have $\widetilde x_k = (\widetilde {\bar
  x}_k,   S_k)$ and $\widetilde v_k = (\widetilde {\bar x}_k,  
S_k/\lambda_k^2)$.

\medskip 

 {\it Claim~4.} We have $\parallel \dot {\widetilde v}_k - X_{F(\widehat s)}(t,  
 \widetilde v) \parallel _\Omega \longrightarrow 0$.
 
 {\it Proof of Claim~4.} We write in this paragraph $F$ for
 $F(\widehat s)$ and we recall that, by condition~(\ref{eq:cond on
   graph of f}), the vector field $X_F$ 
  is vertical (and colinear to the Reeb vector field on the level
  surfaces of $F$). We have
 \begin{equation} \label{eq:homogeneity inequality} 
 \parallel \dot {\widetilde v}_k - X_F(\widetilde v) \parallel _\Omega = \frac 1
 {\lambda_k} \parallel \dot {\widetilde x}_k - X_F(\lambda_k^2
 \widetilde v) \parallel _\Omega  , 
 \end{equation} 
 due to the homogeneity property of the conical metric.  
 We also have   
 \begin{eqnarray*}
  \frac 1 {\lambda_k}    \parallel \dot {\widetilde x}_k -
  X_F(\lambda_k^2 \widetilde v) \parallel
   & \le & 
  \frac 1 {\lambda_k} \big( \parallel \dot {\widetilde x}_k - JY
  ^{\tx{vert}}({\widetilde x}_k)  \parallel
   +  
\parallel  (JY^{\tx{vert}} - X_F) ({\widetilde x}_k)
\parallel \big)\\
& &  + \ 
\frac 1 {\lambda_k} \parallel  X_F(\lambda_k^2 {\widetilde v}_k) -
X_F(\lambda_k^2 \widetilde v) \parallel  .
 \end{eqnarray*}
 
 The first term of the right hand side 
 is bounded by $\frac 1 {\lambda_k} \parallel \mc Y(x)
 \parallel $ and goes to zero as $k$ goes to infinity. The second term
 goes to zero with $k$ by (\ref{homotopie asymptotiquement lineaire}).
 The third term is bounded by $\parallel X_F(\widetilde v _k) -
 X_F(\widetilde v) \parallel
 _\Omega$ and therefore also goes to zero. This proves Claim~4.

  \medskip

  A direct consequence of Claim~4 is that $\dot{\widetilde v} =
  X_{F}(\widetilde v)$. Indeed, we have 
  $\dot{\widetilde v}_k \stackrel {L^2} \longrightarrow
  X_F(\widetilde v)$ and $\widetilde v_k \stackrel {L^2} \longrightarrow \widetilde
  v$. This implies $\dot{\widetilde v} \in L^2$ and $\dot{\widetilde
    v} = X_{F}(\widetilde v)$. 
  
  \medskip 

 {\it Proof of Claim~3.} We choose $\Lambda >1$ such that $\Lambda^2
 \cdot \min   f(s,   \bar u) > \max   f(s,   \bar u)$. 
 There exists $t_0$ such that $S(v(t_0)) = \Lambda^2$. Suppose by
 contradiction that $S(v(t)) = 1$ for some $t$ and let $t_1$ be the
 smallest such  $t$ with $S\circ v > 1$ 
 on $[t_0,   t_1[$ and $S(v(t_1))=1$. The same argument as in
 the proof of Claim~4 can be applied on the interval
 $[t_0,   t_1[$ 
 instead of $\Ss^1$ in order to show that $\dot{\widetilde v} =
 X_F(\widetilde v)$ on 
 this interval. In particular the image of $[t_0,   t_1]$ under
 $\widetilde v$ 
 is located on the level $S(\widetilde v(t_0)) f(\widehat s,  
 \bar{\widetilde v} (t_0))=
 \Lambda ^2 f(\widehat s,   \bar{\widetilde v} (t_0)) > f(\widehat s,
   \bar{\widetilde v} (t_1))$. In 
 particular 
 $S(\widetilde v(t_1)) f(\widehat s,   \bar{\widetilde v} (t_1)) >
 f(\widehat s,   \bar{\widetilde v} (t_1))$, 
 which means $S(v(t_1)) >1$. This is a contradiction and Claim
 3 is proved.

  \medskip

  At this point we have exhibited a $1$-periodic orbit $\widetilde v$
  for $X_F$, living on an arbitrarily large level. This ensures
  by (\ref{croissance stricte si orbite periodique}) that
  $\DP{f}{s}(\widehat s) \ge \epsilon >0$. We shall put to work the
  hypothesis $\int_0^1 \DP{H}{s} \le c$ in
  order to derive a final contradiction and complete the proof of
  Proposition~\ref{borne H 1 2}. We first compute: 
 \begin{eqnarray*}  
 \lefteqn{\DP{H}{s}(s,   t,   \bar{x},   S)} \\
 & = & \int_0^1\frac d {d\gamma} 
 \DP{H}{s} (s,   t,   \bar{x} ,  \gamma S + 1-\gamma)   d\gamma + 
 \DP{H}{s}(s,   t,   \bar{x},   1) \\ 
 & = & \int_0^1 \DP{^2H}{s\partial S}(s,   t,   \bar{x} ,  \gamma S 
 + 1-\gamma)\cdot (S-1)   d\gamma + \DP{H}{s}(s,   t,   \bar{x},   
 1) \\ 
 & = & \int_0^1 \Big( \DP{^2H}{s\partial S}(s,   t,   \bar{x} ,  \gamma S 
 + 1-\gamma) - \DP{^2F}{s\partial S}(s,   t,   \bar{x} ,  \gamma S 
 + 1-\gamma)\Big)\cdot (S-1)   d\gamma  \\  
 & & \qquad + \ \DP{H}{s}  (s,   t,   \bar{x},   
 1) \ + \ \partial _s f(s)\cdot(S-1) . 
 \end{eqnarray*}  
 We have used $\DP{^2F}{s\partial 
 S}=\partial _s f$. On the other hand 
 condition~(\ref{second derivatives s and S}) implies, for
 any $\tau >0$, the existence of a constant  $c_\tau$ such that 
 $$\Big| \DP{^2H}{s\partial S} (s,   t,   \bar{x},   S) - 
 \DP{^2F}{s\partial S} (s,   \bar{x},   S) \Big| \le \tau + \frac 
 {c_\tau} {\sqrt{S}} .$$ 
 For $k$ large enough $s_k$ is close to $\widehat{s}$ and we get:   
 \begin{eqnarray*}  
 \lefteqn{c \ge \int_0^1\DP{H}{s}(s_k,   t,   \bar{x}_k(t),   S_k(t) 
 )   dt } \\ 
 & \ge & -\tau \parallel \sqrt{S_k} \parallel ^2 _{L^2} - c_\tau 
 \parallel \sqrt{S_k} \parallel _{L^1} - \ C + \frac \epsilon 2 \parallel  
 \sqrt{S_k} \parallel ^2 _{L^2} - \ \frac \epsilon 2  .
 \end{eqnarray*}  
 But $\parallel \sqrt{S_k} \parallel _{L^1} \le \parallel \sqrt{S_k} 
 \parallel _{L^2}$ and, for $\tau < \epsilon $, the right hand 
 term of the above inequality goes to $+\infty$ with $k$ because we
 have supposed $\parallel 
 \sqrt{S_k} \parallel _{L^2} \longrightarrow \infty$. This is the
 desired contradiction and Proposition~\ref{borne H 1 2} is proved. 
 \hfill{$\square$}

\medskip 

\begin{lem}[compare~\cite{FH94}, Lemma 9]  
\label{bornes uniformes sur l_action}  
Let $\widehat{J}$, $H$ and $Y$ satisfy conditions  
(\ref{standard S grand} -  \ref{homotopie J}) and 
(\ref{Y is pseudo grad} - \ref{pas d_orbites periodiques pour F et
  G}). There exist  $c_1 \le c_2$ such that every solution of 
\begin{eqnarray*}   
u_s + \widehat{J}(s,   t,   u) u_t -  Y(s,   t ,   u) =0 
 , \\ 
-\infty < \inf_{s \in \RR} A_{H(s)}(u(s)), \ \sup_{s \in \RR} 
A_{H(s)}(u(s)) < + \infty   
\end{eqnarray*}  
satisfies 
\begin{equation}  
A_{H(s)}u(s) \in [c_1,   c_2] . 
\end{equation}  
\end{lem}  
 
\medskip

\demo We have 

$$\frac d {ds} A_{H(s)}u(s) = - \int _0 ^1 \DP{H}{s}(s,   t,   u(s, 
   t))   dt - dA_{H(s)}(u(s)) \cdot  \mc Y(s,   u(s))   ,$$ 
where 
 $$ \parallel \xi \parallel ^2 _{s} \ = \ \int _0 ^1 
 \omega (\xi(t),   \widehat J(s,   t)\xi(t))   dt, \quad x: 
 \Ss^1 \longrightarrow \widehat{E}, \quad \xi \in \Gamma(x^* 
 T\widehat{E}) . $$ 
By (\ref{croissance de H}) we infer that $A_{H(s)}u(s)$ is decreasing
with $s$ and therefore 
 $$\lim_{s\rightarrow -\infty} A_{H(s)}u(s) = \sup_{s\in 
   \RR}A_{H(s)}u(s)  , $$
$$\lim_{s\rightarrow +\infty} A_{H(s)}u(s) = \inf_{s\in 
   \RR}A_{H(s)}u(s) .$$ 
We also have
$\int_{-\infty}^{+\infty} dA_{H(s)}(u(s)) \cdot  \mc Y(s,   u(s)) <
\infty $ and we get a sequence  $s_k \longrightarrow 
\infty$ such that 
$dA_{H(s_k)}(u(s_k))\cdot \mc Y(s_k,   u(s_k)) \longrightarrow 0$. 
We shall prove that $A_{H(s_k)}u(s_k)$ is bounded
from below 
by a universal constant $c_1$. The same argument applied to a sequence
$s_k\longrightarrow -\infty$ will yield the universal upper bound
$c_2$. 

We note at this point that, because
$s_k \longrightarrow \infty$, we can assume by (\ref{homotopie J}) and
(\ref{homotopie H}) that $\widehat J=J_+$, $H=H_+$, $Y=Y_+$ are all
independent of $s$. As usual, we denote $x_k = u(s_k)$, $S_k=S\circ
x_k$. The strong pseudo-gradient property holds therefore with the
uniform constant $a_+>0$ and we infer that 
$\parallel \mc Y(s_k,   x_k) \parallel ^2_{s_k} \longrightarrow 0$.

\medskip 

{\it Claim~1.} $\parallel \sqrt {S_k} \parallel _{L^2}$ is bounded by
a constant depending on $u$. 

\medskip 

Let us assume for a moment that Claim~1 is true. 
Combining it with the fact that 
$\parallel \mc Y(x_k) \parallel _{L^2}$ is bounded we get by Claim
1 of Proposition~\ref{borne H 1 2} a $H^1$-bound on $x_k$
which again depends on $u$. By Lemma~\ref{lem:compact embedding} we
can then find a subsequence still
denoted $x_k$ which converges uniformly to a continuous loop $x$. 

We know that $\mc Y(x_k) \stackrel {L^2} \longrightarrow 0$ and therefore
$x\in H^1$, $\mc Y(x)=0$ and $\dot x_k \stackrel {L^2} \longrightarrow
\dot x$. The fact that $\mc Y$ is a (negative) pseudo-gradient
for the action implies that $x$ satisfies the equation $\dot x =
X_{H_+}(x)$. By Proposition~\ref{borne H 1 2} applied to $H(s)\equiv
H_+$, $\mc Y(s) \equiv \mc Y_+$ we get a {\it
  universal} bound on the $H^1$ norm of $x$. This implies a universal
bound on its $C^0$ norm through Lemma~\ref{lem:compact embedding} and
a universal bound on its $C^1$ norm through the equation $\dot x =
X_{H_+}(x)$. Finally, a $C^1$-bound on a contractible loop implies a
universal bound on the action $A_{H_+}(x)$. Moreover, we have
$A_{H(s_k)}(x_k) \longrightarrow A_{H_+}(x)$. 

It is now clear that we can set 

$$c_1 = \inf \big\{A_{H_+}(x) \ : \ \dot x =X_{H_+}(x) \big\} > -\infty .$$

\medskip 

{\it Proof of Claim~1.} We suppose by contradiction that, up to
considering a subsequence, we have $\parallel \sqrt {S_k} \parallel_{L^2}
\longrightarrow \infty$. We shall derive a contradiction along the
lines of the proof of Proposition~\ref{borne H 1 2}. 

Let $\lambda_k = \parallel \sqrt {S_k} \parallel_{L^2} / \Lambda$, 
$\Lambda > 1$ and $v_k(t) = (\bar x_k(t),   S_k(t)/\lambda_k^2)$. We
are now precisely in the situation of (\ref{eq:hypothesis for
  contradiction}) in Proposition~\ref{borne H 1 2}, with $H\equiv H_+$
and $\mc Y \equiv \mc Y_+$. We get a subsequence still denoted $v_k$
which converges to a continuous limit $v$, giving rise to a
$1$-periodic orbit of $X_{F_+}$ which is located on a level set of
$F_+$ with {\it arbitrarily} large $S$ coordinate. This contradicts
hypothesis (\ref{pas d_orbites periodiques pour F et G}) and concludes
the proof of Claim~1 and of the Lemma. 
\hfill{$\square$}

\medskip 

\begin{lem}[compare~\cite{FH94}, Proposition 11] 
\label{borne norme delta} 
Let $\widehat{J}$, $H$ and $Y$ satisfy  conditions  
(\ref{standard S grand} -  \ref{homotopie J}) and 
(\ref{Y is pseudo grad} - \ref{controle derivee seconde}).
 For any $\delta >0$ 
there is a constant $c_\delta>0$ such that any solution of   
\begin{eqnarray*}   
u_s + \widehat{J}(s,   t,   u) u_t - Y(s,   t ,   u) =0 
 , \\ 
-\infty < \inf_{s \in \RR} A_{H(s)}(u(s)), \ \sup_{s \in \RR} 
A_{H(s)}(u(s)) < + \infty   
\end{eqnarray*}  
satisfies
\begin{equation}  
 [\sqrt{S \circ u}]_\delta \le c_\delta . 
\end{equation}  
\end{lem}  
 
 \demo Let $u:\RR \times \Ss^1 \longrightarrow \widehat E$ 
 satisfy the hypothesis. By Lemma~\ref{bornes 
 uniformes sur l_action} there are constants $c_1 < c_2$  
 such that $A_{H(s)}u(s) \in [c_1,   c_2]$, $s\in \RR$. Let 
 $\bar{c}=c_2-c_1$. We infer that, for any $a < b$, we have

\begin{eqnarray} \label{difference finie pour l_action}  
\lefteqn{\int_a^b \Big(  dA_{H(s)}(u(s)) \cdot \mc Y(s,  u(s))  
+ \int_0^1\DP{H}{s}(s,   t,   u(s,   t))   dt \Big)   ds} \hspace{5.3cm} \\
  & = & A_{H(a)}u(a) - A_{H(b)}u(b) \le \bar{c} . \nonumber
\end{eqnarray}

 Let $\widehat{s}_k = k\delta/4$, $\tau = \delta/16$. Inequality  
 (\ref{difference finie pour l_action}) applied to $a= 
 \widehat{s}_k-\tau$ and $b=\widehat{s}_k+\tau$, together with the
 density of the set $\{s \ : \ a(s) > 0\}$, yields the existence   
 of $s_k \in [\widehat{s}_k-\tau,   \widehat{s}_k +\tau ]$ such that
 $a(s_k)>0$ and such that  
 $x_k = u(s_k)$ satisfies 
 $$ dA_{H(s_k)}(x_k) \cdot 
\mc Y (s_k,   x_k) + \int_0^1\DP{H}{s}(s_k,   t,   x_k(t))   dt
\le \bar{c} \cdot 8/\delta .$$ 
 By the asymptotic behaviour of $a(s)$ there are only a finite number
 of intervals on which $a(s)$ is nonconstant. We deduce the
 existence of some $a_0>0$ such that $a(s_k)\ge a_0$ for all $k\in
 \ZZ$. 
 Proposition~\ref{borne H 1 2} implies the existence of a constant
 $d(\delta)$ such that $\parallel \sqrt{S\circ x_k} \parallel _{H^1}\ < \ 
 d(\delta)$. The compact embedding $H^1(\Ss^1,   \RR) 
 \hookrightarrow C^0(\Ss^1,   \RR)$ gives a constant 
 $c_\delta$ such that $\parallel \sqrt{S \circ x_k } 
 \parallel _{C^0} \le c_\delta$.  
 On the other hand $(s_k) 
 \in \Gamma _\delta$, where $\Gamma_\delta$ is defined by (\ref{Gamma
   delta}), and this implies $[\sqrt{S \circ u } ]_\delta \le c_\delta$.
 \hfill{$\square$}

\medskip 
 
\noindent  
{\small \it Proof of Theorem~\ref{estimation C 0}.} Let  
$\alpha(s,   t) = S\circ u(s,   t)$. By Lemma~\ref{borne norme 
delta} we have $[\alpha]_\delta \le 
(c_\delta)^2$ for any $\delta >0$. 
In view of Proposition~\ref{inequation elliptique}, it
is enough to show that $\alpha$ satisfies an equation of the form 
$$\Delta \alpha \ge -A - B\alpha  ,$$ 
with $A$, $B$ positive constants depending only on $\widehat{J}$, $H$
and $Y$ but not on $u$. Choosing $\delta$ such that 
$\delta^2B<\pi^2$ will yield a $C^0$ estimate on 
$\alpha$.  
 
We first express $\Delta 
\alpha$  in a suitable way. The trick of exhibiting in (\ref{astuce CFH})  
the term $\frac 1 2 (|u_s|_\Omega^2+|u_t|_\Omega^2)$ is borrowed from 
\cite{CFH95}. Let us consider $R\ge 1$ such that $\widehat{J}(s,   t)
= J_V (s,   t) \oplus \widetilde J_B(s,   t)$ is standard split  for $S\ge 
R^2$. Let $\Gamma = \{ (s,   t) : \ \alpha(s,   t) \ge R^2 \}$. 
We have $dS \circ \widehat J (s,   t) = -C(s,   t)S\Theta$ on
$\Gamma$ and the following 
hold.
\begin{eqnarray*}  
\alpha_s(s,   t)  & = & dS\cdot u_s(s,   t) \ = dS\cdot ( 
-\widehat Ju_t+Y(s,  t,   u)) \\  
  & = & C(s,   t)S(u(s,   t) ) \Theta(u_t) + dS \cdot Y(s,   t,   u) . 
\end{eqnarray*}  
\begin{eqnarray*} 
\alpha_t(s,   t)  & = & dS\cdot u_t(s,   t) \ = dS\cdot ( 
\widehat Ju_s-\widehat JY(s,  t,   u)) \\  
  & = & -C(s,   t)S(u(s,   t) ) \Theta(u_s) + C(s,   t)S(u(s,  
  t)) \Theta( Y(s,   t,   u) ) . 
\end{eqnarray*}  
\begin{eqnarray} \label{astuce CFH}  
\lefteqn{ \partial_s\big( S(u(s,   t) ) \Theta(u_t) \big) - \partial
  _t \big(  
S(u(s,   t) ) \Theta(u_s) \big) } \\ 
& = & u_s\big( (S\Theta)(u_t) \big) - u_t\big( (S\Theta)(u_s) \big) 
\nonumber  \\ 
& = & d(S\Theta) (u_s,   u_t) + (S\Theta)\big( [u_s,   u_t] \big) 
\nonumber  \\ 
& = & d(S\Theta) (u_s,   u_t)  \ = \ \Omega(u_s,   u_t) \ = \ \frac  
1 2 \Omega(u_s,   u_t)  + \frac 1 2 \Omega(u_s,   u_t)   \nonumber  \\ 
& = & \frac 1 2 \Omega(u_s,   \widehat Ju_s - \widehat JY) - \frac 1
2 \Omega(u_t,    
-\widehat Ju_t + Y)  \nonumber \\  
& = & \frac  1 2 \big( |u_s|_\Omega^2 + |u_t |_\Omega ^2 \big) - \frac
1 2 \Omega  
(u_s,   \widehat JY) - \frac 1 2 \Omega ( u_t,   Y) .  \nonumber 
\end{eqnarray}  

We get 
\begin{eqnarray*}  
\lefteqn{\Delta\alpha \ = \ 
\frac {C(s,   t)} 2 \big( |u_s|_\Omega^2 + |u_t |_\Omega ^2 \big)
  - \frac {C(s,   t)} 2 \big(  \Omega  
(u_s,   \widehat JY) +  \Omega ( u_t,   Y) \big)} \\  
& & \quad \ + \ C_s (s,   t) S  \Theta(u_t) - C_t (s,   t)
S \Theta(u_s) + C_t(s,   t) S \Theta(Y) \\ 
& &  \quad \ + \ dS\cdot \partial _s Y(s,   t,   u)+ \ dS \cdot
\nabla _{u_s} Y(s,  
  t,   u) \\ 
& &  \quad \ + \ 
C(s,   t) \Big( (dS\cdot u_t) \cdot \Theta(Y) + S\Theta \big(
\partial _t Y(s,   t,   u) \big) + S\Theta  
\big( \nabla _{u_t} Y(s,  t,   u) \big) \Big) . 
\end{eqnarray*} 

\medskip 

We estimate now the terms composing the right side of the
above identity. Condition (\ref{homotopie J}) implies that $C(s,  
t)$ is independent of $s$ for $|s|\ge s_0$. As $C(s,   t)>0$ we get
the existence of strictly positive constants $C_0,   C_1$ such that 
$$0 < C_0 \le C(\cdot, \cdot) \le C_1, \qquad |dC(\cdot, \cdot)| \le C_1 .$$

We get 
$$\frac {C(s,   t)} 2 \big( |u_s|_\Omega^2 + |u_t |_\Omega ^2 \big)
\ge \frac {C_0} 2 \big( |u_s|_\Omega^2 + |u_t |_\Omega ^2 \big) .$$

On the other hand we have 
$$|\Omega(u_s,   \widehat JY) | \le |u_s|_\Omega\cdot |Y|_\Omega,
\qquad |\Omega ( u_t,   Y) | \le |u_t|_\Omega\cdot |Y|_\Omega .$$

We can estimate $|Y|_\Omega$ through the first condition of
(\ref{homotopie asymptotiquement lineaire}), which implies $|Y|_\Omega
\le |Y| _\omega \le c_1(1+\sqrt S) + |X_F|_\omega$. Now a direct
computation shows that $|X_F|_\omega = O(\sqrt S)$ and this gives
$|Y|_\Omega \le \bar c _1 (1+ \sqrt S)$. As a consequence we have 
\begin{eqnarray*} 
|\Omega(u_s,   \widehat JY) | & \le & \bar c_1 (1+ \sqrt \alpha)
|u_s|_\Omega, \\
|\Omega ( u_t,   Y) | & \le &  \bar c_1 (1+ \sqrt \alpha) |u_t|_\Omega
. 
\end{eqnarray*} 

The norm of $\Theta(\bar{u},   S)$ as a linear map is equal to 
$1/\sqrt{S}$ for $S \ge 1$ and we get 
\begin{eqnarray*}
  |C_sS\Theta(u_t)| & \le & C_1(1+\sqrt \alpha)|u_t|, \\ 
|C_tS\Theta(u_s)| & \le & C_1(1+\sqrt \alpha)|u_s|, \\
|C_tS\Theta(Y)| & \le & C_1\bar c_1 (1+\sqrt \alpha)^2 .
\end{eqnarray*}

The norm of $dS(\bar{u},   S)$ as a linear map is equal to 
$\sqrt{S}$ for $S\ge 1$, while\break $|\partial _s
\nabla F(s)|_\omega = O(\sqrt S)$. The second condition in (\ref{homotopie 
asymptotiquement lineaire}) implies $|\partial _s Y |_\omega = O(\sqrt
S)$ and we obtain 
$$| dS\cdot \partial _s Y(s,   t,   u)| \le 
c_2(1+\sqrt{\alpha})^2 .$$ 

We have as well 
$$|  (dS\cdot u_t) \cdot \Theta(Y(s,   t,   
u)) |  \le \bar c_1 (1+\sqrt{\alpha})|u_t|_\Omega  $$ 
and, by (\ref{controle dt}):  
$$|  S\Theta ( \partial _t Y(s,   t,   u) )   | \le 
c(1+\sqrt \alpha)^2 . $$ 

Finally, condition (\ref{controle derivee seconde}) gives   
\begin{eqnarray*} 
|  S\Theta\big( \nabla _{u_t} Y(s,  t,   u) \big)    | & \le &  c  
(1+\sqrt{\alpha} ) |u_t|_\Omega  ,  \\ 
| dS\cdot \nabla _{u_s} Y(s,   t,   u) | & \le &  c (1+\sqrt
\alpha) |u_s |_\Omega .
\end{eqnarray*} 

\medskip 

The modulus of the sum of the terms other than $\frac 1 2 C(s,  
t)(|u_s|_\Omega^2+|u_t|_\Omega^2)$ is therefore bounded by 
$$\bar C_1(1+\alpha) + \bar C_2(1+\sqrt{\alpha})|u_t|_\Omega + \bar C_3 (1+\sqrt{\alpha}) 
|u_s|_\Omega  ,$$
with obvious constants $\bar C_1$, $\bar C_2$ and $\bar C_3$. This
implies 
$$\Delta \alpha \ge -A -B\alpha  ,$$ 
with suitable constants $A$ and $B$. 
This inequality holds for 
$\alpha(s,   t) \ge R^2$. In order to get a global inequality on 
$\RR \times \Ss^1$ we use a trick that we borrow from 
\cite{FH94}. Let $\varphi:\RR_+ 
\longrightarrow \RR_+$ be a smooth function 
such that $\varphi (S) \equiv 0$ for $S \le 
R^2$, $\varphi'(S) =1$ for $S\ge R^2+1$ and $\varphi''(S)>0$ for 
$R^2 < S < R^2+1$. Then $\varphi$ satisfies  
\begin{equation} \label{propriete de phi} 
S\le \varphi(S) +\bar C 
\end{equation}   
for a suitable constant $\bar C$. Let  $\beta(s,   t) = \varphi \circ 
\alpha(s,   t)$. Inequality (\ref{propriete de phi}) gives a
bound on $\alpha$ in terms of a bound on $\beta$. On the other hand we
obviously have  $[\beta]_\delta \le [\alpha]_\delta$, 
$\delta >0$. It is therefore enough to show that $\beta$ satisfies
an inequality of the form $\Delta \beta \ge -A' - B'\beta$, with $A' \ge
0$, $B' \ge 0$ positive constants. Indeed
\begin{eqnarray*} 
\Delta \beta & = & \partial _s \big( \varphi'(\alpha(s,  t))\cdot 
\alpha_s \big) + \partial _t \big( \varphi'(\alpha(s,  t))\cdot 
\alpha_t \big) \\ 
 & = & \varphi''(\alpha(s,  t))\cdot\big( (\alpha_s)^2+(\alpha_t)^2 
\big) + \varphi'(\alpha(s,   t)) \cdot \Delta \alpha \\ 
& \ge & \varphi'(\alpha(s,   t)) (-A-B\alpha) \ \ge \ -A-B\alpha \ 
\ge \ -A -B\beta . 
\end{eqnarray*}  
\hfill{$\square$}

\begin{rmk} {\bf (On the maximum principle)} {\rm \label{rmk:split hamiltonians} 
  A priori $C^0$-bounds for Hamiltonians that are linear at infinity
  can be obtained directly through the maximum principle as
  in~\cite{Viterbo99}. Nevertheless, the interested reader can
  convince himself that such a direct approach is not effective for
  Hamiltonians of the form $h(S) + \widetilde f$. The solutions of
  Floer's equation satisfy in this case a second order elliptic
  equation with zero order term -- an avatar of which has already
  appeared in the previous proof -- and no general maximum principle is
  available in this context. 
}
\end{rmk}


\subsection{Definition of Symplectic cohomology} \label{sec:defiSH}
We define in this mainly expository 
 section the Symplectic cohomology groups $FH^*_{[a,b]}(E)$. We
 combine the philosophy of~\cite{Viterbo99} with the setup of~\cite{FH94}.

 Let $\mc J(\widehat E)$ be the set of {\it admissible almost complex
 structures}, consisting of time-dependent elements $\widehat J$ such that the
 associated constant deformation $\widehat J(s)\equiv \widehat J$
 satisfies~(\ref{standard S grand}) in Definition~\ref{defi:admissible
    deformations a c structures}. Let $\mc J(\RR;   \widehat E)$ 
be the set of {\it admissible deformations of almost complex
    structures} given by Definition~\ref{defi:admissible
    deformations a c structures}. 

Given $\widehat J \in \mc J(\widehat E)$ we denote by $\mc {HY}(E, 
 \widehat J)$ the set of {\it admissible
 Hamiltonians and pseudo-gradient vector fields}, consisting of
 time-dependent pairs $(H,   Y)$ such that the constant deformation
 $(H(s),Y(s))\equiv (H,Y)$ satisfies Definition~\ref{defi:admissible
 Hamiltonian deformation} and $H|_{\Ss^1\times E} <0$. 
Given $\widehat J \in \mc J(\RR;   \widehat E)$ we denote by 
$\mc{HY}(\RR;   E,   \widehat J)$ the set of \emph{admissible deformations
  of Hamiltonians and pseudo-gradient vector fields}, consisting of
pairs $(H,   Y)$ satisfying Definition~\ref{defi:admissible
  Hamiltonian deformation} and such that $H(s)|_{\Ss^1 \times E} <0$,
$s \in \RR$. 

We also define
\begin{equation*}
\mc{JHY}(E) = \bigcup _{\widehat J} \ \{ \widehat J \} \times
\mc{HY}(E,   \widehat J),  
\qquad 
\mc{JHY}(\RR;   E) = \bigcup _{\widehat J} \ \{ \widehat J \} \times
\mc{HY}(\RR;   E,   \widehat J).  
\end{equation*}

\begin{defi}
 Let $p>2$ be fixed. 
 We define $\mc{JHY}_{\tx{reg}}(E)$ to be the set of all triples
 $(\widehat J,  H,   Y)\in \mc{JHY}(E)$ such that the $1$-periodic
 orbits of $H$ 
 are nondegenerate and the linearized operator 
 \begin{equation*}
   D_u : W^{1,p}(u^*T\widehat E) \longrightarrow L^p(u^*T\widehat E) ,
 \end{equation*}
 \begin{equation*}
   D_u\xi = \nabla _s \xi + \widehat J _t(u)\nabla _s \xi + \nabla
   _\xi J_t(u) \partial _t u - \nabla _\xi Y_t(u) 
 \end{equation*}
 is surjective for any finite energy solution $u:\RR \times \Ss^1 
 \longrightarrow \widehat E$ of the equation 
 \begin{equation} \label{eq:Floer equation middle of paper}  
u_s + \widehat J(t,   u(s,  t))u_t = Y(t,   u(s,   t)) .
 \end{equation} 
\end{defi} 

It is proved in~\cite{Floer-SympFixedPoints,SZ92} that, if the
$1$-periodic orbits of $H$ are nondegenerate, the operator 
$D_u$ is Fredholm and its index at a solution $u$ is equal to 
$$
\mathrm{ind} \, D_u = i_{CZ}(x^+)-i_{CZ}(x^-), \qquad
x^\pm=\lim_{s\to\pm\infty}u(s,\cdot). 
$$
We have denoted by $i_{CZ}$ the Conley-Zehnder index as defined
in~\cite{Salamon-lectures}. The arguments in~\cite{FHS} can be adapted
to the pseudo-gradient setting and show that $\mc{JHY}_{\tx{reg}}(E)$
is of the second Baire category in $\mc{JHY}(E)$.   

Let $\mc{P}(H)$ be the set of $1$-periodic orbits of $H$. Given a
  regular triple $(\widehat J,   H,   Y)$ we define the complex   
  \begin{eqnarray}
    \label{eq:nontruncated Floer complex}
    FC^* (\widehat J,   H,   Y) 
  & = & \bigoplus _{\scriptsize \begin{array}{c} x \in \mc P (H) \\
  -i_{CZ}(x)=* 
      \end{array} } \ZZ \langle x \rangle , \\
    \label{eq:differential in Floer complex}
    \partial \langle x \rangle & = & \sum _{y \ : \ \dim \, \mc M(y,  
      x) = 1 } \#
    \big( \mc M (y,   x) / \RR \big) \ \langle y \rangle .
  \end{eqnarray}
Here $\mc M (y,   x)$ stands for the space of solutions $u$ of equation
(\ref{eq:Floer equation middle of paper}) with the limit conditions
$\displaystyle \lim_{s\rightarrow -\infty} u(s, \cdot) = y
(\cdot)$, $\displaystyle \lim_{s\rightarrow \infty} u(s, \cdot) = x
(\cdot)$. The additive group $\RR$ acts on $\mc M(y,x)$ by
reparametrizations, and we denote by $\#(\mc M(y,x)/\RR)$ the
algebraic number of elements of $\mc M(y,x)/\RR$ with respect to a
choice of coherent orientations~\cite{FH-coherent}. We claim that
$\partial ^2 =0$ and hence $\big( FC^*(\widehat J,   H,   Y),
\partial \big)$ is a differential complex.

The first important observation is that the equation $u_s+\widehat
Ju_t=Y(t,u(s,t))$ has the same analytic nature as the ordinary Floer
equation, namely the linearization $D_u$ is a compact perturbation of
the Cauchy-Riemann operator. In particular, a bound on the energy
$E(u)=\int_{\RR\times \Ss^1} |u_s|^2$ implies compactness up to breaking
of trajectories for the relevant moduli spaces. This in turn implies
$\p^2=0$. The second observation
is that such a uniform bound on the energy for the elements of $\mc
M(y,x)$ follows from the strong pseudo-gradient property
$dA_H\cdot \mc Y \ge c \| \mc Y\|^2$, $c>0$. Indeed, if $u_s=-\mc
Y(t,u(s,\cdot))$ and $u(-\infty,\cdot)=y(\cdot)$,
$u(+\infty,\cdot)=x(\cdot)$, we get  
\begin{equation*} 
A_H(y)-A_H(x) = -\int_\RR \frac d {ds} A_H(u(s,\cdot)) = \int_\RR
dA_H\cdot \mc Y \ge c \int_\RR \|\mc Y\|^2 = c  E(u),
\end{equation*} 
so that $E(u)\le (A_H(y)-A_H(x))/c$. 

\begin{rmk} {\bf (Convexity)} \rm Since the space 
of strong pseudo-gradient vector fields is convex, the homology groups
that are computed by means of any strong pseudo-gradient are the same
as the ones computed with the usual gradient of the action
functional. 
\end{rmk} 

The pseudo-gradient property implies that the action decreases along
  solutions of (\ref{eq:Floer equation middle of paper}), so that we
  have subcomplexes 
\begin{equation}
  \label{eq:truncated Floer complex}
  FC^*_{[a,   \infty[} (\widehat J,   H,   Y) =  \bigoplus
  _{A_H(x) > a} \ZZ \langle x \rangle \ \subset \ 
FC^*(\widehat J,   H,   Y)
\end{equation}
defined for $a\le \infty$, as well as quotient complexes 
\begin{equation}
  \label{eq:truncated quotient Floer complex}
  FC^*_{[a,   b]} (\widehat J,   H,   Y) 
  \ = \ FC^*_{[a, 
  \infty[}/FC^*_{[b,  \infty[} \ = \bigoplus
  _{A_H(x) \in ]a,   b]} \ZZ \langle x \rangle
\end{equation}
defined for $-\infty\le a < b \le \infty$. We denote the corresponding
  homology groups by $FH^*_{[a,   b]}(\widehat J,   H,   Y)$. They
  are endowed with natural restriction morphisms induced by inclusions
  of subcomplexes 
  \begin{equation*}
    FH^*_{[a,     b]}(\widehat J,   H,   Y) \longrightarrow   
    FH^*_{[a',     b']}(\widehat J,   H,   Y), \quad a \ge a',
    \quad b \ge b' . 
  \end{equation*}

The above homology groups depend on the asymptotic profile $F$ of the
Hamiltonian $H$ (see Definition~\ref{defi:admissible Hamiltonian
  deformation}). In order to define invariants of $\widehat E$ one needs to use 
an algebraic limit procedure which we describe below, and for which we
need to further restrict the class of admissible deformations
considered in Definition~\ref{defi:admissible Hamiltonian
  deformation} in order to be able to get a priori energy bounds
for solutions of the $s$-dependent Floer equation. 

\begin{defi} \label{defi:good} 
Let $\widehat J\in \mc J(\RR; \widehat E)$. A \emph{good deformation
  of Hamiltonians and pseudo-gradient vector fields} is an admissible
  deformation $(H,Y)$ such that\break 
  $H(s)|_{\Ss^1\times E} <0$, $s\in \RR$ and the following additional
  condition is satisfied.
\begin{itemize}
\item For any compact set
$K\subset \widehat E$ there exists a constant $a_K>0$ such that
\begin{equation} \label{eq:good-strong}
d A_{H(s,\cdot)}(x)\cdot \mc Y(s,x) \ge a_K \| \mc
Y(s,x)\|^2_{\widehat J(s)}
\end{equation} 
for all contractible loops $x:\Ss^1\to K$ and all $s\in\RR$. 
\end{itemize} 
\end{defi} 

\begin{rmk} {\bf (Admissible vs. good deformations)} \rm 
The difference between admissible and good deformations is that, for
the latter, $\mc Y(s,\cdot)$ might still be only a weak
pseudo-gradient for a nowhere dense set of values $s\in\RR$, but this
phenomenon is controlled in a precise way, namely 
$a_K\to 0$ as $K$ exhausts $\widehat E$. 
\end{rmk}

\begin{rmk} {\bf (Example)} \rm Given two triples $T_\pm=(\widehat
  J_{\pm},   H_\pm,   Y_\pm) \in\mc{JHY}(E)$ we can construct a good
  deformation connecting them as follows. We interpolate from $Y_-$ to
  $\nabla^{\widehat J_-}H_-$, then deform $\widehat J_-$ to $\widehat
  J_+$, $H_-$ to $H_+$ and implicitly $\nabla^{\widehat J_-}H_-$ to
  $\nabla^{\widehat J_+}H_+$, and then interpolate from
  $\nabla^{\widehat J^+}H_+$ to $Y_+$. Such deformations admit a
  uniform pseudo-gradient constant $a>0$ depending only on $T_\pm$.   
\end{rmk} 

\begin{rmk} {\bf (Degree of generality)} \rm We actually use in
  Section~\ref{sec:proofs} good 
  deformations which satisfy~\eqref{eq:good-strong}
  with a uniform constant, independent of the compact set
  $K$. Nevertheless, we
  chose to give the more general Definition~\ref{defi:good} in order
  to stress the fundamental role of the strong pseudo-gradient
  inequality. 
\end{rmk}

We denote by 
$$\mc{JHY}(E;  T_-,   T_+) \subset \mc {JHY}(\RR;  E)$$
the space of good deformations connecting $T_-$ and $T_+$. Standard
transversality methods allow one to define a space 
$$\mc{JHY}_{\tx{reg}}(E;   T_-,   T_+)$$
of regular good deformations, which is of second Baire category in 
$\mc{JHY}(E;  T_-,   T_+)$ and has the property that the
spaces of solutions of the equation
\begin{equation}
  \label{eq:Floer eq with param for monotonicity morphisms}
  u_s + \widehat J(s,  t,  u(s,   t))u_t = Y(s,   t,   u(s,  
t))
\end{equation}
are smooth manifolds. We claim that, just like in the case of constant
deformations, these spaces of solutions are compact modulo breaking
of trajectories. The main ingredient is an a priori energy
bound for the elements of the space $\mc M(x^-,x^+)$ of solutions
connecting $x^-\in \mc P(H_-)$ to $x^+\in\mc P(H_+)$, and
Definition~\ref{defi:good} plays a crucial role in obtaining it. Since
a good deformation is admissible, we know that Floer trajectories stay inside
a compact set $K$, hence $\mc Y$ is a strong
pseudo-gradient along Floer trajectories with uniform pseudo-gradient
constant $a_K>0$ and we have 
\begin{eqnarray*}
A_{H_-}(x^-) - A_{H_+}(x^+) & = & -\int_{\RR} \frac d {ds}
A_{H(s)}(u(s,\cdot)) \ = \ \int _{\RR} d A_{H(s)}\cdot \mc Y +
\int_{\RR\times \Ss^1} 
\frac {\partial H} {\partial s} \\ 
& \ge & a_K \int_{\RR} \|\mc Y \|^2 \ = \ a_KE(u).
\end{eqnarray*}
Therefore $E(u)\le (A_{H_-}(x^-) - A_{H_+}(x^+))/a_K$ is a priori
bounded. Note that, although we have used the hypothesis $\partial H/ \partial
s\ge 0$ in the above computation, this is not a crucial assumption at this
point. 

The pseudo-gradient property again implies that the symplectic action
  decreases along solutions
  of~\eqref{eq:Floer eq with param for monotonicity morphisms}, hence a good
  deformation induces chain maps
\begin{equation*}
  \sigma : FC^*_{[a,   \infty[}(\widehat J_+,   H_+,   Y_+)
  \longrightarrow FC^*_{[a,   \infty[}(\widehat J_-,   H_-,   Y_-)
  , 
\end{equation*}
\begin{equation*}
  \sigma \langle x^+ \rangle = \sum_{\scriptsize \begin{array}{c} x^-
  \in \mc P (H_-) \\ \dim \,  \mc M (x^-,   x^+) = 0 \end{array} } 
  \# \mc M(x^-,   x^+)   \langle x^-  \rangle
\end{equation*}
which pass to the quotient as chain maps 
\begin{equation*}
  \sigma : FC^*_{[a,   b]}(\widehat J_+,   H_+,   Y_+)
  \longrightarrow FC^*_{[a,   b]}(\widehat J_-,   H_-,   Y_-),
  \quad -\infty \le a < b \le \infty . 
\end{equation*}
Since the space of good deformations is convex the induced morphisms
in homology 
\begin{equation*}
  \sigma : FH^*_{[a,   b]}(\widehat J_+,   H_+,   Y_+)
  \longrightarrow FH^*_{[a,   b]}(\widehat J_-,   H_-,   Y_-),
  \quad -\infty \le a < b \le \infty
\end{equation*}
do not depend on the choice of the deformation. We call
them \emph{monotonicity morphisms}. We define the \emph{Floer or
  Symplectic cohomology groups} by
\begin{eqnarray}
  \label{eq:defi FH}
  FH^*_{[a,   b]}(E) & = & 
\lim _{\scriptsize \begin{array}{c} \longleftarrow \\
      {\prec} \end{array}} FH^*_{[a,   b]} (\widehat J,   H,   Y),
  \\
 \label{eq:last defi FH} 
FH^*(E) & = & \lim _{\scriptsize \begin{array}{c} \longleftarrow \\ b
   \end{array} } FH^*_{[-\infty,   b]}(E).
\end{eqnarray}
Here the partial order $\prec$ on $\mc {JHY}(E)$ is defined by 
\begin{equation}
  (J_-,   H_-,   Y_-) \prec (J_+,   H_+,   Y_+) \mbox{ iff }
  H_-(t,   x) \le H_+(t,   x), t\in \Ss^1,   x\in \widehat E , 
\end{equation}
and it makes $\mc{JHY}(E)$ into a directed set. 

We conclude this section with an invariance statement. 

\begin{thm}[\cite{Ci02}, Lemma 3.7, \cite{Viterbo99}, Theorem 1.7] 
\label{thm:invariance under def of symp form}
  Let $\beta_t$, $\Omega_t$, $Z_t$, $t\in [0,   1]$ be a deformation of
  the symplectic structure on the fibration $E$ such that 
  $$(E,   \pi,   B,   F,   \Omega_t,   Z_t,   \beta_t)$$ 
  defines a negative
  symplectic fibration in the symplectically aspherical   
  category for each $t\in
  [0,   1]$. We then have a natural isomorphism 
  $$FH^*(E;   \beta_0,   \Omega_0,   Z_0) \simeq FH^*(E;  
  \beta_1,   \Omega_1,   Z_1) .$$ 
\end{thm}


\subsection{Properties of Symplectic cohomology} \label{sec:properties of
  Floer homology} 

In this section all symplectic manifolds are assumed to be symplectically
aspherical. We recall the definition of {\it positive contact type
  boundary} given in the Introduction. 

\begin{defi} \label{defi:positive contact type}
  Let $(M,  \omega)$ be a symplectically aspherical manifold with
  boundary of contact type and Liouville form $\lambda$ defined in a
  neighbourhood of $\partial M$. We say that $\partial M$ is of 
  \emph{positive contact type} if every positively oriented
  closed contractible characteristic $\gamma$ has positive action
  $A_\omega(\gamma)$ bounded away from zero i.e. there exists $T_0 >0$
  such that 
   $$ 
   A_\omega(\gamma,\bar \gamma) = A_\omega(\gamma) :=
  \int_{D^2}\bar \gamma ^* \omega   \ge T_0,
   $$ 
  where $\bar \gamma : D^2 \longrightarrow M$ is any map satisfying $\bar \gamma
  |_{\partial D^2} = \gamma$. 
\end{defi}

\begin{rmk}{\bf (Examples)} {\rm 
1. Restricted contact type implies positive contact type in
view of the equality $\int_{D^2} \bar \gamma ^* \omega  
= \int _{\Ss^1} \gamma^* \lambda$. Boundaries of Stein domains are in
particular of positive contact type. 

2. If the boundary $\partial M$ is of contact type and 
has no closed contractible characteristics then it trivially satisfies the 
positive contact type property. 

3. Negative unit disc bundles satisfy the positive contact type
   property. In that case we have $\omega = \pi^*\beta + \Omega$ and
   $\lambda = (1+r^2)\theta$, where $\theta$ is the transgression
   $1$-form (see Example~\ref{sec:neg line bdles}). 
   The closed characteristics are contractible in the fibers 
   and we have, for each of them, $A_\omega(\gamma) = \int_{\Ss^1}
   \gamma^*\theta = \frac 1 2 \int_{\Ss^1} \gamma^* \lambda $. 

4. I do not know any example of a symplectically aspherical
manifold whose boundary is of contact type but not of positive contact
type. One should note that the definition only makes sense in a
symplectically aspherical manifold. If $\langle \omega,\pi_2(M)
\rangle \neq 0$ we can always glue a sufficiently negative sphere to
any filling disc $\bar \gamma$ so that $A_\omega(\gamma,   \bar
\gamma)$ becomes negative. 
}
\end{rmk} 

\begin{rmk} {\bf (Computation)} \rm If $\partial E$ has the positive contact type
    property as defined in Section~\ref{sec:properties of Floer
      homology} we can compute $FH^*(E)$ with a cofinal family
of Hamiltonians such that their $1$-periodic orbits are
      either constant with negative action close to zero, or
      nonconstant with positive action~\cite{Viterbo99}. Given any
      $a<0$ we then have  
$$
FH^*(E) = \displaystyle \lim_{\scriptsize \begin{array}{c} \longleftarrow \\ b \end{array} }
FH^*_{[a,   b]}(E).
$$
\end{rmk} 

Under the positive contact type assumption on the fibers of $E$, the
proofs of~\cite{Viterbo99}, based on manipulations 
of energy levels, apply verbatim
in order to show that our fibered version of symplectic homology
has the following properties. 

\renewcommand{\theenumi}{\alph{enumi}}
\begin{enumerate}
\item (\cite{Viterbo99}, Prop. 1.4) If $\mu>0$ is small enough we
  have 
  $$
FH^*_{]-\infty,    \mu]}(E) \simeq FH^*_{[-\mu,   \mu]}(E)
\simeq H^{*+m}(E,   \partial E),
\qquad 2m=\dim \, E.
$$
 In particular, there is a canonical
  morphism $c^*: FH^*(E) \to H^{*+m}(E,   \partial E)$ induced by the
  truncation of the range of action. 
\item (\cite{Viterbo99}, Thm. 3.1) 
  Any codimension $0$ embedding $j:W \hookrightarrow
  E$ of a domain $W$ such that $\partial W$ is of positive contact type
  and $W$ satisfies condition (A) of~\cite{Viterbo99} induces a
  transfer morphism $Fj^!: FH^*(W) \to
  FH^*(E)$ which makes the following diagram commutative
$$
\xymatrix{FH^*(W) \ar[r]^{Fj^!} \ar[d]_{c^*} & FH^*(E) \ar[d]^{c^*} \\
H^{*+m}(W,   \partial W) \ar[r]^{j^!} & H^{*+m}(E,   \partial E).
}
$$
The bottom arrow is the Poincar\'e dual of $H_{m-*}(W) \stackrel {j_*}
\longrightarrow H_{m-*}(E)$. The requirement that $\partial W$ be of
positive contact type can be 
  relaxed to the weaker assumption that $\partial W$ 
  be the boundary of a negative symplectic fibration whose fibers
  satisfy the positive contact type condition. Moreover, if there
  exists no closed characteristic on $\partial W$ or if $\partial W$
  is of restricted contact type \emph{in} $M$, then $Fj^!$ is
  defined in the symplectically aspherical case without reference to
  the additional condition~(A) of~\cite{Viterbo99}.
\item (\cite{Viterbo99}, Thm. 4.1) If the map $FH^*(E) \longrightarrow
  H^{2m}(E,   \partial E)$ is not surjective, then any contact type
  hypersurface which bounds a domain $W$ in $M$ carries a closed
  characteristic. The same conclusion holds if $\partial W$ is the
  boundary of a negative symplectic fibration.
\end{enumerate}

\begin{rmk} \rm {\bf (Homological vs. cohomological formalism)} \label{rmk:homology}
As announced in the introduction, one can build symplectic {\it homology}
groups based on the same chain groups as the cohomological ones, but
with dual differential
$$\delta \langle y \rangle = \sum _{x \ : \ \dim \,  \mc M(y,  x)=1} \#
\big( \mc M(y,  x)/\RR \big)   \langle x \rangle .$$
This formula is to be compared with~(\ref{eq:differential in
    Floer complex}). The main difference between cohomology
  and 
homology is that the first involves an inverse limit, whereas the
second involves a direct limit. The latter is always an exact functor
and therefore the homological spectral sequence holds with
integer coefficients, whereas the cohomological one holds with field
coefficients. We chose to work with cohomology in order to 
respect the setting of~\cite{Viterbo99} on which we base our applications. 
\end {rmk}


\section{Pseudo-gradient vector fields} \label{sec:main section pseudo
  gr vector fields}

We construct in this section a special family of almost
complex structures, Hamiltonians and
pseudo-gradient vector fields $(J_\nu, K_\nu, Y_\nu)$, $\nu \to \infty$
which is cofinal for the previously defined order
$\prec$ and which satisfies the following two properties. 

\renewcommand{\theenumi}{\Alph{enumi}}
  \begin{enumerate}
  \item \label{first geometric condition for filtration} 
    The $1$-periodic orbits of $K_\nu$ are located
    in the fibers over the critical points of a function 
    $f:B\longrightarrow \RR$ which is $C^2$-small and whose gradient
    flow is Morse-Smale. The gradient is computed with
    respect to the metric induced by a generic time-independent 
    almost complex structure $J_B$. 
  \item \label{second geometric condition for filtration} 
    Floer trajectories for the Floer complex $FC^*(J_\nu,  
    Y_\nu)$ project on gradient trajectories of the Morse complex
    $FC^*(J_B,   c _\nu f)$ on $B$, for some $c_\nu >0$.
\end{enumerate}

  Conditions~(\ref{first geometric condition for
  filtration}) and~(\ref{second geometric condition for filtration})
  will be used in \S\ref{sec:construction of the sp seq} in order to filter 
  each Floer complex $FC^*(J_\nu,   Y_\nu)$ by the Morse index
  of the projections of the $1$-periodic orbits of $K_\nu$.

  
\subsection{A model pseudo-gradient property} 
We begin by introducing horizontal distributions which are
time-dependent. These play an important role in order to achieve
transversality within the class of split almost complex structures
(see Section~\ref{sec:main section transv}). We denote by  
$$
H_0=\tx{Vert}^{\perp _\Omega}
$$
the horizontal
distribution determined by $\Omega$, and by $H=(H_t)$, $t\in \Ss^1$
an arbitrary time-dependent horizontal distribution such that $H\equiv
H_0$ on $\{ S \ge 1 \}$. Every such horizontal distribution can be
described as a loop of graphs of linear maps 
$$L_t : H_0 \longrightarrow \tx{Vert}, \ t\in \Ss^1,$$
with $\tx{graph}(L_t) \subset H_0 \oplus \tx{Vert} = T\widehat E$ and
$L_t$ supported in $E$. 
Given a horizontal distribution $H$ we mark
the horizontal lift of objects on $B$ by the symbol $\widetilde{\
}$. In order to emphasize the horizontal distribution $H$ with respect
to which we construct the lift, we shall sometimes use
the superscript $H$. We point out that time-dependent horizontal
distributions produce time-dependent lifts, even if the objects on $B$
are time-independent. On the other hand, the property
of being vertical is independent of the choice of horizontal distribution.
We denote by 
$$
X=X^\tx{h}+ X^{\tx{v}}
$$
the decomposition according to the splitting $T\widehat E = H_0\oplus
\tx{Vert}$. 
We denote by 
$$
X=X^{\tx{horiz}}+X^{\tx{vert}}
$$
the decomposition of a vector $X\in T\widehat E$ according to the
splitting $T\widehat E = H \oplus \tx{Vert}$. We have
$X^{\tx{horiz}}_t = X^{\tx{h}} + L_t (X^{\tx{h}})$, $X_t^{\tx{vert}} =
X^{\tx{v}} - L_t (X^{\tx{h}})$.  


In the next statement we make the following notations: 
\begin{itemize} 
\item $J_V$ is an almost complex structure on
  $\tx{Vert}$ which is compatible with $\Omega|_{\tx{Vert}}$ and which
  is standard on $\partial E \times [1,\infty[$; 
\item $J_B$ is an almost complex structure on $B$ compatible with
  $\beta$ and such that $\Omega(\cdot,  \widetilde J_B  \cdot)$ is
  positive on the horizontal distribution
  $H_0$ for $S\ge 1$;  
\item $f:B\longrightarrow \RR$ is a Morse function with 
critical points $\{p_1,\ldots,  p_\ell\}$ and $U_i$, $i\in
\{1,\ldots,\ell\}$ are mutually disjoint open neighbourhoods of the $p_i$'s;  
\item $H$ is a horizontal distribution given by a loop $L=(L_t:H_0
     \longrightarrow \tx{Vert})$, $t\in \Ss^1$ supported in
     $E\setminus \bigcup _{i=1}^\ell \pi^{-1}(U_i)$; 
\item $h: \widehat E\longrightarrow \RR$ is a Hamiltonian with
vertical Hamiltonian vector field, linear for $S\ge 1$ with slope
$\lambda_{\tx{max}}$ (typically of the form $h=h(S)$).
\end{itemize} 

We use a superscript $\epsilon$ in order to emphasize that the
Hamiltonian vector fields or the Hamiltonian action are computed with
respect to $\omega_\epsilon=\pi^*\beta+\epsilon\Omega$. 

\begin{prop} \label{conj:ps gr f tilde}  
  Assume the almost complex structure $J_B\in J(B,  \beta)$ is
  time-independent, and assume that the maximal slope of $h$ satisfies
  the condition  
$$
\lambda_{\tx{max}}\notin \tx{Spec}(\partial E).
$$

There exist constants $\epsilon_0,\delta_0,\rho_0,\alpha_0>0$ such that,
  for $\epsilon\le \epsilon_0$, $\alpha\le\alpha_0$, $\delta \le
  \delta_0$ and $\parallel L \parallel _{C^0} \le \rho_0\delta$, the
  following statements hold true: 
\begin{itemize} 
\item the form $\omega_\epsilon$, $0<\epsilon\le \epsilon_0$ is
  nondegenerate on $\widehat E$ and tames $J_V\oplus\widetilde J_B^H$;
\item the vector field defined on the space of $1$-periodic loops by 
$$\mc Y^\epsilon(x) = 
J(\dot x - X^\epsilon_{\epsilon h} \circ x - \widetilde
{X_{\epsilon\delta f}} \circ x) 
$$
satisfies the strong pseudo-gradient inequality
  \begin{equation}
    \label{eq:first pseudo-gradient result}
    dA^\epsilon_{\epsilon(h+\delta \widetilde f)}(x) \cdot \mc Y^\epsilon(x)
 \ge \alpha
\parallel \mc Y^\epsilon(x) \parallel^2_{\omega_\epsilon},
  \end{equation}
  with equality iff $x$ is a $1$-periodic orbit of
  $X_h=X^\epsilon_{\epsilon h}$ in a critical fiber of $\widetilde
  f$. 
\end{itemize} 
\end{prop}

\begin{rmk} {\bf (The trivial case)} \rm
If $X_f\equiv 0$ the pseudo-gradient property is 
clearly satisfied with $\alpha=1$ since $\mc Y^\epsilon = \nabla
A^\epsilon_{\epsilon h}$.    
\end{rmk} 

\begin{rmk} {\bf (Comparing $\mc Y^\epsilon$ and $\nabla
    A^\epsilon_{\epsilon(h+\delta \widetilde f)}$)} \rm  
Let us assume in order to simplify notation that
$\epsilon=\delta=1$. Then $\mc Y(x)-\nabla A_{h+\widetilde
  f}(x)=(X_{\widetilde f}-\widetilde {X_f})\circ x$. On the other
hand, if $\Omega|_{\partial E}$ is nondegenerate on
$H_0$ we have $|X_{\widetilde f}|_\beta\to 0$ as $S\to\infty$, whereas
$|\widetilde {X_f}|_\beta$ stays constant as $S\to \infty$. The vector field
$\mc Y$ is thus a ``big'' perturbation of $\nabla A_{h+\widetilde f}$
and inequality~\eqref{eq:first pseudo-gradient result} should come as
a pleasant surprise. 
\end{rmk} 

\demo The statement concerning $\omega_\epsilon$ follows immediately
  from the {\sc (negativity)} assumption and from the fact that $L$ is
  supported in $E$, so that we are left to prove the statement concerning
  $\mc Y^\epsilon$. We note that $\mc Y^\epsilon(x)=
J(\dot x - X_h \circ x - \epsilon\widetilde
{X_{\delta f}} \circ x)$ and, in order not to burden the notation, we
  give the proof for $\epsilon=\epsilon_0=1$. The
  reader can easily convince himself that the proof holds verbatim for
  an arbitrary value $0<\epsilon<\epsilon_0$. The intuitive reason is
  that, as $\epsilon$ decreases, the factor in front of $\widetilde f$
  is allowed to vary in the smaller interval $]0,\epsilon\delta_0]$
  and the vector field $\mc Y^\epsilon$ gets closer to $\nabla
  A^\epsilon_{\epsilon(h+\delta\widetilde f)}$.  

Let $E=dA_{h+\widetilde f}(x)\cdot \mc
Y(x) = dA_{h+\widetilde f}(x)\cdot J(\dot x
- X_h - \widetilde{X_f})$. We have 
\begin{eqnarray*} 
  E & = & \int\omega(\dot x,   J \dot x - JX_h - J\widetilde
  {X_f} ) \ - \ \int (d\widetilde f + dh)\cdot (J\dot x - JX_h -
  J\widetilde {X_f} )  \\
 & = & \parallel \dot x - X_h \parallel
 ^2 _\omega - \int \Big( \omega(\dot x,   J \widetilde{X_f}) +
 \pi^*\beta(\widetilde{X_f},   J\dot x - J \widetilde{X_f}) 
  - \omega(X_h,   J\widetilde{X_f}) \Big) \\
 & = & \parallel \dot x - X_h - \widetilde {X_f} \parallel ^2 _\omega
 + \int \Big( \Omega(\widetilde{X_f},   J\dot x) - \Omega(\widetilde
  {X_f},   J \widetilde {X_f}) + \Omega(JX_h,   \widetilde {X_f}) \Big) .
\end{eqnarray*} 

We distinguish three cases: either the loop $x$ is contained in
$\partial E \times [1,   \infty[$, either it is contained in the
compact region $S \le \Lambda$, either it intersects both $S=1$ and $S
= \Lambda$. The real number $\Lambda \ge 1$ will be suitably chosen
below (we will see that $\Lambda=4$ is a convenient choice). 

\medskip 

{\it Case 1.} We suppose that $x$ is contained in $\partial E \times [1,  
\infty[$. Because $H=H_0$ on $\partial E \times [1,   \infty[$ 
the term $\Omega(JX_h,   \widetilde {X_f})$ vanishes and 
inequality~(\ref{eq:first pseudo-gradient result}) becomes  
\begin{eqnarray} \label{eq:1 for x outside boundary}
\lefteqn{(1-\alpha^2) \parallel (\dot x - X_h - \widetilde
 {X_f})^{\tx{h}} \parallel ^2 _\beta +} \nonumber \\
& + & \Big( (1-\alpha^2) \parallel \dot x - X_h - \widetilde {X_f}
 \parallel ^2 _\Omega 
 - \parallel \widetilde {X_f} \parallel ^2_\Omega + \int
 \Omega(\widetilde{X_f},   J\dot x) \Big)  \ \ge \ 0 . 
\label{eq:2 for x outside boundary}
\end{eqnarray} 
It is enough to prove that the term in~(\ref{eq:2 for x outside
  boundary}), which we denote by $E_1$, is positive. 
Let $x(t) = (\bar x (t),   S(t))$. 
Because $X_h$ has no $1$-periodic orbits on $\partial E\times
[1,  \infty[$ and $X_h$ is contained in a nondegeneracy subspace of
$\Omega$, there exists $c>0$ such that any $1$-periodic loop $\bar
x:\Ss^1\longrightarrow \partial E$ satisfies 
\begin{equation} \label{eq:c}
\int_{\Ss^1} |\dot{\bar x} - X_h |^2_\Omega \ge c .
\end{equation} 
We denote $\delta = \max _E | \widetilde {X_f} |_\Omega$. We let
$\eta>0$ be a positive number. The condition $E_1\ge 0$ is invariant 
under homotheties and we may
therefore assume that $\min_{\Ss^1}   S(t)=1$. 
We have
\begin{eqnarray*} 
  E_1 & = & 4(1-\alpha^2)\parallel {\sqrt S}^{  \prime} \parallel _{L^2}
  ^2 + \ (1-\alpha^2) \int_{\Ss^1} S(t) |\dot {\bar x} -
  X_h|^2_\Omega - \ 
  \alpha^2 \int_{\Ss^1} S(t) |\widetilde {X_f} \circ \bar x|^2_\Omega
  \nonumber \\ 
 & & - \ (2-2\alpha^2) \int_{\Ss^1} S(t)  \langle \dot {\bar x},  
 \widetilde {X_f} \rangle _\Omega  \ + \ \int _{\Ss^1} S(t)
 \Omega(\widetilde {X_f},  J\dot {\bar x}) \nonumber 
\end{eqnarray*}
\begin{eqnarray*}
& \ge &  4(1-\alpha^2)\parallel {\sqrt S}^{  \prime} \parallel _{L^2}
  ^2 \ + (1-\alpha^2)\int_{\Ss^1} S(t) |\dot {\bar x} - X_h|_\Omega^2 
  \ - \ \alpha^2\delta^2\parallel \sqrt S \parallel
  _{L^2} ^2 \\
& & - \ (1-\alpha^2 + \frac 1 2 \parallel \Omega \parallel _{E, 
  \Omega}) 
\int_{\Ss^1} S(t) \big(\eta|\dot{\bar x} - X_h|^2_\Omega
+ \frac 1 \eta |\widetilde {X_f}|^2_\Omega \big)  \\
& \ge & 4(1-\alpha^2)\parallel {\sqrt S}^{  \prime} \parallel _{L^2}
  ^2 + \ (1-\alpha^2 - (1-\alpha^2)\eta - \eta \parallel \Omega \parallel _{E, 
  \Omega}/2 ) \cdot c \\
& & - \ \delta^2 \big(\alpha^2 + 
  (1-\alpha^2  + \frac 1 2 \parallel \Omega \parallel _{E, 
  \Omega}) /\eta \big) \parallel \sqrt S \parallel
  _{L^2} ^2.
\end{eqnarray*}
This last expression is strictly positive for $\eta$ and
$\delta$ small enough due to the Poincar\'e inequality which, for
$\min_{\Ss^1}  S(t)=1$, writes  
$$\parallel \sqrt S \parallel _{L^2} \le 1+\parallel {\sqrt S}^{ 
  \prime} \parallel _{L^2} .$$
We note that, because we do not assume 
$\widetilde J_B$ to be compatible with $\Omega$ on $H$, the norm 
$\parallel \Omega \parallel _{E,\Omega}$ of
$\Omega$ as a bilinear map in the induced (possibly degenerate) metric
may be arbitrarily large. One can construct explicit examples
for this phenomenon.

\begin{rmk} {\bf (Slope)} \rm The above argument crucially uses the
  hypothesis that the maximal slope of $h$ does not belong to
  $\mathrm{Spec}(\partial E)$, through inequality~\eqref{eq:c}. 
\end{rmk}

{\it Case 2.} We suppose now that $x$ intersects both regions $S<1$ and
$S>\Lambda$, where $\Lambda$ is to be chosen later ($\Lambda=4$ is 
a suitable choice).  Let $J = \{t\in \Ss^1   :  
S(x(t))\ge 1 \}$ and $J^c=\Ss^1 \setminus J$. 
We can 
assume without loss of generality that $x$ has transverse
intersection with $\partial E$, in which case
$J$ is a finite union of intervals $J_k$, $k\in \{1, \ldots,   N\}$
and $S\circ x|_{\partial J_k} \equiv 1$. Let $\delta^2 = \displaystyle 
{\max_B |\beta( {X_ f},   J_B  {X_f})|}$. 

We must prove that $E' \ge 0$, where 
\begin{eqnarray*}
  \label{eq:5 loop x crossing boundary high up}
  E' & = & 
  (1-\alpha^2) \parallel \dot x - X_h - \widetilde {X_f}
 \parallel ^2 _\omega + \ \int
 \Omega(\widetilde{X_f},   J\dot x - JX_h - J\widetilde {X_f}) \\
& = & (1-\alpha^2)\int_{J^c}|\dot x - X_h - \widetilde {X_f}|
 ^2 _\omega + \int_{J^c} 
 \Omega(\widetilde{X_f},   J\dot x - JX_h - J\widetilde {X_f}) \\ 
   & & + \ (1-\alpha^2) \int_{J}|\dot x - X_h - \widetilde {X_f}|
 ^2 _\omega + \int_{J} 
 \Omega(\widetilde{X_f},   J\dot x - JX_h - J\widetilde {X_f}) .
\end{eqnarray*}

For $\epsilon_0$ small enough we have $\parallel \Omega
\parallel _\beta \le 1$ and $|\cdot|_\omega\ge \frac 1 2 
|\cdot|_\beta$ on $E$, hence $|\Omega(\widetilde{X_f},   J\dot x - J
X_h - J \widetilde {X_f})| \le \delta |\dot x - X_h -
\widetilde{X_f}|_\beta \le 2\delta |\dot x - X_h - \widetilde {X_f} |_\omega$. 
We therefore obtain 
\begin{eqnarray*}
  \lefteqn{(1-\alpha^2) \int_{J^c} |\dot x - X_h - \widetilde {X_f}|
 ^2 _\omega +  \int_{J^c} 
 \Omega(\widetilde{X_f},   J\dot x - JX_h - J\widetilde {X_f})}  \\
& \ge & \int_{J^c} \Big( \sqrt{1-\alpha^2} |\dot x - X_h -
 \widetilde{X_f}|_\omega -  \delta / \sqrt{1-\alpha^2} \Big)
 ^2 - \int_{J^c} \delta^2 / (1-\alpha^2) \\
& \ge & - \delta^2 /(1-\alpha^2) .
\end{eqnarray*}
On the other hand we have $x(t) \in \partial E \times [1,   \infty[$
for $t\in J$ and $H=H_0$ on $\partial E \times [1,   \infty[$. If
$\eta$ is a small enough positive real we get
\begin{eqnarray*}
  \lefteqn{(1-\alpha^2) \int_{J}|\dot x - X_h - \widetilde {X_f}|
 ^2 _\omega + \int_{J} 
 \Omega(\widetilde{X_f},   J\dot x - JX_h - J\widetilde {X_f})} \\
& \ge & 4(1-\alpha^2)\parallel {\sqrt S}^{  \prime} \parallel _{L^2(J)}
 ^2 + (1-\alpha^2 - \eta \parallel \Omega \parallel _{E,  \Omega} /2)
 \int_J  S(t) |\dot{\bar x} - X_h - \widetilde {X_f} |^2_\Omega \\
& & - \ \frac 1 {2\eta} \parallel \Omega \parallel _{E,  \Omega} \cdot
 \int_J S(t) |\widetilde {X_f} \circ \bar x |^2_\Omega \\
& \ge &  4(1-\alpha^2)\parallel {\sqrt S}^{  \prime} \parallel _{L^2(J)}
 ^2 - \frac {\delta^2 } {2\eta} \parallel \Omega \parallel _{E, 
 \Omega} \cdot \parallel \sqrt S \parallel ^2_{L^2(J)} .
\end{eqnarray*}
We denote $A=4(1-\alpha^2)$, $B=1/(1-\alpha^2)$, $C=\parallel \Omega
\parallel _{E,  \Omega}/2\eta$ and we have obtained
\begin{equation} \label{expr:intermediate for pseudo gradient}
E' \ge A\parallel \sqrt{S}^{\ \prime} \parallel _{L^2(J)} -
 B\delta^2 - 
C \delta^2 \parallel \sqrt S \parallel ^2_{L^2(J)} .
\end{equation}

The Poincar\'e inequality for a positive function $f$ 
defined on an interval
$I$ of length $a$ gives 
$$\parallel f \parallel ^2_{L^2(I)} \le
a(m+\parallel f ' \parallel _{L^1(I)} )^2 \le a(m+\sqrt a \parallel f
' \parallel _{L^2(I)} )^2  ,$$
where $m= \displaystyle{\min_I   f} \ge 0$. If $a_k$ denotes the
length of $J_k$, with $\sum_k a_k < 1$, we obtain 
\begin{eqnarray*}
  \parallel \sqrt S \parallel ^2_{L^2(J)} & \le & \sum_k a_k
  (1+\sqrt{a_k}\parallel {\sqrt S}^{\ \prime}\parallel  _{L^2(J_k)}
  )^2 \\
  & \le & 1 + 2\parallel  {\sqrt S}^{\ \prime}\parallel   _{L^2(J)} + 
\parallel  {\sqrt S}^{\ \prime}\parallel ^2  _{L^2(J)}  \ = \ (1+
\parallel  {\sqrt S}^{\ \prime}\parallel  _{L^2(J)}) ^2 .
\end{eqnarray*}

At this point we exploit the hypothesis on $x$ in order to produce a
lower bound 
on $\parallel{\sqrt S}^{\ \prime}\parallel_{L^2(J)}$. 
Let $I=[t_0,   t_1]
\subset J$ be an interval such that $S\circ x|_I\ge 1$,
$S(x(t_0))=1$ and $S(x(t_1))=M$, with $M = \max   S\circ x   >  
\Lambda$. We have
$$\parallel {\sqrt S}^{\ \prime}\parallel
_{L^2(J)} \ \ge \ \parallel {\sqrt S}^{\ \prime}\parallel _{L^1(J)} 
\ \ge \
\parallel {\sqrt S}^{\ \prime}\parallel _{L^1(I)} \ \ge \ 
\sqrt M - 1 \ > \ \sqrt \Lambda -1.$$
We claim that (\ref{expr:intermediate for pseudo gradient}) is
strictly positive if we choose $\Lambda = 4$, $\alpha\le \frac 1 2$ 
and $\delta$ such that 
$ \delta^2 B \le  1 $ and $\delta^2 C \le
\frac 1 4$. Indeed, 
we obtain $\parallel \sqrt S \parallel ^2_{L^2(J)} < 4 \parallel
{\sqrt S}^{\ \prime}\parallel  _{L^2(J)} ^2$ and the expression in 
(\ref{expr:intermediate for pseudo gradient}) is bigger than 
$\parallel{\sqrt S}^{\ \prime}\parallel  _{L^2(J)} ^2$, hence bigger
than $1$.

\medskip

{\it Case 3.} We suppose $x$ is contained in the region $\{ S \le
\Lambda \}$. We choose $\epsilon_0$ small enough so that 
$|\cdot |^2_\omega \ge \frac 2 3 |\cdot | ^2_\beta$ on $E$ and
$\parallel \Omega|_{H_0}\parallel _{\{S\le \Lambda\},  \beta} \le 
1/2$. We write 
\begin{eqnarray}
  E' & = & (1-\alpha^2) \parallel \dot x - X_h \parallel^2_\omega -
  2(1-\alpha^2) \langle \dot x - X_h,   \widetilde {X_f} \rangle
  _\omega + \int \Omega(\widetilde {X_f},   J\dot x) 
   \label{eq:first line for case 3}\\
& & + \ \parallel \widetilde {X_f}^{\tx{h}} \parallel^2_\beta -
  \alpha^2 \parallel \widetilde {X_f} \parallel ^2_\omega 
   + \int
  \Omega(JX_h,   \widetilde {X_f} ) .
  \label{eq:second line for case 3}  
\end{eqnarray}
 
  The choice of a small enough constant $\rho = \rho(f,  h)$ ensures 
$|\Omega(JX_h,   \widetilde {X_f})| \le \alpha ^2 | \widetilde {X_f}^{\tx{h}}
|^2 _\beta$ and $|\Omega(X_h ,   \widetilde {X_f} )| \le \alpha^2 |
\widetilde {X_f}^{\tx{h}} | ^2 _\beta$ pointwise. Let us argue for
$\Omega(JX_h,   \widetilde {X_f})$. The inequality is clearly true on
$\bigcup_{i=1}^k \pi^{-1}(U_i) \bigcup \partial E \times [1,   \infty[$,
where $\Omega(JX_h,   \widetilde {X_f}) = 0$. On $E \setminus \bigcup
_{i=1}^k \pi^{-1}(U_i)$ there exists $\eta>0$ such that $|\widetilde
{X_f}^{\tx{h}}|^2_\beta \ge \eta$. On the other hand 
$$|\Omega(JX_h,   \widetilde {X_f}) | = |\Omega(JX_h,   L(\widetilde
{X_f}^{\tx{h}} ))| \le \parallel \Omega \parallel_{\infty,E} \parallel X_h
\parallel _{\infty,E} \parallel L \parallel _{\infty,E} |\widetilde
{X_f}^{\tx{h}}|_{\beta},$$
where $\parallel \cdot \parallel _{\infty,E} = \max_E |\cdot
|_\beta$. By choosing $\rho < \alpha^2\eta/\parallel \Omega \parallel
_{\infty,E}\parallel X_h \parallel _{\infty,E}$ we get 
$|\Omega(JX_h,   \widetilde {X_f})| \le \alpha^2| \widetilde
{X_f}^{\tx{h}} |^2_\beta$. Because $\rho$ depends on $\eta$, 
the actual constant that we get
in the statement of the Proposition 
is of the type $\rho(f,h,\epsilon)\delta$. 
The same argument applies to $\Omega(X_h,   \widetilde {X_f})$. 

 By further diminishing $\rho$ we can achieve $|\widetilde {X_f}
 |^2_\beta \le 2|\widetilde {X_f} ^{\tx{h}} |^2_\beta$. Indeed, we have
 $$|\widetilde {X_f}|^2_\beta = |\widetilde {X_f}^{\tx{h}}|^2_\beta +
 |L(\widetilde {X_f}^{\tx{h}})|^2_\Omega \le (1+\parallel
 \Omega\parallel_{\{S\le \Lambda\},\beta} \parallel L \parallel^2
 _{\infty,E}) |\widetilde {X_f}^{\tx{h}}|^2_\beta \le 2|\widetilde
 {X_f}^{\tx{h}}|^2_\beta . $$
  We infer that, for $\epsilon_0$ and
  $\rho$ small enough, the expression in~(\ref{eq:second line
  for case 3}), which we denote by $E'_2$, satisfies 
  $$E'_2 \ge (1-4\alpha^2) \parallel \widetilde {X_f} ^{\tx{h}}
  \parallel ^2_\beta .$$

  We denote the expression in~(\ref{eq:first line for case 3}) by
  $E'_1$. It satisfies
  \begin{eqnarray*}
    E'_1 & = & (1-\alpha^2) \parallel (\dot x - X_h )
    ^{\tx{v}}\parallel ^2_\Omega - (1-2\alpha^2) \langle (\dot x - X_h
    )^{\tx{v}} ,   \widetilde {X_f} ^{\tx{v}} \rangle _\Omega \\ 
   & & + \ (1-\alpha^2) \parallel \dot x ^{\tx{h}} \parallel ^2_\omega -
    2(1-\alpha ^2) \langle \dot x^{\tx{h}},   \widetilde {X_f}
    ^{\tx{h}} \rangle _\omega + \int \Omega (\widetilde {X_f}^{\tx{h}}
    ,   J \dot x ^{\tx{h}}) \\ 
   & = & (1-\alpha^2) \parallel (\dot x - X_h)^{\tx{v}} - \frac
    {1-2\alpha^2} {2(1-\alpha^2)} \widetilde {X_f}^{\tx{v}} \parallel
    ^2_\Omega - \frac {(1-2\alpha^2)^2} {4(1-\alpha^2)} \parallel
    \widetilde {X_f} ^{\tx{v}} \parallel ^2_\Omega \\
   & & + \ \frac {1-\alpha^2} 2 \parallel \dot x ^{\tx{h}} \parallel
    ^2_\omega - 2(1-\alpha^2) \langle \dot x ^{\tx{h}},   \widetilde
    {X_f} ^{\tx{h}} \rangle _\beta \\
   & & + \  \frac {1-\alpha^2} 2 \parallel \dot x ^{\tx{h}} \parallel
    ^2_\omega - (1-\alpha^2) \int \Omega(\dot x ^{\tx{h}},   J
    \widetilde {X_f} ^{\tx{h}}) + \alpha^2 \int 
    \Omega (\widetilde {X_f} ^{\tx{h}},
      J\dot x ^{\tx{h}}), 
\end{eqnarray*}
hence
\begin{eqnarray*}
 E'_1  & \ge & - \frac {\rho^2} {4(1-\alpha^2)}  \parallel
    \widetilde {X_f} ^{\tx{h}} \parallel ^2_\beta \ + \ \frac {1-\alpha^2}
    3 \parallel \dot x ^{\tx{h}} \parallel ^2 _\beta - 2(1-\alpha^2)
    \langle \dot x ^{\tx{h}},   \widetilde {X_f} ^{\tx{h}} \rangle
    _\beta \\ 
   & & + \ \frac {1-\alpha^2} 3 \parallel \dot x ^{\tx{h}} \parallel
    ^2_\beta -  \frac 1 2 \int |\dot x ^{\tx{h}} |_\beta
    |\widetilde{X_f} ^{\tx{h}} |_\beta  \\
   & \ge & -\big( \frac 1 4 + \frac {\rho^2} 4 \big) 
\parallel \widetilde {X_f} ^{\tx{h}} \parallel
    ^2_\beta + \Big(\frac {1-\alpha^2} 3 - \frac 1 4  \Big) \parallel \dot x
    ^{\tx{h}} \parallel _\beta ^2 \ \ge \ -\big( \frac 1 4 + \frac
    {\rho^2} 4 \big) \parallel 
    \widetilde {X_f} ^{\tx{h}} \parallel ^2_\beta .
  \end{eqnarray*}

We have used for the first inequality that
$|\widetilde{X_f}^{\tx{v}}|_\Omega^2\le\rho^2|\widetilde{X_f}^{\tx{h}}|_\beta^2$,  
$\parallel\dot x^{\tx{h}}\parallel _\omega^2\ge \frac 2 3 {\parallel
\dot x^{\tx{h}}\parallel_\beta^2}$, and $\parallel
\Omega\parallel_{\{S\le\Lambda\},\beta}\le \frac 1 2$. 
We have used for the second inequality the fact that
$(1-\alpha^2)/3 \parallel \dot x ^{\tx{h}} \parallel ^2_\beta -
2(1-\alpha^2) \langle \dot x ^{\tx{h}},   \widetilde {X_f} ^{\tx{h}}
\rangle _\beta \ge 0$ if $\delta$ is small enough, by
Lemma~\ref{lem:ineq} below applied with $(W,g)=(B,g_{_B})$,
$g_{_B}(\cdot,\cdot)=\beta(\cdot,J_B\cdot)$ and $\eta=1/12$. We have
also used the inequality $|\dot x ^{\tx{h}} |_\beta |\widetilde{X_f}
^{\tx{h}} |_\beta \le \frac 1 2 (|\dot x ^{\tx{h}} |^2_\beta + |\widetilde{X_f}
^{\tx{h}} |^2_\beta)$. Finally, the last inequality holds if
$(1-\alpha^2)/3 - 1/4 \ge 0$. As a conclusion we obtain
$$E' \ge \big(\frac 3 4 -4\alpha^2-\frac {\rho^2} 4\big) \parallel \widetilde {X_f}
^{\tx{h}} \parallel ^2_\beta,$$
hence $E'\ge 0$ if $\alpha$ and $\rho$ are small enough. The equality case is
readily characterized from this last inequality.
\hfill{$\square$}

\begin{lem} \label{lem:ineq}
Let $(W,g)$ be a compact Riemannian manifold. If $\eta>0$ is small
enough, any time-independent 
$C^1$-vector field $X$ on $W$ with $\|X\|_{C^1}<\eta$ satisfies 
  \begin{equation}
    \label{eq:second fundam ineq on loops}
    \parallel \dot x \parallel _{L^2} ^2 \ \ge  \ \frac 1 {2\eta}
    \langle \dot x,   X\circ x \rangle _{L^2}  
  \end{equation}
  for any $1$-periodic loop $x$ of class $C^1$. Equality is achieved
  if and only if $x$ is a constant loop.
\end{lem} 

\proof We embed isometrically $(W,g)$ into Euclidean space and, for
$\eta$ small enough, we can extend $X$ to a
vector field supported near the image of $W$ and whose derivative
is pointwise bounded by $2\eta$. We can thus assume without loss of
generality that $(W,g)$ is Euclidean space with the standard metric,
and $X$ is a compactly supported vector field such that
$\|dX\|_\infty\le 2\eta$. Since $\langle \dot x , \mathrm{ct.}\rangle
=0$ and since we only impose a hypothesis on $\|dX\|_\infty$, we can
further assume that $X(0)=0$ and $x(0)=0$. By the Poincar\'e
inequality we obtain 
$$|\langle \dot x,   X \circ x \rangle_{L^2}| \le \parallel \dot x
\parallel _{L^2} \parallel X\circ x \parallel _{L^2} \le \parallel dX
\parallel _\infty \cdot \parallel \dot x  \parallel _{L^2}^2 \le 2\eta
\parallel \dot x \parallel ^2_{L^2} .$$
The equality case is readily characterized. \hfill{$\square$}

\begin{rmk} {\bf (Small orbits)} \rm
Inequality~\eqref{eq:second fundam ineq on loops} implies in
particular $\|\dot x - X\circ x \|_{L^2} \ge \|X\circ x\|_{L^2}$. It
should therefore be seen as a quantitative expression of the
well-known fact that a vector field which is small-enough in
$C^1$-norm has no nonconstant $1$-periodic orbits. In particular,
inequality~\eqref{eq:second fundam ineq on loops} does not hold if $X$
is time-dependent.   
\end{rmk} 

\begin{rmk} {\bf (Asphericity)} \label{rmk:why symp asph is important} 
\rm
Case 3. in the proof of Proposition~\ref{conj:ps gr f tilde} is the one
crucially involving the fact that $J_B$ does not depend on time through
the use of 
Lemma~\ref{lem:ineq}. This is the main reason
for requiring the base to be symplectically aspherical in the
construction of the spectral sequence: it is the only case where
transversality for Floer's equation can be achieved within the class
of time-independent almost complex structures. 
\end{rmk}


\subsection{Definition of the geometric Hamiltonians} 
\label{sec:the geometric Hamiltonians}
We construct now the
Hamiltonians $K_\nu$ announced in the introduction of the
present section. We recall that the vertical coordinate $S$ was
defined only in a neighbourhood of $\partial E$ as $S \in
[1-\delta,   1]$. We abandon in this section the notation that we
have used in \S\ref{sec:proof of C0 bounds} and do not consider
anymore the function $S$ as being extended over $E$.

We first define Hamiltonians $H_\nu : \widehat E \longrightarrow \RR$,
$\nu \in \NN$ by the formula  
$$H_\nu = h_\nu + c _\nu \widetilde f .$$
Here $f:B \longrightarrow \RR$ is a Morse function on $B$, its lift to
$\widehat E$ is 
$\widetilde f = f\circ \pi$, while $(c_\nu)_\nu$ is a decreasing sequence of
strictly positive real numbers converging to zero, to be chosen below.

We define the function $h_\nu$ as follows. 
We fix a strictly increasing sequence $\lambda_\nu$ of
positive real numbers which do not belong to the period spectrum of
$\partial E$ and verify $\lambda_\nu \longrightarrow \infty$, $\nu
\longrightarrow \infty$. We choose $h_\nu$ to be a smooth function
which is constant on 
$\widehat E \setminus  \partial E \times [1- \frac \delta 2,   \infty
[$ and depends only on $S$ on $\partial E \times [1- \frac\delta
2,\infty[$. We shall use the notation 
$h_\nu(S)$ both for the function $h_\nu$ and for the
corresponding function on $[1-\frac \delta 2,   \infty[$. We denote by
$T_0$ the smallest element of $\tx{Spec}(\partial E)$ and use the
convention $T_0=\infty$ if $\tx{Spec}(\partial E) = \emptyset$. We impose
the following conditions on $h_\nu$: 
\begin{itemize} 
  \item $h_\nu (S)$ is strictly convex on $]1-\frac \delta 2,   1[$
    and $h'_\nu(1-\frac \delta 4) < T_0$;
  \item $h'_\nu(S) = \lambda_\nu$ for $S \ge 1-\frac \delta 8$;
  \item $h_\nu < 0$ on $E$.  
\end{itemize} 

We note that the smoothness assumption on $h_\nu$ implies that
$h'_\nu$ vanishes at infinite order at $S=1-\frac \delta 2$. 
We moreover require that 
\begin{itemize} 
  \item for every $\nu \in \NN$ we have $h_\nu +c_\nu
    \widetilde f <0$ on $E$;
  \item for $\nu < \nu'$ we have $h_\nu < h_{\nu'}$, and therefore 
    $H_\nu < H_{\nu'}$ for $c_\nu$ small enough.
\end{itemize} 

Because $X_{h_\nu}$ is vertical, 
Proposition~\ref{conj:ps gr f tilde} implies that the
vector field $(J_V\oplus \widetilde J_B)\big( \dot x - X_{h_\nu} \circ
x - \epsilon c _\nu \widetilde {X_f} \circ x \big)$, $x\in
\Lambda_0\widehat E$  is a strong
negative pseudo-gradient for the action functional
$A^\epsilon_{\epsilon H_\nu}$ if $\lambda_{\tx{max}} \notin
\tx{Spec}(\partial E)$ and $c_\nu$ is small enough.  
The lift $\widetilde {X_f}$ is considered here with respect to a horizontal
distribution which is close enough to $H_0$.

The $1$-periodic orbits of $H_\nu$ are all degenerate
and fall in two classes: 

\begin{enumerate} 

\item constants in the critical fibers
of $\widetilde f|_{\widehat E \setminus \partial E \times ]1-\frac
\delta 2,   \infty [}$;

\item nontrivial orbits in the critical fibers of
  $\widetilde f$, appearing in the region $\{ 1-\frac \delta 4 < S
  < 1-\frac \delta 8 \}$ and corresponding 
  to closed characteristics with period smaller than $\lambda_\nu$.
  The characteristics are understood to be parametrized by
  $X_\tx{Reeb}$. In the best of the situations, they are transversally
  nondegenerate. 

\end{enumerate}

We construct now Hamiltonians $K_\nu$ with nondegenerate
$1$-periodic orbits by perturbing $H_\nu$. We  need two kinds of
perturbations, corresponding to the above two types of orbits:

\begin{enumerate}
\item \label{perturbation 1} 
  a time-independent perturbation localized in a neighbourhood of the critical
  fibers of $\widetilde f|_{\widehat E \setminus \partial E \times ]1-\frac
\delta 2,   \infty [}$;

\item \label{perturbation 2} 
  a time-dependent perturbation localized in a neighbourhood of
  the nonconstant $1$-periodic orbits. 
\end{enumerate}

\medskip 

(\ref{perturbation 1}) Let $\{p_1, \ldots,   p_\ell \}$ be the
critical points of the Morse function $f: B \longrightarrow \RR$. We
denote $F_i = F_{p_i}$ and $S_i = S|_{F_i}$. 
We choose mutually disjoint open sets $U_i \ni p_i$ admitting
trivializations $\Psi_i :\pi^{-1}(U_i) 
\stackrel \sim \longrightarrow  U_i\times F_i$ such that
$S_i\circ \tx{pr}_2 \circ \Psi_i = S |_{\pi^{-1}(U_i)}$ and
$\Psi_* H_y=\bar H_y$ for all $y\in F_i$, where $\bar H$ is the trivial
horizontal distribution on $U_i \times F_i$. Such a
trivialization can be constructed by parallel transport along the radii
of a geodesic ball centered at $p_i$. This type of trivialization is
even a symplectic diffeomorphism in the fibers, but we shall not use
this fact. 
$$
\xymatrix{  
\Psi _i : \pi^{-1} ( U_i) \ar[rr]^\sim \ar[dr]_\pi && U_i\times       
F_i\ni (x,y) \ar[dl]^{\tx{pr}_1} \ar[dr]_{\tx{pr}_2} \ar@{-->}[rr] & & \RR \\
& U_i & & F_i \ar[ur]_{\varphi_i} &
}       
$$

We choose now functions $\varphi_i: F_i \longrightarrow \RR$ subject
to the following conditions: 

\begin{itemize}
\item $\varphi_i = 0$ for $1-\frac \delta 4 \le S \le 1$;
\item $\varphi_i = \varphi_i(S)$ for $1-\frac \delta 2 \le S < 1-\frac
  \delta 4$, a strictly concave function satisfying $|\varphi_i'(1-\frac
   \delta 2)| < T_0$; 
\item $\varphi_i$ is a Morse function on $E \setminus \partial E
  \times [1-\frac \delta 4,   1]$.
\end{itemize}

We denote by $f_i$ the composition 
$\varphi_i \circ \tx{pr}_2 \circ \Psi_i :
\pi^{-1}(U_i) \longrightarrow \RR$.  
We choose relatively compact open subsets $V_i \Subset U_i$,
$p_i \in V_i$  and smooth compactly supported 
cut-off functions $\rho_i : U_i \longrightarrow \RR$ such that
$\rho_i|_{V_i}=1$ and $0 \le \rho_i \le 1$.

We now define the first perturbation $\widetilde K_\nu$ of $H_\nu$
to be
\begin{equation} \label{eq:the first perturbation}
\widetilde K_\nu = H_\nu + c'_\nu \sum_{p_i \in \tx{Crit}(f)}
\widetilde \rho _i \cdot f_i .
\end{equation} 
Here $c'_\nu$ is a decreasing sequence of strictly positive
real numbers, with $c'_\nu$  small enough
such that the only critical points of $\widetilde K_\nu$ inside $E$ are
the critical points of $f_i$ in $F_i$.

\medskip 

(\ref{perturbation 2}) The Hamiltonian $\widetilde K_\nu$
has nontrivial $1$-periodic orbits in the critical fibers $F_i$ in the
region $1-\frac \delta 4 < S < 1$. For any $\delta>0$ there exists a 
time-dependent Hamiltonian $\chi^\delta_{\nu,  i}: \Ss^1 \times F_i
\longrightarrow \RR$, with $\parallel {\chi}^\delta_{\nu,   i}
\parallel_{C^2}\le \delta$ and supported in an arbitrarily small neighbourhood
of the nontrivial $1$-periodic orbits of $\widetilde K_\nu|_{F_i}$, 
such that the $1$-periodic orbits of 
$\widetilde K_\nu|_{F_i} + {\chi}_{\nu,   i}^\delta $ are
nondegenerate. We denote
$G_{\nu,   i}^\delta = {\chi}_{\nu,   i}^\delta \circ \big(\tx{id}
  \times   (\tx{pr}_2 \circ \Psi_i)\big) : \Ss^1 \times \pi^{-1}(U_i)
\longrightarrow \RR$. We define 
\begin{equation} \label{eq:the second perturbation}
K_\nu = \widetilde K_\nu + c '_\nu \sum_{p_i \in
  \tx{Crit}(f)} \widetilde \rho_i \cdot G_{\nu,   i}^\delta .
\end{equation}

The property $\Psi_*H_y=\bar H_y$, $ y \in F_i$ ensures
that $X_{G_{\nu,   i}^\delta}$ and $X_{f_i}$ are vertical along
the fibers $F_i$. We infer the existence of a constant $m>0$ depending
on all the choices made before such that, for a given $\nu$, we have 
\begin{eqnarray*} 
| X_{G_{\nu,   i}^\delta} ^{\tx{h}} |_\beta  & \le & 
m | \widetilde {X_f}^{\tx{h}} |_\beta , \\ 
| X_{f_i} ^{\tx{h}} |_\beta  & \le & m |
\widetilde {X_f}^{\tx{h}} |_\beta .
\end{eqnarray*} 
Moreover, by multiplying $G_{\nu,   i}^\delta$ and $f_i$ by
sufficiently small positive reals we can achieve that the above
two inequalities hold for all $\nu$ and $i$ with a uniform constant
$m>0$, which can moreover be chosen arbitrarily small.


\subsection{Geometric pseudo-gradient vector fields}

The objects that we consider in the next statement are those of
Proposition~\ref{conj:ps gr f tilde}, namely an almost complex structure
$J=J_V\oplus \widetilde J_B$, a Morse function $f:B\to\RR$, a horizontal
distribution $H$ given by a loop $L=(L_t:H_0\to\mathrm{Vert})$, and a
Hamiltonian $h$ with vertical Hamiltonian vector field. We consider in
addition a time-dependent pertubation $G$ supported in $E \cap
\bigcup_{i=1}^\ell \pi^{-1}(U_i)$ satisfying the inequality
  \begin{equation} \label{eq:condition for the small perturbation} 
  | X_G^{\tx{h}} |_\beta \le m | \widetilde
  {X_f}^{\tx{h}} | _\beta
  \end{equation}
for some $m>0$. 

\begin{prop}\label{prop:general ps gr}
  Assume the almost complex structure $J_B\in\mc J(B,\beta)$ is
  time-independent, and assume that the maximal slope of $h$ satisfies
  the condition 
$$
\lambda_{\tx{max}} \notin \tx{Spec}(\partial E).
$$
There exist constants $\epsilon_0,\delta_0,\rho_0,\alpha_0,m_0>0$ such that, for
$\epsilon\le\epsilon_0$, $\alpha\le \alpha_0$, $\delta\le \delta_0$,
$\parallel L \parallel _{C^0} \le \rho_0\delta$ and $m\le m_0$, the
vector field defined on the space of $1$-periodic loops by 
$$\mc Y^\epsilon(x) 
= J(\dot x - X^\epsilon_{\epsilon h}\circ x - \widetilde {X_{\epsilon
    \delta f}} \circ x - 
(X^\epsilon_{\epsilon G})^{\tx{v}} \circ x )
$$
satisfies the strong pseudo-gradient inequality
\begin{equation}
    \label{eq:pseudo-gradient result in general}
    dA^\epsilon_{\epsilon(h+\delta \widetilde f+G)}(x) \cdot \mc Y^\epsilon(x)
\ge \alpha \parallel \mc Y^\epsilon(x) \parallel
    ^2_{\omega_\epsilon}. 
\end{equation}
  Equality holds iff $x$ is a periodic orbit of
  $X_{h+G}=X^\epsilon_{\epsilon(h+G)}$ in a critical fiber of
  $\widetilde f$.  
\end{prop}

\demo We follow the proof of Proposition~\ref{conj:ps gr f tilde} and
we again assume without loss of generality that $\epsilon=1$ in order
not to burden the notation. 
Because $G$ is supported in $E$ the proof of Case 1. remains
unchanged. The specific feature of Case 2. is the Poincar\'e
inequality for loops and this remains unchanged as well, although some
new estimates are needed in the preliminary computations. These
estimates appear also in the proof of Case 3. and we give full details
only for this last case. 
We denote $E=dA_{h+\widetilde f + G}\cdot \mc Y(x)$ 
and assume that $x$ is contained in
$\{ S \le \Lambda \}$. We have 
\begin{eqnarray*} 
E & = & \int \omega\big(\dot x ,  J(\dot x - X_h - \widetilde {X_f} -
X_G^{\tx{v}})\big) \ - \ d(h + \widetilde f + G)\cdot J(\dot x - X_h -
\widetilde {X_f} - X_G^{\tx{v}})  \\
& =  & \parallel \dot x - X_h - \widetilde {X_f} - X_G^{\tx{v}}
 \parallel ^2 _\omega -
 \int \Omega(\widetilde {X_f},   J\widetilde {X_f})  + \int
 \Omega(\widetilde{X_f},   J\dot x - J X_h)  \\ 
& & - \ \int \omega(X_G^{\tx{h}},   J\dot x - J \widetilde {X_f}) .
\end{eqnarray*} 

Inequality (\ref{eq:pseudo-gradient result in general}) is equivalent
to  
\begin{eqnarray*}
\lefteqn{(1-\alpha^2)\parallel \dot x - X_h - \widetilde {X_f} -
 X_G^{\tx{v}} 
 \parallel ^2 _\omega -
 \int \Omega (\widetilde {X_f},   J\widetilde {X_f})}  \\
 & & + \ \int
 \Omega(\widetilde{X_f},   J\dot x - J X_h) - \int
 \omega(X_G^{\tx{h}},   J\dot x - J \widetilde {X_f}) \ \ge \ 0. 
\end{eqnarray*}

The left hand side of the above inequality can be written as
\begin{eqnarray} 
\label{eq:1 for final ps grad} 
E'& = & (1-\alpha^2)\parallel\dot x-X_h-X_G^{\tx{v}}\parallel^2_\omega   
- 2(1-\alpha^2) \langle \dot x-X_h -X_G^{\tx{v}}, 
  \widetilde {X_f} \rangle_\omega
\\
& & + \ \int \Omega(\widetilde {X_f},   J\dot x) - \int \omega
(X_G^{\tx{h}},   J\dot x)  
\label{eq:2 for final ps grad} \\ 
& & + \ (1-\alpha^2) \parallel \widetilde {X_f}^{\tx{h}} \parallel
^2_\beta - \alpha^2 \int \Omega(\widetilde {X_f},   J\widetilde
{X_f}) 
\label{eq:3 for final ps grad} \\
& & - \ \int \Omega(\widetilde {X_f},   JX_h) + \int \omega(X_G^{\tx{h}},
  J \widetilde {X_f}) .
\label{eq:4 for final ps grad}
\end{eqnarray} 

For $\epsilon_0$ and $\rho$ small enough, the expression $E''$
obtained by summing up (\ref{eq:3 for final ps grad}) and (\ref{eq:4
  for final ps grad}) satisfies 
$$E'' \ge (1- 3\alpha^2 - 2m) \parallel \widetilde {X_f}^{\tx{h}}
\parallel ^2_\beta .$$

Again for $\epsilon_0$ and $\rho$ small enough, 
we break the expression obtained by summing up 
(\ref{eq:1 for final ps grad}) and (\ref{eq:2 for final ps grad}) 
as a sum $E''_1+E''_2$ as follows. 

\begin{eqnarray*} 
E''_1 & = & (1-\alpha^2)  \parallel \dot x^{\tx{v}} - X_h - X_G^{\tx{v}}
\parallel ^2_\Omega \\
& &  - \ (2-2\alpha^2) \langle \dot x^{\tx{v}} - X_h -
X_G^{\tx{v}},   \widetilde {X_f}^{\tx{v}} \rangle _\Omega + \int
\Omega(\widetilde {X_f} ^{\tx{v}},   J\dot x ^{\tx{v}})  \\
& = & (1-\alpha^2) \parallel \dot x^{\tx{v}} - X_h - X_G^{\tx{v}} -
\frac {1-2\alpha^2} {2(1-\alpha^2)} \widetilde {X_f} ^{\tx{v}}
\parallel ^2_\Omega \\ 
& & - \ \frac {(1-2\alpha^2)^2} {4(1-\alpha^2)}
\parallel \widetilde {X_f} ^{\tx{v}} \parallel ^2_\Omega + \int
\Omega(X_h,   \widetilde {X_f} ^{\tx{v}}) \\
& \ge & -2\alpha^2 \parallel \widetilde{X_f} ^{\tx{h}} \parallel
^2_\beta .
\end{eqnarray*}
\begin{eqnarray*}
E''_2 & = & (1-\alpha^2) \parallel \dot x ^{\tx{h}} \parallel ^2_\omega
- (2-2\alpha^2) \langle \dot x ^{\tx{h}},   \widetilde {X_f}
^{\tx{h}} \rangle _\omega \\
& & + \ \int \Omega (\widetilde {X_f} ^{\tx{h}},
  J \dot x ^{\tx{h}})  - \int \omega (X_G^{\tx{h}},   J\dot x
^{\tx{h}}) \\ 
& \ge & \frac {2(1-\alpha^2)} 3 \parallel \dot x ^{\tx{h}} \parallel
^2_\beta - (2-2\alpha^2) \langle \dot x ^{\tx{h}},   \widetilde
{X_f} ^{\tx{h}} \rangle _\beta - \int \pi^*\beta (X_G^{\tx{h}}, J \dot
x ^{\tx{h}}) \\
& & - \ (1-\alpha^2) \int \Omega(\dot x ^{\tx{h}},   \widetilde
{X_f}^{\tx{h}}) 
+  \alpha^2 \int \Omega (\widetilde {X_f} ^{\tx{h}},  
\dot x ^{\tx{h}}) - \int \Omega (X_G^{\tx{h}},   J\dot x
^{\tx{h}}) \\
& \ge & \frac {1-\alpha^2} 3 \parallel \dot x^{\tx{h}}
\parallel^2_\beta - (2-2\alpha^2) \langle \dot x ^{\tx{h}},   \widetilde
{X_f} ^{\tx{h}} \rangle _\beta - 6m^2 \parallel \widetilde {X_f}
^{\tx{h}} \parallel ^2_\beta - \frac 1 {24} \parallel \dot x ^{\tx{h}}
\parallel ^2_{\beta} \\ 
& & + \ \frac {1-\alpha^2} 3 \parallel \dot x^{\tx{h}}
\parallel^2_\beta - \frac 1 2 \int |\dot x ^{\tx{h}} | _\beta
|\widetilde {X_f} ^{\tx{h}} | _\beta - \frac m 4 \int  |\dot x
^{\tx{h}} | _\beta |\widetilde {X_f} ^{\tx{h}} | _\beta \\
& \ge & \Big( \frac {1-\alpha^2} 3 - \frac 1 4 - \frac 1 {24} - 
\frac m 4 \Big)
\parallel \dot x ^{\tx{h}} \parallel ^2_\beta - \Big(\frac 1 4 + \frac
m 4 + 6m^2 
\Big) \parallel \widetilde {X_f} ^{\tx{h}} \parallel ^2_\beta \\
& \ge &  - \Big(\frac 1 4 + \frac
m 4 + 6m^2 \Big) \parallel \widetilde {X_f} ^{\tx{h}} \parallel ^2_\beta .
\end{eqnarray*} 

The inequalities involving $E''_2$ 
hold if $(1-\alpha^2)/3 - 1/4 - \frac 1 {24} - m/4 \ge 0$ 
and if $\epsilon_0$ is
small enough so that 
$\parallel \Omega |_{H_0} \parallel _{\infty,\{S\le 4\}} \le 1/4$ and
$|\cdot | ^2_\omega \ge \frac 2 3 |\cdot | ^2_\beta$ on $E$. 
The inequality involving $E''_1$ holds if $\rho$ is small enough 
(determined by $\alpha$). 

We finally obtain 
$$E' \ge \Big(\frac 3 4 - 5\alpha^2 - 2m - \frac m 4 - 6m^2\Big) \parallel
\widetilde {X_f}^{\tx{h}} \parallel ^2_{\beta} \ge 0 .$$
The last inequality holds if $m$ is small enough,
provided $\alpha$ is also small enough.

The fact that equality in (\ref{eq:pseudo-gradient result in
  general}) is attained only if $x$ is a $1$-periodic orbit of $h+G$
in a critical fiber of $\widetilde f$ is obvious from the fact that
all the above inequalities have to be equalities. In particular we
must have $\widetilde {X_f} \circ x \equiv 0$. 
\hfill{$\square$}


\section{Transversality for split almost complex structures}
\label{sec:main section transv}

One crucial ingredient in the construction of the Floer complex associated
to a vector field $Y$ and to an almost
complex structure $J$ is the possibility to choose the pair $(Y,  
J)$ such that the linearized operator 
$$
D_u: W^{1,p}(\RR \times \Ss^1,   u^*T\widehat E) \longrightarrow
L^p(\RR \times \Ss^1,   u^*T\widehat E),
$$ 
$$
\xi \longmapsto \nabla_s \xi + J(u) \nabla_t \xi + \nabla_\xi J(u)
\cdot u_t - \nabla _\xi Y(u)
$$
is surjective for every finite energy solution of the equation 
\begin{equation} \label{eq:Floer intermediate in transv}
u_s + J(u)\cdot u_t = Y(u) .
\end{equation} 
In Floer's original
setting one has $Y = J X_H$, with $H$ a given Hamiltonian, whereas in
our setting $Y$ is a ($s$-independent) vector field satisfying
Definition~\ref{defi:admissible Hamiltonian deformation}. 
 
Our definition of the Floer homology groups makes use
of the fibered structure on $\widehat E$ only in order to prove the a priori
$C^0$-bounds on the finite energy solutions of~(\ref{eq:Floer
  intermediate in transv}). The arguments developed
in~\cite{FH94,FHS,HS95,SZ92} apply in order to show that,
for a fixed 
choice of $Y$, transversality can  be achieved by a generic
choice of $J$ provided one
allows the use of almost complex structures that are time-dependent.
This is sufficient in order to {\it define} Floer homology in the
setting of the present paper. Nevertheless, in order to {\it compute}
it by  
constructing a spectral sequence, one
needs to establish transversality inside the smaller class of 
split almost complex structures
whose horizontal component is time-independent. This requires a
refinement of the above mentioned arguments, by allowing
not only variations of
the vertical almost complex structure, but also of the horizontal
distribution as in~\cite[\S8.2 sq.]{MS04} and~\cite{Seidel97}. The purpose
of this section is to prove this refined version of transversality.

We denote by $V$ or $\tx{Vert}$ the vertical subbundle $\ker \, \pi_*
\subset T \widehat E$, whereas the horizontal
subbundle $V^{\perp_\Omega}$ is denoted by $H_0$ or $\tx{Hor}$.  
We denote by $\mc J^{\tx{vert}}_\tau$ the space of smooth
$\tau$-periodic almost complex structures on $V$ 
which are time-independent and 
standard outside a compact set and which are compatible with $\Omega$ 
(see Definition~\ref{defi:standard almost complex structure on
    symplectic cone}). We denote by $\mc J_B$ the space of smooth {\it
    time-independent} almost complex structures on $B$ which are
  compatible with $\beta$ and which satisfy the {\sc (negativity)}
  property of Definition~\ref{defi:main}. We endow $\mc
  J^{\tx{vert}}_\tau$ and $\mc J_B$ with the
  $C^\infty$-topology. 
Given $J_B \in \mc J_B$ we denote by $\widetilde J_B^H$ its lift with
respect to a given horizontal distribution $H$. 

The connexion $2$-form $\Omega$ can be perturbed while preserving at the
same time closedness and keeping it unchanged along the fibers. 
We describe here a method  borrowed from~\cite[\S8.2]{MS04}, with the
significant difference that we need to allow time-dependent
perturbations of the horizontal distribution in order to achieve
transversality.  

The starting point is to consider on $B$ 
a $1$-form $H$ with values in the bundle $\mc
C^\infty_0(\widehat E)$ whose fiber at $z\in B$ is the space
$C^\infty_0(\widehat E_z)$ of compactly
supported smooth functions on $\widehat E_z$. We assume that $H$ is
$\tau$-periodic and we use from now on the notation $H^t$ in order to
express the dependence on $t$, with $H^t = H^{t+\tau}$. 
We denote the action of $H$ by  
$$
T_z B \longrightarrow C^\infty_0(\Ss^1_{\tau} \times 
\widehat E_z),\quad z\in B,
$$
$$
\zeta \longmapsto H^t_\zeta,
$$
where $\Ss^1_\tau$ is the circle of length $\tau$. 
We define a time-dependent $1$-form 
$\sigma^t_H \in \Omega^1(\widehat E)$ by
$$
\sigma^t_H(x;   v) = H^t_\zeta(x), \qquad  \zeta = \pi_* v ,
$$
where $x\in \widehat E$ and $v\in T_x\widehat E$. The $1$-form
$\sigma^t_H$ vanishes on the fibers by definition. The connection
$2$-form associated to $H^t$ is defined to be 
\begin{equation}
  \label{eq:conn 2 form perturbed}
  \Omega^t_H = \Omega - d\sigma^t_H .
\end{equation}
The horizontal subspace of $\Omega^t_H$ at $x\in \widehat E$ is 
\begin{eqnarray*}
  \tx{Hor}_{H^t;  x} & = & \{ v - X_{H^t_\zeta}(x) \ : \ v\in \tx{Hor}_x,
  \ \zeta = \pi_*v \} \\
   & = & \tx{graph} \big( -X_{H^t}: \tx{Hor}_x \longrightarrow \tx{Vert}_x ,
   \quad v\longmapsto -X_{H^t_{\pi_*v}}(x) \big) .
\end{eqnarray*}
In the above notation $X_{H^t_\zeta}$ represents, for a given $z\in B$, 
the Hamiltonian vector field of the
function $H^t_\zeta$ defined on $\widehat E_z$. 
We denote $\mc H_\tau = \Omega^1(B,  \mc C^\infty_0(\Ss^1_{\tau} \times 
\widehat E))$ and, for
a given compact set $N \subset \widehat E$, we let 
$$
\mc H_\tau(N) \subset \mc H_\tau
$$ 
be the subspace of all those forms with support
contained in $N$. 
The connection $2$-form $\Omega^t$ is 
time-dependent but nevertheless constant on the fibers. For each curve
in $B$ parallel transport along $H^t$ defines a path of
symplectomorphisms between the fibers.

We use the shorthand notation $\widetilde J_B^H$ for the
(time-dependent~!) lift
$\widetilde J_B^{\tx{ Hor}_{H^t}}$ of an almost complex structure $J_B$ on
the base. 
Any triple $(J_V,   J_B,   H) \in \mc{J}^{\tx{vert}}_\tau \times \mc J_B
\times \mc H_\tau$ gives rise to an almost complex structure $J$ on
$\widehat E$ defined as 
$$
J = J_V \oplus \widetilde J_B^H .
$$
Let $J_0$ be the almost complex
structure corresponding to the fixed triple $(J_V,   J_B,   0)$. The
action of the almost complex structures corresponding to triples
$(J_V,   J_B,   H)$, $H\in \mc H_\tau$ can be explicitly described
(see~\cite{MS04}) as 
\begin{equation} \label{eq:action of a c str under perturb of hor distrib}
J_tv = J_0v + J_{V,t} X_{H^t_{\pi_*v}}(x) - X_{H^t_{J_{_B}\pi_*v}}(x),
\quad v\in 
T_x\widehat E .
\end{equation}

  Before stating our transversality result, 
  we recall the following theorem of Salamon and Zehnder, which
  implies in particular transversality for all moduli spaces of Floer
  trajectories in the time-independent setting for symplectically
  aspherical manifolds.  

  \begin{thm}[\cite{SZ92}, Thm. 7.3] \label{thm:transversality on the base} 
    Let $(B,   \beta)$ be a closed symplectic manifold such that 
$$
\langle [\beta],   \pi_2(B) \rangle =0 . 
$$
  Let $f:B\longrightarrow \RR$ be a Morse function and $J_B$ a
  time-independent almost complex structure compatible with $\beta$,
  such that the flow of $\nabla ^{J_{B}} f$ is Morse-Smale. Let 
  $$
   D_{u,\tau} : W^{1,p}(\RR \times \Ss^1_\tau,   u^*TB)
   \longrightarrow L^p(\RR \times \Ss^1_\tau,   u^*TB),
  $$
  $$
   \xi \longmapsto \nabla_s\xi + J_B(u) \nabla_t\xi + \nabla_\xi J_B
   (u) \cdot u_t - \nabla_\xi \nabla f
  $$
  be the linearization of Floer's equation 
  \begin{equation} \label{eq:Floer Morse eq for transv}
  u_s + J_B(u) u_t = (\nabla ^{J_B}f) \circ u,
  \end{equation} 
  defined for $\tau$-periodic
  maps $u:\RR \times \Ss^1_\tau \longrightarrow B$. 
  The following assertions hold if $\tau$ is small enough. 
  \renewcommand{\theenumi}{\alph{enumi}}
  \begin{enumerate}
  \item The operator $D_{u,\tau}$ is surjective for any
    solution $u:\RR \longrightarrow B$
    of~(\ref{eq:Floer Morse eq for transv}) which is independent of
    $t$. 
  \item Every finite energy solution of~(\ref{eq:Floer Morse eq for
      transv}) is independent of $t$. 
  \end{enumerate}
  \end{thm}

\begin{rmk} {\bf (Reparametrizations)} \rm 
(i) The norm of 
$D_{u,\tau}$ does not depend on the parameter
$\tau\in]0,\tau_0]$ because $D_{u,\tau}\xi=\nabla _s \xi
-\nabla_\xi\nabla f$.
 
(ii) The statement of the above theorem
remains true if we fix the period and allow the coefficient in front
of $f$ to go to zero. The reason is that any $\tau$-periodic
solution $u(s,t)$ gives rise to a $\tau_0$-periodic solution
$u_0(s,t)=u(\frac 
\tau {\tau_0}s,\frac \tau {\tau_0} t)$, which in turn satisfies the
equation $\p_s u_0 + J_B \p_t u_0 = \frac \tau {\tau_0} \nabla f$.
\end{rmk}

We fix from now on a Morse function $f:B \longrightarrow \RR$, an 
almost complex structure $J_B \in \mc J_B$ and a
period $\tau > 0$ such that the conclusions of
Theorem~\ref{thm:transversality on the base} hold true.

If the almost complex structure $J_B$ satisfies the {\sc (negativity)}
assumption of Definition~\ref{defi:main}, then for every $(J_V,  
H)\in \mc J^{\tx{vert}} _\tau \times \mc
H_\tau$   
there exists $\epsilon_0 > 0$ 
such that 
$$
\omega_\epsilon = \pi^*\beta + \epsilon \Omega
$$ 
tames 
$$
J=J_V\oplus \widetilde J_B^H
$$ 
for all $0 < \epsilon \le \epsilon _0$. 
Following the previous section, we consider  Hamiltonians of the form 
$$
K=h+\widetilde f + G .
$$
The function $h=h(S)$ is convex and linear for $S \ge 1$ with slope
$\lambda_{\max}$ satisfying 
$$
\frac \tau {\epsilon_0}\cdot \lambda_{\tx{max}} \notin
\tx{Spec}(\partial E). 
$$ 
The function $G$ is a $\tau$-periodic 
perturbation localized in a neighbourhood of the
critical fibers of $\widetilde f = \pi \circ f$. More precisely, we
denote the critical points of $f$ by
$p_i$, $1\le i \le \ell$, we fix open
neighbourhoods $p_i \in V_i \Subset U_i$ such that $U_i \cap U_j =
\emptyset$, $i\neq j$ and we require that 
$\tx{supp}(G) \subset E \cap \bigcup_{p_i \in \tx{Crit}(f)}
\pi^{-1}(V_i)$. We denote 
$$
N = E  \ \cap \ ^c\pi^{-1}\big(\bigcup U_i\big),
$$
so that $N \cap \tx{supp}(G) = \emptyset$. 
Given $h$ we choose $G$ and $\epsilon_0$ 
such that the $\tau$-periodic
orbits of $K$ with respect to $\omega_{\epsilon_0}$ 
are nondegenerate and lie in the critical fibers of 
$\widetilde f = \pi \circ f$, while the vector field  
$$
\mc Y^{\epsilon_0} (t,  x) = J_t \dot x - Y^{\epsilon_0}(t,  x),
$$
$$
Y^{\epsilon_0} (t,   x)= J_t\big(X^{\epsilon_0}_{\epsilon_0 h} +
\epsilon_0 \widetilde {X_f}  +
(X^{\epsilon_0}_{\epsilon_0 G})^{\tx{v}}(t,  x)\big)
$$
is a negative pseudo-gradient for the action functional
$A_{\epsilon_0 K}^{\epsilon_0}$ 
defined on the space of contractible $\tau$-periodic loops in
$\widehat E$. The superscript $\epsilon_0$ indicates, as usual, the
fact that the Hamiltonian vector fields and the symplectic action  
are computed with respect to the form $\omega_{\epsilon_0}$.

\begin{rmk} {\bf (Uniform upper bound for $\epsilon$)} 
\label{rmk:choice of epsilon 0}
\rm The rescaling parameter $\epsilon$ is allowed to vary in some
  interval $]0,\epsilon_0]$ with $\epsilon_0$ small enough in order to
  ensure nondegeneracy and taming for $\omega_\epsilon=\pi^*\beta
  +\epsilon \Omega$.  
The only other point where  
we use the rescaling $\Omega \rightsquigarrow \epsilon \Omega$ 
is in the proof of the 
pseudo-gradient property for loops contained in $E$ in
Proposition~\ref{prop:general ps gr}, where one might 
need to further diminish the constant
$\epsilon_0$, depending on  $\| G|_E \|_{C^1}$. Since the latter
  quantity can be uniformly bounded independently of the choice of $h$,
  we conclude that we can construct an admissible cofinal family of
  Hamiltonians $K$ admitting a uniform constant $\epsilon_0$.
\end{rmk}

Given a Hamiltonian $K$ and a parameter $0 < \epsilon \le \epsilon_0$ as in
Remark~\ref{rmk:choice of epsilon 0}  
we define the space of \emph{regular} vertical almost complex
structures and Hamiltonian 
perturbations
$$
(\mc J^{\tx{vert}}_\tau \times \mc H_\tau )^{\tx{reg}}(\epsilon,   K)
\subset \mc J^{\tx{vert}}_{\tau} 
\times \mc H_{\tau}(N)
$$
as consisting of pairs $(J_V,   H) \in \mc
J^{\tx{vert}}_{\tau} 
\times \mc H_{\tau}
(N)$ such that, for  
every finite energy solution $u:\RR \times \Ss^1_{\tau}
\longrightarrow 
\widehat E$ of the equation 
\begin{equation} \label{eq:Floer eq upstairs} 
u_s + J u_t = Y^\epsilon \circ u,   
\end{equation} 
the linearized operator 
$$
D_{u,\tau}^\epsilon : W^{1,p}(\RR \times
\Ss^1_{\tau} ,   u^*T\widehat E) \longrightarrow 
L^p(\RR \times \Ss^1_{\tau} ,  u^* T^* \widehat E),
$$
$$
\xi \longmapsto \nabla_s^\epsilon \xi + J(u) \nabla_t^\epsilon \xi +
\nabla_\xi^\epsilon J(u) 
\cdot u_t - \nabla _\xi^\epsilon Y^\epsilon(u)
$$
is surjective. 
The connexion $\nabla^\epsilon$  
is the Levi-Civita connexion associated to the metric defined by
$\omega_\epsilon$ and $J=J_V \oplus \widetilde J_B^H$.

\medskip 

\noindent {\bf Notation.} 
Given a map $u: \RR \times \Ss^1 _{\tau} \longrightarrow
\widehat E$ we denote its projection by 
$$
v =\pi \circ u .
$$
Given $J_V\in\mc J^{\tx{vert}}_\tau$ we denote by $\mc
J^{\tx{vert}}_\tau(J_V)$ the space of vertical almost complex
structures which coincide with $J_V$ outside a compact set.

\begin{rmk} {\bf (Geometric property of the pseudo-gradient equation)}
\rm The fundamental property of equation~(\ref{eq:Floer eq upstairs}) is
that the projected solutions $v=\pi \circ u$ satisfy the equation
$$
v_s + J_B (v)v_t = \epsilon \nabla ^{J_B} f (v) ,
$$
for which transversality is ensured by Theorem~\ref{thm:transversality
  on the base}. This geometric property plays a crucial role not only
in the construction of the spectral sequence, but also in the proof of 
transversality within the class of split almost complex structures of
the type $J_V\oplus \widetilde J_B^H$. 
\end{rmk} 

The aim of this section is to prove the following result. 

\begin{thm}[split transversality]
\label{thm:transversality for split almost complex structures} 
 Let $K$, $\epsilon_0$ be as above and $J_{V,0} \in \mc J^{\tx{vert}}_\tau$. 
\renewcommand{\theenumi}{\alph{enumi}}
\begin{enumerate}
  \item There exist a positive constant
    $\epsilon(K)\in]0,\epsilon_0]$ and an open neighbourhood
    $\mc O$ of $(J_{V,0},   H_0)$ in $\mc
    J^{\tx{vert}}_{\tau}(J_{V,0}) \times \mc H_{\tau}(N)$ such that
    the operator   
$$
F_{u,\tau}^\epsilon : W^{1,p}(\RR \times
\Ss^1_{\tau},   v^*TB) 
\longrightarrow L^p(\RR \times \Ss^1_{\tau},   v^*TB),
$$
$$
F_{u,\tau}^\epsilon (\xi) = \pi_* \big(
D_{u,\tau}^\epsilon \cdot 
\widetilde \xi \ \big)
$$ 
is surjective for any finite energy solution $u$ of~(\ref{eq:Floer eq
  upstairs}) with $0 < \epsilon \le \epsilon(K)$ and  
$(J_V,   H) \in \mc O$;
 \item For every $0< \epsilon \le
    \epsilon(K)$ the set 
   $(\mc J^{\tx{vert}}_\tau \times \mc H_\tau)^{\tx{reg}}(\epsilon,K) \ \cap \
    \mc O$ is dense and of second Baire category in $\mc O$.  
 \end{enumerate}
\end{thm}

\begin{rmk} {\bf (On the parameter $\epsilon(K)$)} 
\label{rmk:why epsilon nu}
\rm The parameter $\epsilon(K)$ depends actually on the asymptotic
    slope of $K$. Indeed, it will be appearant from the proof that
    $\epsilon(K)$ depends on the $C^0$-bound for solutions of
    $u_s+Ju_t=Y^\epsilon\circ u$ through the use of
    Lemma~\ref{lem:transv up through trans down} below. These bounds,
    in turn, depend via Lemma~\ref{borne norme delta} on the maximal
    difference between the actions of two closed orbits, i.e. on the
    asymptotic slope of $K$. 
\end{rmk}

\demo a) Let $\mc O$ be a neighbourhood of $(J_{V,0},H_0)$ such that
the following hold for $(J_V,H)\in\mc O$: 
\begin{itemize} 
\item the almost complex structures $J=J_V\oplus \widetilde J_B^H$ are
  tamed by $\omega_\epsilon$, $0< \epsilon \le \epsilon_0$;
\item finite energy solutions of $u_s+Ju_t=Y^\epsilon\circ u$ admit a
    common uniform $C^0$-bound for $0 < \epsilon \le
    \epsilon_0$. This can be achieved as a consequence of the following
    two observations which slightly generalize the proof of
    Theorem~\ref{estimation C 0}. Firstly, it is clear that, for
    $\epsilon>0$ fixed, one 
    can allow the almost complex structure to slightly vary inside a
    compact set. Secondly, as $\epsilon>0$ varies the vertical part of
    the vector field $Y^\epsilon$ remains unchanged, whereas the
    horizontal part is rescaled by $\epsilon$, so that
    assumptions~(\ref{eq:horizontal bound on Y}-\ref{controle derivee
    seconde}) in Definition~\ref{defi:admissible Hamiltonian
    deformation} are still satisfied.
\end{itemize}

\begin{lem} \label{lem:transv up through trans down}
  There exist constants $c(\epsilon)$, $0 < \epsilon \le \epsilon_0$
  with $c(\epsilon)\longrightarrow 0$, $\epsilon \longrightarrow 0$
  such that, for any finite energy solution $u$ of
  equation~(\ref{eq:Floer eq 
    upstairs}) with $(J_V,   H) \in \mc O$, we have
  \begin{equation*}
    ||| F_{u,\tau}^\epsilon - D_v ||| \le c(\epsilon) .
  \end{equation*} 
\end{lem}

\demo The Levi-Civita connection $\nabla$ associated to a metric
$\langle \cdot, \cdot \rangle$  
can be expressed as follows:
\begin{eqnarray}
\langle \nabla_Y X ,   Z \rangle & = & \frac 1 2 \big\{ X\langle Y,
  Z \rangle + Y \langle Z,   X \rangle - Z \langle X,   Y \rangle
\\ 
& & - \langle [X,   Z],   Y \rangle - \langle [Y,   Z],   X
\rangle - \langle [X,   Y],   Z \rangle \big\} . 
\end{eqnarray}  
By applying the above formula to the connection $\nabla^\epsilon$
associated to the (time-dependent) metric $g_\epsilon (v,   w) = \frac 1 2
(\omega_\epsilon(v,   Jw) + \omega_\epsilon(w,   Jv))$ one sees
that, for any two 
vector fields $X\in \mc X (\widehat E)$ and $Y\in \mc X(B)$, we have 
\begin{eqnarray*}
&  (\nabla^\epsilon_{\widetilde Y} X )^{\tx{h}} \longrightarrow
\widetilde{\nabla _Y {\pi_* X} }, & \\  
& (\nabla^\epsilon_X \widetilde Y )^{\tx{h}} \longrightarrow
\widetilde{\nabla _{\pi_* X} Y}, & \epsilon \longrightarrow 0, 
\end{eqnarray*}
where $\nabla$ is the Levi-Civita connection on $B$ corresponding to
the metric $g_{_B}(\cdot,   \cdot) = \beta(\cdot,  J_B \cdot)$. The
convergence is uniform on every compact set. More precisely, for any
compact set $\mc K \subset \widehat E$ we have   
\begin{equation} 
\label{eq:estimate on cov deriv of vectors}
|| (\nabla^\epsilon_{\widetilde Y} X )^{\tx{h}} -  \widetilde{\nabla
  _Y {\pi_* X}} ||_{g_{_B}} \le c_1(\epsilon,   \mc K) || X ||
_{C^1(g_{\epsilon_0})} ||Y || _{C^0}, 
\end{equation} 
\begin{equation}
\label{eq:estimate on cov deriv of vectors 2}
|| (\nabla^\epsilon_X \widetilde Y )^{\tx{h}} -  \widetilde{\nabla
  _{\pi_* X} Y} ||_{g_{_B}} \le c_1(\epsilon,   \mc K) || X ||
_{C^0(g_{\epsilon_0})} ||Y || _{C^1}, 
\end{equation}
with $c_1(\epsilon,   \mc K) \longrightarrow 0$, $\epsilon
\longrightarrow 0$.  
Similarly we have 
\begin{equation} \label{eq:estimate on cov deriv of J}
|| (\nabla^\epsilon _{\widetilde Y} J)^{\tx{h}} - \widetilde
{\nabla _Y J_B} || _{g_{_B}} \le c_2(\epsilon,  \mc K) || Y ||
_{C^0(g_{_B})},  
\end{equation}
with $c_2(\epsilon,   \mc K) \longrightarrow 0$, $\epsilon
\longrightarrow 0$. The estimates (\ref{eq:estimate on cov deriv of
  vectors} -- \ref{eq:estimate on cov deriv of J}), together with the
explicit form of the operators involved and the existence of the
uniform $C^0$-bound on solutions of Floer's equation~(\ref{eq:Floer eq
  upstairs}) for $(J_V,   H) \in \mc O$ and $0<
\epsilon \le \epsilon_0$,  
 imply the conclusion.
\hfill{$\square$}

\begin{lem} \label{lem:uniform upper bound} 
The operators $D_v$, where $v$ runs over all Floer
(time-independent) trajectories on $B$ corresponding to $f$ and $J_B$,
admit uniformly bounded right inverses.
\end{lem}

\demo
This is a reformulation of
the gluing theorem for Floer trajectories in the transverse case (see
for example~\cite{Salamon-lectures} for the latter). The key step in
the gluing construction is to prove that the linearized operator is
surjective along preglued curves and that it admits a right inverse
which is uniformly bounded for large enough values of the gluing
parameter. This implies that one can find uniformly bounded right
inverses for the operator $D_v$ when $v$ belongs to a small
neighbourhood of the boundary of the moduli space of 
trajectories. Such a uniform bound can clearly be found on the
remaining relatively compact set contained in the interior of the
moduli space of trajectories, and gives the existence of
uniformly bounded right inverses for $D_v$ for any choice of $v$.  
\hfill{$\square$}

\medskip 

We prove now assertion a) in the theorem. Let $C$ be the
  uniform upper bound provided by Lemma~\ref{lem:uniform upper bound}
  and choose $\epsilon(K)$ small enough so that the constant
  $c(\epsilon)$ in Lemma~\ref{lem:transv up through trans down}
  satisfies $c(\epsilon)<1/2C$ for $0<\epsilon\le \epsilon(K)$. Given 
  a finite energy solution $u$, let $Q_v$ be a right inverse for
  $D_v$, $v=\pi \circ u$ such that $\|Q_v\|\le C$. Then
  $$
\|F^\epsilon_{u,\tau}Q_v - \textrm{Id}\| =
\|F^\epsilon_{u,\tau}Q_v - D_vQ_v\| \le 1/2,
  $$
hence the operator $F^\epsilon_{u,\tau}Q_v$ is invertible and
the norm of its inverse is $\le 2$. Then 
$Q_v(F^\epsilon_{u,\tau}Q_v)^{-1}$ is a right inverse for 
$F^\epsilon_{u,\tau}$ of norm $\le 2C$, and in particular 
$F^\epsilon_{u,\tau}$ is surjective.

\medskip 

We prove now assertion b). We follow the proof of Theorem 5.1
in~\cite{FHS} -- in particular, it is enough to obtain the conclusion
when  $\mc J^{\tx{vert}}_{\tau}(J_{V,0}) \times \mc H_{\tau}(N)$ is
endowed with the $C^\ell$-topology, $\ell \ge 1$.

We fix
$0<\epsilon \le \epsilon(K)$. The key step is to
prove that, for any $\tau$-periodic orbits $x^-,   x^+$ of $K$, 
the universal moduli space 
$$ 
\mc M(x^-,   x^+,   \mc O) = \{ (u,   J_V,   H) \ :
\ \bar \partial _{J_V \oplus \widetilde J_B^H,   Y^\epsilon} \ u =
0\} 
$$
is a Banach manifold. The universal moduli space is naturally the zero
set of the section  
$$
\mc F: \mc B \times \mc O \longrightarrow \mc E, \qquad
\mc F(u,   J_V,   H) = \bar \partial _{J_V \oplus \widetilde J_B^H,
    Y^\epsilon} \ u .  
$$
Here $\mc B = \mc B(x^-,   x^+)$ is the space of continuous maps $u:
\RR \times \Ss^1_{\tau} \longrightarrow \widehat E$ which are
locally of class $W^{1,p}$ and which converge to $x^-$, $x^+$ as
$s\longrightarrow \pm \infty$ with a suitable exponential decay
condition at infinity~\cite{FHS}, while $\mc E$ is the Banach
bundle whose fiber at $(u,   J_V,   H)$ is $L^p(u^*T \widehat E)$.  

We need to prove that $\mc F$ is transverse to the zero section of
$\mc E$. We denote by $\pi$ the vertical projection $T_{((u,   J_V,
    H),  0)}\mc E \longrightarrow \mc E_{(u,   J_V,   H)}$. The
vertical differential  
$$
D\mc F(u,   J_V,   H) = \pi \circ d\mc F(u,   J_V,   H)
$$ 
is given at a solution $u$ of~(\ref{eq:Floer eq
  upstairs}) by 
\begin{equation} \label{eq:vert diff for split a c str}
D\mc F(u,   J_V,   H) \cdot (\xi,   Z,   h) = D_u \xi + Z_t(u) J_V
u_s^{\tx{vert}} + X_{h_{\pi_* u_s^{\tx{horiz}} }} . 
\end{equation}
Here $\xi \in T_u \mc B = W^{1,p}(u^*T\widehat E)$ and $h\in T_H \mc
H_{\tau}(N) = \mc H_{\tau}(N)$. 
The tangent vector $Z\in T_{J_V}\mc J
^{\tx{vert}}_{\tau}(J_{V,0})$  is a $C^\ell$-map
$\Ss^1_{\tau} \times T\widehat E \longrightarrow T\widehat E$ which
has compact support in $\widehat E$ and satisfies   
$$
J_V^t Z_t + Z_t J_V^t = 0, \qquad \Omega(Z_t v,   w) + \Omega (v,  
Z_t w) = 0, \ v, w \in \tx{Vert} .  
$$

Let us explain the term $X_{h_{\pi_* \partial _s u}}$ in~(\ref{eq:vert
  diff for split a c str}). We need to study the change in  
$$
\bar \partial _{J_V \oplus \widetilde J_B^H,   Y^\epsilon} \ u = u_s
+ (J_V \oplus \widetilde J_B^H) u_t - (J_V \oplus \widetilde J_B^H)
(X^{\tx{v}} + \epsilon {\widetilde {X_f}}^H) 
$$
as we replace $H$ by $H+h$. By~(\ref{eq:action of a c str under
  perturb of hor distrib}) we obtain 
\begin{eqnarray*}
 (J_V \oplus \widetilde J_B^{H+h}) \cdot \widetilde {X_f} ^{H+h} & = &
 (J_V \oplus \widetilde J_B^{H+h}) \cdot (\widetilde {X_f}^H -
 X_{h_{X_f}}) \\ 
 & = & -J_V X_{h_{X_f}} +  (J_V \oplus \widetilde J_B^H)\cdot
 \widetilde {X_f}^H + J_V   X_{h_{X_f}} -  X_{h_{J_B X_f}} \\ 
 & = & (J_V \oplus \widetilde J_B^H)\cdot \widetilde {X_f}^H - X_{h_{J_B X_f}} 
\end{eqnarray*} 
and 
\begin{eqnarray*} 
 (J_V \oplus \widetilde J_B^{H+h}) \cdot u_t & = & (J_V \oplus
 \widetilde J_B^H) \cdot u_t + J_V X_{h_{\pi_* u_t}} - X_{h_{J_B \pi_*
     u_t}} . 
\end{eqnarray*} 

Since $J_B\pi_*u_t - \epsilon J_B X_f = \pi_* u_s$ we infer that the
variation in the direction $h$ of $\bar \partial 
_{J_V \oplus \widetilde J_B^H,   Y^\epsilon}  u$ is  
$$ 
\delta \bar \partial \cdot h = J_V X_{h_{\pi_* u_t}} + X_{h_{\pi_* u_s}} .
$$
On the other hand $\pi_* u_t$ vanishes because the projected
trajectories on the base are time-independent, hence $\delta \bar
\partial \cdot h =  X_{h_{\pi_* u_s}}$.

We need to show that $D\mc F$ is onto for any $(u,   J_V,   H) \in
\mc M (x^-,   x^+,   \mc O)$. The operator $D_u$ is
Fredholm and the same holds for $D\mc F$. In particular $\tx{im} \, (D\mc
F)$ is closed and, in order to prove surjectivity, it is enough to
prove that $D\mc F$ has a dense range. Equivalently, we have to show
that the annihilator  
$$
A = \{ \eta \in L^q(u^*T\widehat E) \ : \ \iint_{\RR \times
  \Ss^1_{\tau}} \langle \eta,   D\mc F (u,   J_V,   H)
\cdot (\xi,   Z,   h) \rangle   ds   dt = 0, \quad \forall \ \xi,
  Z,   h \} 
$$
is zero, where $1/p + 1/q =1$. For $\eta \in A$ we must have
\begin{equation} \label{eq:first annihilator condition} 
\iint \langle \eta,   D_u\xi \rangle   ds   dt  = 0, 
\end{equation} 
\begin{equation} \label{eq:second annihilator condition}
\iint \langle \eta,   Z_t(u) J_V u_s^{\tx{vert}} \rangle   ds   dt  =  0,
\end{equation}
\begin{equation} \label{eq:third annihilator condition}
\iint \langle \eta,   X_{h_{\pi_* u_s^{\tx{horiz}}}} \rangle   ds   dt  =  0 
\end{equation}
for all $\xi \in W^{1,p}(u^*T\widehat E)$, $Z\in T_{J_V} \mc
J^{\tx{vert}}_{\tau}(J_{V,0})$, $h \in \mc H_\tau (N)$.  
Condition~(\ref{eq:first annihilator condition}) states that $\eta$ is
a weak solution of $D_u^*\eta=0$, where $D_u^*$ is the formal adjoint
of $D_u$ which is obtained by formally replacing in $D_u$ the term
$\nabla_s$ with $-\nabla_s$.  Elliptic regularity implies that $\eta$
is of class $C^\ell$ and is a strong solution of $D_u^*\eta=0$, while
unique continuation (see~\cite{FHS}) ensures that it is enough to show
that $\eta$ vanishes on an open set in order to obtain global
vanishing. We denote $\eta=
\eta^{\tx{h}} + \eta ^ {\tx{v}}$ the decomposition of $\eta$
into horizontal and vertical parts with respect to the splitting
$T\widehat E = \tx{Vert} \oplus \tx{Vert}^{\perp_\Omega}$. In
general, the (time-dependent) 
decomposition of a vector $X$ with respect to the
splitting $T\widehat E = \tx{Vert} \oplus \tx{Hor}_H$ 
will be denoted $X= X^{\tx{vert}} +
X^{\tx{horiz}}$.

We first show that $\eta^{\tx{v}} \equiv 0$ on some open set $U$. 
The special form of
our Floer equation separates nonconstant trajectories $u$ in two
classes: those entirely contained in a fiber, and those satisfying
$\pi_*u_s=v_s \neq 0$ running from one fiber to another.

Let $u$ be a nonconstant trajectory contained in a fiber. One of the
fundamental results of~\cite{FHS} states that the set 
$$
R(u) = \big\{ (s,   t) \in \RR \times \Ss^1_{\tau} \ : \ u_s(s,
  t) \neq 0, \ u(s,   t) \neq x^{\pm}(t), \ u(s,   t) \notin u(\RR
- \{ s \} ,   t) \big\}  
$$
of {\it regular} points is open and dense in $\RR \times
\Ss^1_{\tau}$. We claim that $\eta^{\tx{v}}\equiv 0$ on
$R(u)$, hence by density $\eta^{\tx{v}}\equiv 0$. Assume by
contradiction that this is not the case and $\eta^{\tx{v}}(s,  t)
\neq 0$ for some $(s,  t) \in R(u)$, hence on a small neighbourhood
of $(s,  t)$. Because $u_s^{\tx{v}}= u_s = u_s^{\tx{vert}}$ 
we can choose, according to~\cite{FHS},
a time-dependent tangent vector $Z_t$ such that $\iint_{\RR
\times \Ss^1_{\tau}} \langle \eta,  
Z_t(u)J_Vu_s^{\tx{vert}} \rangle >0$, contradicting~(\ref{eq:second
  annihilator condition}).

Let now $u$ be a nonconstant trajectory running from one fiber to
another. In this case we have $R(u)=\RR \times
\Ss^1_{\tau}$. 
We distinguish two situations: 
either $u$ crosses $N$, or not. 
In the second case the part of $u$ projecting onto $B\setminus
\bigcup_i U_i$ must live in $\{S\ge 1\}$, and we deduce the existence
of a point $(s,  t)$
such that $u(s,  t) \in E \setminus N$ and $u_s^{\tx{vert}}(s,  t) =
u_s^{\tx{v}}(s,   t) \neq 0$. Then we can use again
condition~(\ref{eq:second annihilator condition}) in order to 
show that $\eta^{\tx{v}}$ has to vanish in a neighbourhood $U$ of
$(s,  t)$. In 
the first case we use~(\ref{eq:third annihilator condition}) in order
to show that $\eta^{\tx{v}}$ vanishes on $U=u^{-1}(\tx{im} \,  u \cap
\tx{int}  N)$. By contradiction, let $(s,   t)$ be a point where
$\eta^{\tx{v}}(s,   t) \neq 0$. We know that $\pi_*u_s = v_s \neq 0$,
hence we can choose a time-dependent
tangent vector $h$ with support in $N$ such that $\iint _{\RR \times
  \Ss^1_{\tau}} \langle \eta,   X_{h_{\pi_*u_s}} \rangle
>0$, contradicting~(\ref{eq:third annihilator condition}).

Let us now prove that $\eta^{\tx{h}}$ also vanishes on some nonempty open
set $V\subset U$, knowing that $\eta^{\tx{v}}$ vanishes on $U$. Let $\beta:\RR
\times \Ss^1_{\tau} \longrightarrow \RR_+$ be a positive 
smooth function supported in $U$ which is not identically zero. We
claim that $\eta^{\tx{h}}$ vanishes on the open set $\{z\in U \ : \
\beta(z) > 0 \}$. By contradiction, let us assume that this is not the case.   
The projection $\pi_*(\beta\eta)  = \pi_*(\beta\eta^{\tx{h}})$ is an $(s, 
t)$-dependent vector field along $v = \pi_*u$, supported in $U$ (note that $v$
may as well be constant, but this does not interfere with the argument). We
have seen that the operator $F^\epsilon_{u,   \tau} = \pi_*
\circ D_u \circ \widetilde{\quad}$ is surjective, hence there exists
$\xi \in W^{1,p}(v^*TB)$ such that $(D_u\cdot \widetilde
\xi)^{\tx{h}} = \beta \eta^{\tx{h}}$. We now use $\eta^{\tx{v}}\equiv
0$ on $U$ in order to obtain
$$\iint_{\RR \times \Ss^1_{\tau}} 
\langle \eta,   D_u \widetilde \xi \rangle = \iint_{\RR \times
  \Ss^1_{\tau}}  \langle \eta ,   \beta\eta \rangle = \iint
_U \langle \eta^{\tx{h}},   \beta
\eta^{\tx{h}} \rangle >0,  $$  
which is a contradiction with~(\ref{eq:first annihilator condition}).

We have therefore proved that $\eta$ vanishes on a nonempty open set
$V\subset\RR \times \Ss^1_{\tau}$. By unique continuation
$\eta$ vanishes identically on $\RR \times \Ss^1_{\tau}$ 
and this finishes
the proof of Theorem~\ref{thm:transversality for
  split almost complex structures}.
\hfill{$\square$}


\section{The spectral sequence}
\label{sec:construction of the sp seq}


\subsection{General formalism} \label{sec:general formalism sp seq}

We recall some relevant notions concerning spectral sequences in order
to fix notation, and we refer to~\cite{BT,McC} for details. For our
purposes a spectral sequence is a sequence of 
bigraded differential modules $(E_r^{p,q},d_r)$ associated to a
graded differential complex $(C=\oplus_{k\ge 0} C^k,\partial)$ endowed
with a filtration  
$$ 
C^k= F_0 C^k \supset F_1 C^k \supset \ldots \supset F_n C^k \supset 0. 
$$
Saying that $F_pC$, $p\ge 0$ defines a filtration means that $\partial
(F_p C^k) \subset F_p C^{k+1}$, where we assume that the differential
$\partial$ has degree $+1$. The differential $d_r$ has bidegree
$(r,-r+1)$. The main feature of a spectral sequence is that
$E_{r+1}=H(E_r,d_r)$, $r\ge 0$, with
$E_0^{p,q}=F_pC^{p+q}/F_{p+1}C^{p+q}$. Within the above
setup the groups $E_r^{p,q}$ stabilize and the limit $E_\infty^{p,q}$
satisfies $E_\infty^{p,q}=F_pH^{p+q}/F_{p+1}H^{p+q}$ for some
filtration $F_pH^k$, $p\ge 0$ on the cohomology
$H=H(C,\partial)$. Spectral sequences are functorial in the sense that
a morphism of filtered complexes induces a morphism between the
associated spectral sequences. 

\begin{ex} \label{example:spectral sequence} \rm
  We assume $C^k=\oplus _{p=0}^n C^k_p$ and $\partial = \partial _0 +
  \partial _1 +\ldots + \partial _n$, with $\partial _r : C^k_p
  \longrightarrow C^{k+1}_{p+r}$. We denote $C_p=\oplus _{k\ge
  0}C^k_p$ and $F_p C^k = \oplus _{s\ge p} C^k_s$, so that $F_pC^k$
  defines a filtration. We then have $E_0^{p,q}=C_p^{p+q}$ and
   $$E_1 = \oplus_p H(C_p,  \partial _0),$$
   $$d_1([\alpha_p]) = [\partial _1 \alpha_p] \in H(C_{p+1}, 
   \partial _0), \qquad [\alpha_p] \in H(C_p,  \partial _0).$$
  More generally, let $D_r=\partial _0 + \partial _1 +\ldots +
  \partial _r$, $r\ge 0$. An element $\alpha_p \in C_p$ defines a class in
  $E_r$ if and only if there exist $\beta_{p+i}\in C_{p+i}$, $1\le i
  \le r-1$ such that $D_{r-1}\alpha_p + D_{r-2}\beta_{p+1} + \ldots +
  D_0 \beta_{p+r-1}=0$, and the differential $d_r$ acts as 
  \begin{eqnarray*} 
   d_r[\alpha_p] & = & [\partial _r\alpha_p + \partial _{r-1} \beta_{p+1}
   + \ldots + \partial _1 \beta_{p+r-1}] \\
   & = & [D _r\alpha_p + D _{r-1} \beta_{p+1}
   + \ldots + D _1 \beta_{p+r-1}].
  \end{eqnarray*} 
\end{ex}


\subsection{Morse homology and local systems of coefficients} 
\label{sec:local systems}

\subsubsection{Morse homology with values in a local system of
  coefficients} \label{sec:MHlocsyst}

The formalism of local cofficients for singular or cellular homology 
was introduced by Steenrod~\cite{Steenrod}. This is our main
reference for this section, together with McCleary~\cite{McC}. The
adaptation to Morse homology is straightforward but, to our knowledge,
has not appeared previously in the literature. 

\begin{defi}
  Let $R$ be a ring and $M$ be an $R$-module. A \emph{local system of
    coefficients with fiber $M$} on a topological space $B$ consists
  of the following data: 
  \renewcommand{\theenumi}{\arabic{enumi}}
  \begin{enumerate}
  \item one copy $M_x$ of $M$ for each $x\in B$, called \emph{the
      fiber at $x$};
  \item a family of isomorphisms $\big(\Phi_\alpha:M_{\alpha(0)}
    \stackrel \sim \longrightarrow M_{\alpha(1)}\big)_{\alpha \in \mc
      P(B)}$, where $\mc P(B)$ is the set of continuous paths in $B$,
    such that: 
    \begin{itemize}
    \item if $\alpha \simeq \beta$ are homotopic with fixed endpoints,
      then $\Phi_\alpha = \Phi_\beta$;
    \item if $\alpha,\beta \in \mc P(B)$ satisfy $\alpha(1)=\beta(0)$,
      then $\Phi_{\alpha \cdot \beta} = \Phi_\alpha \circ \Phi_\beta$.
    \end{itemize}
   We call $\Phi_\alpha$ \emph{the parallel transport along $\alpha$}. 
  \end{enumerate}
\end{defi}

\noindent {\bf Remarks.} 1. Isomorphism classes of local systems of
coefficients on a manifold are in one to one correspondence with
isomorphism classes of locally constant sheaves. 

2. If $B$ is simply-connected then all local systems having the same
   fiber are isomorphic. More generally, the choice of a basepoint
   $x_0 \in B$ and of a collection of paths connecting $x_0$ to $x\in
   B$, $x\neq x_0$ determines a one to one correspondence between
   local systems with fiber $M$ and representations $\pi_1(B,   x_0)
   \longrightarrow \tx{Aut} (M_{x_0})$. We call the representation
   associated to a local system $\mc S$ the \emph{monodromy
   representation of $\mc S$}. It is well defined up to conjugation by
   an element of $\tx{Aut} (M_{x_0})$.

\medskip 

We define the \emph{cohomology groups of $B$ with values in the local system
$\mc S$}, denoted by $H^*(B;\mc S)$, as the cohomology groups of $B$
with values in the associated locally constant sheaf $\mc S$, and we
refer to~\cite{Steenrod} for a description in terms of singular
cochains with coefficients.  


\begin{ex} \label{ex:Leray Serre local system}
  \rm
  {\bf (The second term in the Leray-Serre spectral sequence).} Let $F
  \hookrightarrow E \stackrel \pi \longrightarrow B$ be a locally
  trivial fibration. For any $q\ge 0$ we define a local system $\mc H
  ^q(F,   \partial F)$ with fiber $H^q(F,   \partial F)$ as
  follows. The fiber at $x\in 
   B$ is $H^q(F_x,   \partial F_x)$. For a path $\alpha$ contained in
   a contractible open set $U \subset B$ we define $\Phi_\alpha =
   i_{\alpha(1)}^* {i_{\alpha(0)}^*}^{-1}$, where $F_{\alpha(0)}
     \stackrel {i_{\alpha(0)}} \hookrightarrow \pi ^{-1}(U) \stackrel
     {i_{\alpha(1)}} \hookleftarrow F_{\alpha(1)}$ are the inclusions,
     inducing isomorphisms in cohomology. This isomorphism is
     independent of $U$ as long as the latter is contractible. For a
     path $\alpha \in \mc P(B)$ we consider a subdivision $0 = t_0 <
     t_1 < \ldots < t_N < t_{N+1} = 1$ such that $\alpha_i = \alpha|_{[t_i, 
       t_{i+1}]}$ is contained in a contractible open set and define 
  $\Phi_\alpha = \Phi_{\alpha_N} \circ \Phi_{\alpha_{N-1}} \circ
  \ldots \circ \Phi_{\alpha_0}$. If $B$ is closed the second term of
  the Leray-Serre spectral sequence ${_{_{LS}}}E_r^{p,q} \Rightarrow
  H^{p+q}(E, \partial E)$ is 
$$
{_{_{LS}}}E _2^{p,q} \simeq H^p(B~;   \mc{H}^q(F,   \partial F)).
$$
\end{ex}

\medskip 

We define now \emph{Morse cohomology of a closed manifold $B$ with
coefficients in a local system $\mc S$}. Let $f: B \longrightarrow \RR$
be a Morse function and $Y$ be a 
Morse-Smale negative pseudo-gradient vector field. Pick an
orientation of the unstable manifolds of $Y$ and define the
cohomological Morse complex with values in $\mc S$ as 
$$
C^k(B;Y,\mc S) = \bigoplus _{\scriptsize \begin{array}{c} x
    \in \tx{Crit}(f) \\ \tx{ind}_{\tx{Morse}}(x) = k \end{array}} M_x,
$$
with differential $\partial : C^k \longrightarrow C^{k+1}$ given by  
\begin{equation} \label{eq:differential in Morse complex}
\partial (m\langle x \rangle ) = \sum_{\tx{ind}_{\tx{Morse}}(y)=k+1} \Big( \sum _{\gamma \in
\mathcal M(y,  x)} n_\gamma \Phi_\gamma ^{-1}(m) \Big) \langle y
\rangle .
\end{equation} 
Here $x\in \tx{Crit}(f)$, $m\in M_x$, 
$\tx{ind}_{\tx{Morse}}(x)=k$ and $n_\gamma$ is the
sign which is associated to the trajectory $\gamma$.
The fundamental identity $\partial ^2=0$ is proved by the usual gluing
argument, taking into account that the cancelling pairs of
trajectories form the boundary of a two-disc and therefore parallel
transport is the same along the two ``half-circles'' forming its
boundary. Similarly, one shows by a continuation argument that the
resulting cohomology groups do not depend on the choice of Morse
function, nor on the choice of pseudo-gradient vector field. 

Any proof showing that Morse cohomology with constant coefficients is
isomorphic to singular cohomology carries over to the case of locally
constant coefficients. The approach that is most convenient for us is
also the most geometric and uses \emph{cellular cohomology}
$H^*_{\tx{Cell}}(B;\mc S)$ as an intermediate device. Given a 
CW-decomposition of $B$, let $B^k$ be the $k$-skeleton and define the
cellular complex by $\tx{Cell}^k(B;\mc S)=H^k(B^k,B^{k-1};\mc S)$ with
differential 
$$
\partial_{\tx{Cell}}:H^k(B^k,B^{k-1};\mc S)\to H^{k+1}(B^{k+1},B^k;\mc S)
$$
given by the connecting homomorphism in the long exact sequence of the
triple $(B^{k+1},B^k,B^{k-1})$. It is shown
in~\cite[Appendix~A.4]{MS} that we have a canonical isomorphism
$H^*(\tx{Cell}^*(B;\mc S),\partial_{\tx{Cell}})\simeq H^*(B;\mc S)$ if
$\mc S$ is constant, but the proof carries over verbatim to an
arbitrary local system. 

The connection with Morse homology is realized by expressing the
cellular differential in an alternative way, using the incidence
numbers of the cells $e^k_i$ of the CW-decomposition. Let
$x^k_i$ be the center of the cell $e^k_i$. Each choice of orientation
of the cells determines an isomorphism $\tx{Cell}^k(B;\mc S) \simeq
\bigoplus_i M_{x^k_i}$ and incidence numbers $[e^k_i : e^{k+1}_j]$. 
The
differential $\partial _{\tx{Cell}}$ is then equal to  
\begin{equation} \label{eq:dcellinc}
\partial_{\tx{Cell}}(m_i)=\sum_j[e^k_i : e^{k+1}_j]\Phi_{ij}(m_i),
\qquad m_i\in M_{x^k_i},
\end{equation} 
where $\Phi_{ij}$ is parallel transport along a path from $x^k_i$ to
$x^{k+1}_j$ contained in the closure of $e^{k+1}_j$. In order for
parallel transport to be independent of the path it is enough that the
closure of $e^{k+1}_j$ be simply-connected. This can
fail only if $k=0$ and the endpoints of $e^1_j$ coincide with some
$e^0_i$, in which case the term $[e^0_i : e^1_j]\Phi_{ij}(m_i)$ has to
be replaced by $\langle e^0_i\rangle
(\Phi^+_{ij}(m_i)-\Phi^-_{ij}(m_i))$. Here $\Phi^\pm_{ij}$ is the path
running from $x^0_i$ to $x^1_j$ and having the same/the opposite
orientation as $e^1_j$, and $\langle e^0_i\rangle$ is the sign
(orientation) of $e^0_i$. We refer to~\cite[IV.10]{Br} for
a proof of~\eqref{eq:dcellinc} in the case where $\mc S$ is constant,
which carries over verbatim to the case of an arbitrary local system. 

We now use the fact proved by Laudenbach~\cite{Laudenbach92} that, if
the vector field  $Y$ is equal near its zeroes to the negative
gradient of a quadratic form with respect to the Euclidean metric,
then its unstable manifolds provide a CW-decomposition of $B$. It then
follows directly from the definitions that
$$
[W^u(x) : W^u(y)]=\sum_{\gamma\in\mc M(y,x)}n_\gamma.
$$
In particular, once an orientation of the unstable manifolds has been
chosen, the Morse complex $(C^*(B;Y,\mc S),\partial)$ is
canonically identified with the cellular complex
$(\tx{Cell}^*(B;\mc S),\partial_{\tx{Cell}})$, and we have canonical
isomorphisms 
$$
H^*(B;Y,\mc S)\simeq H^*_{\tx{Cell}}(B;\mc S)\simeq H^*(B;\mc S).
$$

\subsubsection{Local subsystems and extensions} \label{sec:localsubsystems}

The motivation for introducing \emph{local subsystems} 
is that parallel transport is
defined in the Floer setting only along certain paths in $B$. The
question arises whether such a system of isomorphisms can be
extended to a local system, and if yes, in how many non-isomorphic
ways. The notion of a local subsystem is a convenient way to
organize the available data.

\begin{defi} \label{defi:loc subsys}
  Let $R$ be a ring and $M$ an $R$-module. A \emph{local subsystem with
    fiber $M$} on the topological space $B$ consists of the following
  data. 
  \renewcommand{\theenumi}{\arabic{enumi}}
  \begin{enumerate}    
       \item a subset $C \subset B$ and one copy $M_x$ of $M$   
         for each $x \in C$;   
       \item \label{cond 2 in loc subsys} a subset $\mc{P}\subset \mc{P}(B)$ such that
             \begin{itemize}    
                \item if $\alpha \in \mc{P}$, then $\alpha(0),    
                  \alpha(1) \in C$;   
                \item for any $x\in C$ the constant path
                  $\underline{x}(\cdot) \equiv x$ belongs to $\mc{P}$;   
                \item $\alpha \in \mc{P}$ if and only if 
                  $\alpha ^{-1} \in \mc{P}$;   
                \item if $\alpha,   \beta \in \mc{P}$ and $\alpha(1)   
                  =\beta(0)$ then $\alpha \cdot \beta \in \mc{P}$;
             \end{itemize}    
       \item a family of isomorphisms
              $\Phi=\big( \Phi_{\alpha}~:M_{\alpha(0)}   
             \stackrel{\sim}{\longrightarrow}    
               M_{\alpha(1)} \big) _{\alpha \in \mc{P} }$ such that
           \begin{itemize}    
              \item if $\alpha \simeq \beta$ are homotopic in $B$ with fixed
                endpoints, then $\Phi_\alpha =   
                \Phi_\beta$;   
              \item if $\alpha,   \beta \in \mc{P}$ and
                 $\alpha(1) =\beta(0)$ then $\Phi_{\alpha \cdot   
                  \beta} \ = \ \Phi_\beta \circ \Phi_\alpha$.   
           \end{itemize}    
  \end{enumerate}   
\end{defi}

\begin{defi} 
We call the pair $\tx{Supp}(\mc S) = (C,   \mc P)$ \emph{the support}
of the local 
subsystem $\mc S=(C,   \mc P,   \Phi)$. 

The \emph{connected component $\mc S(x_0) = \big(C(x_0),   \mc
  P(x_0),   \Phi(x_0) 
\big)$ of $x_0$ in $\mc S$} is 
defined as 
$$
C(x_0) = \{ x\in C \ : \ \exists \ \alpha \in \mc P, \
\alpha(0)=x_0, \ \alpha(1)=x \},
$$
$$
\mc P(x_0) = \{ \alpha\in \mc P \ : \ \alpha(0), \ \alpha(1) \in
C(x_0) \},
$$
$$
\Phi(x_0) = \{ \Phi_\alpha \ : \ \alpha \in \mc P(x_0) \} .
$$

We say that $\mc S$ has \emph{connected support} if $C(x_0)=C$ for
some (and hence for any) $x_0\in C$. 

The \emph{fundamental group $\pi_1(\mc S,   x_0)$ of $\mc S$ at
  $x_0$} is defined as the 
set of homotopy classes in $B$ relative to $x_0$ of based loops in $\mc
P(x_0)$. 
Multiplication is given by the catenation of loops. 

Given two local subsystems $\mc S = (C,   \mc P,   \Phi)$ and $\mc
S' = (C',   \mc
P',   \Phi')$ we say that $\mc S'$ is an \emph{extension of $\mc S$},
and write $\mc S \prec \mc S'$, if $C\subset C'$, $\mc P \subset \mc
P'$ and $\Phi \subset \Phi'$.
\end{defi} 

Every local subsystem gives rise to a representation 
$$
\chi_{\mc S, x_0} : \pi_1(\mc S,   x_0) \longrightarrow
\tx{Aut}(M_{x_0}).
$$
Given $\mc S \prec \mc S'$ we have an obvious inclusion 
$\pi_1(\mc S,   x_0) \hookrightarrow \pi_1(\mc S',   x_0)$
which fits into the commutative diagram 
$$
\xymatrix{
\begin{array}{c}   
\pi_1(\mc S,\ x_0) \qquad  \end{array}
 \ar@{^(->}[dr] \ar[rr]^{\chi_{\mc S,x_0}} & & \tx{Aut}(M_{x_0})
\\
& \pi_1(\mc S',   x_0) \ar[ur]_{\chi_{\mc S',x_0}} & 
}
$$
In view of the fact that local systems are in one-to-one
correspondence with representations $\pi_1(B,  x_0) \longrightarrow
\tx{Aut}(M_{x_0})$, the following statement is tautological.

\begin{prop} \label{prop:extending_local_subsystems}
  Let $\mc S$ be a local subsystem having \emph{connected support}.
  Each extension of $\mc S$ to a local system corresponds to one and
  only one factorization 
$$
\xymatrix{
\begin{array}{c}   
\pi_1(\mc S,\ x_0) \qquad  \end{array}
 \ar@{^(->}[dr] \ar[rr]^{\chi_{\mc S,x_0}} & & \tx{Aut}(M_{x_0})
\\
& \pi_1(B,   x_0) \ar@{.>}[ur] & 
}
$$
\end{prop}

We see in particular that, if a local subsystem $\mc S$ with connected
support is such that the inclusion  $\pi_1(\mc S,   x_0)
\hookrightarrow \pi_1(B,   x_0)$ is an isomorphism, then $\mc S$
admits a unique extension to a local system.


\subsection{Filtered Floer complexes} \label{computation of the E 1 term} 

We construct in this section filtered Floer complexes
using appropriate Hamiltonians, pseudo-gradient vector
fields and almost complex structures. 

%
%
Let $f:B \longrightarrow \RR$ be a $C^2$-small Morse function, let
$J_B$ be a time-independent almost complex structure on $B$, and
assume that the gradient $\nabla ^{J_B} f$ is Morse-Smale and is equal
to the gradient of a quadratic form with respect to the Euclidean
metric near $\mathrm{Crit}(f)$. This last condition ensures
by~\cite{Laudenbach92} that the unstable manifolds of $-\nabla ^{J_B}
f$ define a CW-decomposition of $B$. By
Theorem~\ref{thm:transversality on the base} there exists $\tau >0$
such that all solutions of Floer's equation 
$$
u_s + J_B u_t = J_B X_f
$$
with period less than $\tau$ are time-independent and cut  the
defining equation  transversally. We note the fact that, upon
multiplying the function $f$ by a constant $c>0$, the corresponding
$\tau$ transforms as $\tau \mapsto \tau/c$. As the function $f$ gets
multiplied in the sequel by constants $0<c<1$, the bound $\tau>0$ can
be taken uniform with respect to $c$.

The constructions in Sections~\ref{sec:main section pseudo gr vector
  fields} and~\ref{sec:main section transv} provide a bound
  $\epsilon_0>0$, sequences $0<\epsilon_\nu\le\epsilon_0$ and
  $c_\nu>0$, $\nu\in\NN$, as well as sequences of $\tau$-periodic
  Hamiltonians $K_\nu$, split almost complex structures $J_\nu$, and
  vector fields $Y^\epsilon_\nu$, $0 < \epsilon \le \epsilon_0$,
  such that

\begin{enumerate}
\item \label{item:A} The $\tau$-periodic orbits of $\epsilon K_\nu$, $0 <
  \epsilon \le \epsilon_0$ are nondegenerate and located in the fibers
  lying over the critical points of $f$;
\item \label{item:B} The constant $\tau$-periodic orbits of $\epsilon K_\nu$
  are critical points of $K_\nu$, the vector field $Y^\epsilon_\nu$,
  $0<\epsilon\le\epsilon_0$ is a negative pseudo-gradient for $K_\nu$
  on $E$, and $Y^\epsilon_\nu$ is equal to the gradient of a quadratic
  form with respect to the Euclidean metric near
  $\mathrm{Crit}(K_\nu)$;
\item \label{item:C} The vector field 
$$
   \mc Y^{\epsilon}_\nu(x) = J_\nu \dot x - Y^\epsilon _\nu
   \circ x 
$$ 
defined on contractible $\tau$-periodic loops in
$\widehat E$ is a strong pseudo-gradient for
the action functional $A^\epsilon_{\epsilon K_\nu}$ for
$0<\epsilon\le\epsilon_0$; 
\item \label{item:D} For $0<\epsilon\le\epsilon_\nu$ the $\tau$-periodic solutions 
of Floer's equation 
\begin{equation} \label{eq:FloerY}
 u_s + J_\nu u_t = Y^\epsilon_\nu(t,   u(s,   t))
\end{equation} 
cut the equation transversally. The resulting Floer
complexes are denoted
$$
  FC^*(\epsilon,   \nu) = FC^*(Y^\epsilon_\nu,   J_\nu) ;
$$ 
\item \label{item:E} Solutions of~\eqref{eq:FloerY} project on gradient
  trajectories of $\epsilon c_\nu f$. We denote 
  $$
   C^*(\epsilon,\nu) = C^*(B;-\epsilon c_\nu \nabla^{J_B}f);
  $$
\item \label{item:F} The Hamiltonians $K_\nu$ form a cofinal sequence. 
\end{enumerate}

Statements~\eqref{item:A} and~\eqref{item:C} follow from
Proposition~\ref{prop:general ps gr}. Property~\eqref{item:B} is realized
by choosing the perturbations $f_i$ in Section~\ref{sec:the geometric
  Hamiltonians} quadratic near their critical
points. Statement~\eqref{item:D} follows from
Theorem~\ref{thm:transversality for split almost complex structures}
and by noticing that the perturbations of almost complex structures can
be taken trivial near the critical points of $K_\nu$, so that
property~\eqref{item:B} is preserved. Statement~\eqref{item:E} follows
from the special form of the vector fields
$Y^\epsilon_\nu$. Property~\eqref{item:F} is realized by choosing a
cofinal ``stem sequence'' $h_\nu$ in order to construct the
Hamiltonians $K_\nu$. 

\begin{rmk} {\bf (Rescaling)}
\rm Cofinality does not depend on the rescaling 
  $\Omega \rightsquigarrow  \epsilon \Omega$, because the coordinate $S$ depends
  only on the  
  vertical vector field $Z$, which remains unchanged. On the other hand,
  the action of the orbits corresponding to closed characteristics
  gets multiplied by a factor $\epsilon$ under the rescaling $\Omega
  \rightsquigarrow \epsilon \Omega$. 
\end{rmk} 

\begin{rmk} {\bf (Parameters $\epsilon_\nu\to 0$ are unavoidable)}
\rm The use of the constants $\epsilon_\nu$ cannot be avoided because of
transversality issues, and more precisely because of
Lemma~\ref{lem:transv up through trans down}. The constants $c_\nu$
have more of a formal role, mainly in order to ensure $K_\nu\le
K_{\nu'}$ for $\nu\le\nu'$ (cf. Section~\ref{sec:the geometric
  Hamiltonians}).   
\end{rmk}

 Let $R$ be a ring and $-\infty \le a < b \le \infty$ such that
 $\tau a,\tau b\notin\tx{Spec}(\partial E)$. We define
\begin{equation}
  \label{eq:filtration in constr of sp seq}
  F_pC^*_{[\epsilon a,   \epsilon b]}(\epsilon,   \nu)  \qquad = 
\bigoplus _ { \scriptsize
{\begin{array}{c} \alpha \in \mc{P}(K_\nu)
\\  
A^\epsilon_{\epsilon K_\nu}(\alpha) \in [\epsilon a,   \epsilon b ] \\
- i_{CZ} (\pi (\alpha) ) \ge p 
\end{array} 
} }
R \langle \alpha \rangle 
\qquad \subseteq \qquad  
FC^*_{[\epsilon a,  \epsilon b]}(\epsilon,   \nu) . 
\end{equation}
Here $\mc{P}(K_\nu)$ denotes the set of $\tau$-periodic orbits of
$K_\nu$. 

\begin{prop} \label{prop:filtration and differential}
  Formula~\eqref{eq:filtration in constr of sp seq} defines a
  filtration, i.e.
   $$
   \partial \big( F_pC^k_{[\epsilon a,   \epsilon b]}(\epsilon,   \nu)
   \big) \subset 
   F_pC^{k+1}_{[\epsilon a,   \epsilon b]}(\epsilon,   \nu) .
   $$
\end{prop}

\demo For $\alpha\in \mc P(K_\nu)$ such that $A^\epsilon_{\epsilon
  K_\nu}  (\alpha) \in [\epsilon a,   \epsilon b]$ and $-i_{CZ}(\pi(\alpha)\ge
  p$ we have  
\begin{equation*}
  \partial \langle \alpha\rangle = \sum 
    _{\scriptsize{\begin{array}{c} \beta \ : \
    -i_{CZ}(\beta)+i_{CZ}(\alpha)=1 \\
    A^\epsilon_{\epsilon K_\nu}(\beta)\in[\epsilon a,   \epsilon b]
  \end{array} }} \# \left(\mc M(\beta,  \alpha) /\RR \right) \ \langle \beta \rangle .
\end{equation*}
For each orbit $\beta$ appearing in the above sum there exists a solution
of Floer's equation $u_s+J_\nu u_t=Y^\epsilon_\nu(t,   u(s, 
t))$ running from $\beta$ to $\alpha$. It follows that $v=\pi\circ u$ is a
solution of the equation $v_s+J_Bv_t=\epsilon c_\nu J_B X_f(v)$ running
from $\pi(\beta)$ to $\pi(\alpha)$. By transversality we must have
$\tx{ind}_{\tx{Morse}}(\pi(\beta);    f) \ge
\tx{ind}_{\tx{Morse}}(\pi(\alpha);    f)$, with equality if and only
if $v$ is constant i.e. $u$ is entirely contained in the critical
fiber lying over the point $\pi(\beta)=\pi(\alpha)$. Equivalently, we get 
$-i_{CZ}(\pi(\beta) \ge -i_{CZ}(\pi(\alpha)\ge p$.
\hfill{$\square$}

\medskip 

The filtration~(\ref{eq:filtration in constr of sp seq}) defines a spectral sequence 
\begin{equation} \label{eq:Floer spectral sequence}
E_r^{p,q}(\epsilon,   \nu,   a,   b)
\end{equation} 
such that 
$$
E_r^{p,q}(\epsilon,   \nu,   a,  b) \Longrightarrow
FH^{p+q}_{[\epsilon a,  
\epsilon b]} (\epsilon,   \nu).
$$


\subsection{The terms $E_1$ and $E_2$}
The filtered complex $FC^*_{[\epsilon a,   \epsilon b]}(\epsilon,   \nu)$ 
falls under the formalism of Example~\ref{example:spectral
  sequence}. The proof of Proposition~\ref{prop:filtration and
  differential} shows that the component $\partial _0$ of the differential
is given by the Floer differential in the fibers, hence 
\begin{equation} \label{E1}   
E_1 \ \simeq \ \bigoplus_{p \in \ZZ}    
\bigoplus_{\left. \begin{array}{c}      
                                   _{ p_i\in \tx{Crit}(f)}\\     
                                   _{ -i_{CZ}(p_i)=p }    
                      \end{array} \right.}     
FH^*_{[\epsilon a,   \epsilon b]}(\widehat E_{p_i},  
J_\nu|_{\widehat E_{p_i}}, 
  \epsilon K_\nu|_{\widehat E_{p_i}}, \epsilon \Omega),
\end{equation}   
and $d_1: E_1 \longrightarrow E_1$ acts by
$$    
\underset{    
     \left. \begin{array}{c}    
             _{-i_{CZ}(\alpha)=k} \\    
             _{-i_{CZ}(\pi(\alpha))=p}    
           \end{array} \right.    
}    
{\underbrace{[\alpha]}}    
\longmapsto     
\Big[ \sum_    
{ \left. \begin{array}{c}    
           _{ -i_{CZ}(\beta)=k+1} \\    
           _{ -i_{CZ}(\pi(\beta)) =  p+1}     
         \end{array} \right. }    
   \# (\mc M(\beta, \alpha)/\RR)   \langle \beta \rangle  
   \Big] .$$     
Here the representative $\alpha$ of the cohomology class 
is an $R$-linear combination of periodic
orbits $\alpha_i$ with $-i_{CZ}(\alpha_i)=k$ and lying in the same
   fiber over a critical point of index $p$, while  $ \# (\mc M(\beta,
   \alpha)/\RR) $ represents, for some $\beta$ with $-i_{CZ}(\beta) = k+1$,    
the extension by linearity in the second argument of the quantity 
$ \# (\mc M(\beta, \alpha_i)/\RR) $.

The formula for  $d_1$ becomes transparent if we define the
{\it parallel transport map} 
\begin{equation} \label{eq:parallel transport map Floer} 
\Phi_\gamma^{p_j,p_i}(\epsilon,   \nu) : FC^q({\widehat
  E}_{p_i},   \epsilon,   \nu)
\longrightarrow FC^q({\widehat E}_{p_j},   \epsilon,   \nu)  ,
\end{equation}     
$$\langle \alpha \rangle \longmapsto    
\sum_{-i_{CZ}^{\tx{fiber}}(\beta)=q}
\#\mc M_{\gamma}^{\epsilon,\nu}(p_j\otimes
\beta,  p_i\otimes \alpha) \langle \beta \rangle  ,$$     
where  $p_i,p_j\in \tx{Crit}(f)$ with $-i_{CZ}(p_j)>-i_{CZ}(p_i)$,
$\gamma \in \mc M(p_j,   p_i)$, $FC^q({\widehat E}_{p_i},   \epsilon,
\nu)$ and $FC^q({\widehat E}_{p_j},   \epsilon,   \nu)$ are the Floer
complexes in the fibers ${\widehat E}_{p_i}$, ${\widehat E}_{p_j}$ for
the restrictions of $J_\nu$, $\epsilon K_\nu$ and $\epsilon
\Omega$. We have denoted  
\begin{eqnarray*}
\mathcal{M}_{\gamma}^{\epsilon,\nu}(p_j\otimes \beta,   p_i \otimes
\alpha) & = & 
\Bigg\{ u:\RR \times \Ss^1_\tau \longrightarrow \widehat E \
: \ \begin{array}{l} 
     u_s + J_\nu u_t = Y^\epsilon_\nu \circ u, \\
     \pi \circ u(s,  t) = \gamma(s),   \forall \ s,  t,  \\
     \displaystyle \lim_{s\rightarrow -\infty} u(s,\cdot)  =
       p_j\otimes \beta (\cdot), \\
     \displaystyle \lim_{s\rightarrow +\infty} u(s,\cdot)  =
       p_i\otimes \alpha (\cdot)    
    \end{array}
 \Bigg\}.
\end{eqnarray*} 
The notation $p_i \otimes
\alpha$ stands for the orbit $\alpha$ in the fiber ${\widehat
  E}_{p_i}$ viewed as an orbit in $\widehat E$, and $p_j \otimes
\beta$ has a similar meaning. It is easy to see that we have 
$$
\dim\,
\mathcal{M}_{\gamma}^{\epsilon,\nu}(p_j\otimes \beta,   p_i \otimes 
\alpha)=-i_{CZ}^{\tx{fiber}}(\beta) + i_{CZ}^{\tx{fiber}}(\alpha)
$$ 
regardless of the difference of indices between $p_j$ and $p_i$. 

We prove in Lemma~\ref{lem:first properties of Morse
  parallel transport} below that $\Phi_\gamma^{p_j,p_i}(\epsilon,\nu)$
is a morphism of differential complexes which induces an isomorphism
in homology and preserves the action filtration. As a consequence, the
differential $d_1$ can be rewritten
\begin{equation*}
  d_1(\epsilon,   \nu): \bigoplus _{p_i \in \tx{Crit}(f)}
  FH^*_{[\epsilon a,
    \epsilon b]}(\widehat 
  E_{p_i},   \epsilon,   \nu) \longrightarrow 
 \bigoplus _{p_i \in \tx{Crit}(f)} FH^*_{[\epsilon a,   \epsilon
  b]}(\widehat E_{p_i},   \epsilon,   \nu),
\end{equation*}
\begin{equation}
  \label{eq:differential d 1}
  p_i \otimes [\alpha] \longmapsto \sum _{-i_{CZ}(p_j) = -i_{CZ}(p_i)
    + 1} p_j \otimes \sum_{[\gamma] \in \mc M(p_j,   p_i)/\RR}
  n_\gamma \Phi_\gamma^{p_j,   p_i}(\epsilon,   \nu) \cdot [\alpha] .
\end{equation}
Here $\alpha \in FC^*_{[\epsilon a,   \epsilon 
b]}(\widehat E_{p_i},   \epsilon,  
\nu)$ is such that $\partial _{\tx{fiber}}(\alpha)=0$. We shall
moreover prove in Proposition~\ref{prop:Floer local system} below that
there is a unique local system $\mc{FH}^q_{[\epsilon a,\epsilon
  b]}(\widehat F,\epsilon,\nu)$ on $B$ with fiber $FH^q_{[\epsilon
  a,\epsilon b]}(\widehat F,\epsilon,\nu)$ so that the maps 
$\Phi_\gamma^{p_j,p_i}$ described above are the parallel transport
maps with respect to this local system. It then follows from the
definition~\eqref{eq:differential in Morse complex} of Morse
cohomology with values in a local system that 
\begin{equation} \label{eq:E2}
E_2^{p,q}(\epsilon,\nu,a,b)\simeq H(E_1,d_1) = H^p(B;\mc{FH}^q_{[\epsilon a,\epsilon b]}(\widehat
F,\epsilon,\nu)).
\end{equation}


\subsection{The Floer local system} 

We prove in this section that the parallel transport maps
$\Phi_\gamma^{p_j,p_i}(\epsilon,\nu)$ are chain morphisms which induce
isomorphisms in homology, and moreover they can be incorporated into a
uniquely determined local system on $B$. 

\begin{lem} \label{lem:first properties of Morse parallel transport}
  The map $\Phi_\gamma^{p_j,p_i}(\epsilon,  \nu)$
  in~\eqref{eq:parallel transport map Floer} is a morphism of
  differential complexes:  
$$
\Phi_\gamma^{p_j,p_i} \circ \partial_{p_i} + \partial _{p_j} \circ
\Phi_\gamma^{p_j,p_i} = 0,
$$
where $\partial _{p_i}$, $\partial _{p_j}$ are the Floer differentials
in the fibers $\widehat E _{p_i}$, $\widehat E_{p_j}$ respectively.
Moreover, it induces an isomorphism in Floer cohomology and preserves
the action filtration on Floer complexes. 
\end{lem}

\demo The idea is to identify the moduli spaces $\mc
M_\gamma^{\epsilon,\nu}(p_j\otimes \beta,   p_i\otimes \alpha)$ with
the moduli spaces $\mc M(\beta,   \alpha)$ corresponding to a
deformation of the Floer 
equations on the fibers $\widehat E_{p_i}$, $\widehat E_{p_j}$. More
precisely, let us choose a symplectic trivialization 
$$\Psi: \widehat E|_{\tx{im} \, \gamma} \longrightarrow \RR \times
\widehat F,$$
$$\widehat E_{\gamma(s)} \longrightarrow \{ s \} \times \widehat F$$
which, under the projection $\gamma(s)\longmapsto
s$ and after having chosen isomorphisms $\widehat F \simeq \widehat
E_{p_i} \simeq \widehat E_{p_j}$, 
coincides with the trivializations $\Psi_i,\Psi_j$ of Section~\ref{sec:the
  geometric Hamiltonians} on $\tx{im} \,  \gamma \cap U_i$,
respectively on $\tx{im} \,  \gamma \cap U_j$. In this 
trivialization we can interpret 
$$
\widehat K = K_\nu \circ \Psi^{-1}
$$ 
as an $s$-dependent deformation from $K_\nu \circ
\Psi_{p_i}^{-1} = \varphi+\chi^\delta_\nu +c_-$ to
$K_\nu\circ \Psi_{p_j}^{-1}=\varphi+\chi^\delta_\nu + c_+$,
where $\chi^\delta_\nu$ is the perturbation described in Section~\ref{sec:the
  geometric Hamiltonians} and $c_\pm$ are arbitrarily small
constants. Writing $\widehat 
K=(\widehat K_s)$, $s\in \RR$ with $\widehat K_s:\Ss^1_\tau
\times \widehat F\longrightarrow \RR$ and considering on $\{s\}\times
\widehat F$ the almost complex structure $J^{\tx{vert}}_s$ 
induced from $\widehat E_{\gamma(s)}$, the vector field
$Y^\epsilon_\nu$ takes the form  
$$
\widehat Y ^\epsilon_\nu (s,   p) = \nabla^{J^{\tx{vert}}_s}
\epsilon \widehat K_s + \DP{}{s}, 
$$
with $\nabla^{J^{\tx{vert}}_s}$ the Levi-Civita connection
for $\epsilon\Omega$. 
We infer that, if $\alpha$ and $\beta$ are two periodic orbits of
$\varphi+\chi^\delta_\nu$, the moduli space $\mc M_\gamma^{\epsilon,\nu}
(p_j\otimes \beta ,   p_i\otimes \alpha)$ is isomorphic to the
moduli space
$$\mc M(\beta,   \alpha; \widehat K) = \Big\{ u:\RR \times
\Ss^1_\tau \longrightarrow \widehat F \ : \ 
\begin{array}{l}
u_s + J_s^{\tx{vert}} u_t = \nabla ^{J^{\tx{vert}}_s} \epsilon \widehat K_s \\
u(-\infty)=\beta, \ u(+\infty)=\alpha
\end{array}
\Big\} .
$$
It is a standard fact in Floer theory that the count of the elements
of these moduli spaces gives rise to chain morphisms which induce
isomorphisms in homology. Moreover, parallel transport preserves the
filtration by the action if the constants $c_\pm$ are chosen
small enough.   
\hfill{$\square$}

\medskip 

In the next Lemma we omit the indices $\epsilon$ and $\nu$ to improve
readability.

\begin{lem} \label{lem:gluing property for parallel transport}
  Let $x,y,z \in \tx{Crit}(f)$ and $u\in \mc M(x,  y)$, $v\in \mc
  M(y,  z)$. Let $\gamma \in \mc M(x,  z)$ be such that 
  $[\gamma]\in \mc M(x,   z)/\RR$ belongs to a component whose
  boundary contains $([u],   [v]) \in 
(\mc M(x,   y)/\RR) \times (\mc M(y,   z)/\RR)$.
  The following equality holds:
 \begin{equation*}
   \Phi_{\gamma}^{x,z} =  \Phi_u^{x,y} \circ \Phi_v^{y,z} . 
 \end{equation*}
\end{lem}

\demo We first remark that $\Phi_\gamma^{x,z}$ is independent of
$\gamma$ as long as $[\gamma]$ varies in one
component of $\mc M(x,  z)/\RR$. Indeed, components are path
connected and we can choose an embedded path $[\gamma_\tau]$, $\tau\in
[0,   1]$ between any two given points
$[\gamma_0]$ and $[\gamma_1]$. A choice of symplectic trivialization of
$\widehat E$ over $\bigcup _\tau \tx{im} \, \gamma_\tau$ induces a
homotopy between the homotopies $(\widehat K_0, \widehat Y_0)$ and
$(\widehat K_1,   \widehat Y_1)$ corresponding to $\gamma_0$ and
$\gamma_1$ through Lemma~\ref{lem:first properties of Morse parallel
  transport}. The corresponding homomorphisms are then chain
homotopic and coincide at the level of homology. 

It is therefore enough to prove the claim if $[\gamma]$ 
lies in the image of the
gluing map 
$$\# : K \times [R_0,   \infty[
\longrightarrow \mc M(x,   z)/\RR,$$
where $K$ is a relatively compact neighbourhood 
of $([u],  [v])$ in $(\mc M(x, 
y)/\RR) \times (\mc M(y,   z)/\RR)$ and $R_0>0$ is large enough. 
The key point, borrowed from the proof of the invariance of
Floer homology~\cite[Lemma~3.11]{Salamon-lectures}, 
is that the morphism $\Phi_u^{x,y} \circ \Phi_v^{y,z}$ 
is induced by the gluing of the
two homotopies $(\widehat K_0,   \widehat Y_0)$, $(\widehat K_1,  
\widehat Y_1)$ into 
$$(\widehat K_{01},   \widehat Y_{01}) := \left\{ \begin{array}{ll} 
(\widehat K_0,   \widehat Y_0)(s+R), & s\le 1\\
(\widehat K_1,   \widehat Y_1)(s-R), & s\ge -1 .
\end{array} \right.
$$
Here $R$ is chosen large enough so that the resulting homotopy
is regular. Ignoring the additive constants as in 
Lemma~\ref{lem:first properties of Morse parallel
  transport}, this homotopy is a deformation of the constant
Hamiltonian $\varphi+\chi_\nu^\delta$.

For $R$ large enough the image of
$[\gamma] = \#\big(([u],   [v]),   R \big)$ lies in a contractible
neighbourhood $\mc U$ of $\tx{im} \, u   \cup   \tx{im} \,  v$ and we
choose a symplectic trivialization of $\widehat E$ over $\mc U$ which
extends the trivialization over $\tx{im} \, u   \cup \tx{im} \,  v$ which was
implicit in the construction of $(\widehat K_{01},   \widehat
Y_{01})$. The morphism
$\Phi^{x,z}_\gamma$ now arises from a homotopy $(\widehat K_\gamma,  
\widehat Y_\gamma)$ which, in the given trivialization, is also a
deformation of the constant Hamiltonian
$\varphi+\chi_\nu^\delta$. 

As a conclusion, both homotopies induce the identity in the given
trivialization, and in particular $\Phi_\gamma^{x,z}=\Phi_u^{x,y}\circ
\Phi_v^{y,z}$.  
\hfill{$\square$}

\medskip

In order to incorporate the maps $\Phi_\gamma^{p_j,p_i}(\epsilon,\nu)$
in local systems for the various values of $q\in\ZZ$, we exhibit local
subsystems $(C,\mc P,\Phi^q)$
(cf. Section~\ref{sec:localsubsystems}) determined by these maps. We
recall that we denote by $\mc P(B)$ the space of continuous paths in
$B$. We define 
$$
C=\tx{Crit}(f), \quad \mc P = \langle \tx{ negative gradient
  trajectories of } f \ \rangle_{\tx{Crit}(f)},
$$
where, for $\mc R \subset \mc P(B)$ and $C \subset B$, the notation
$\langle \mc R \rangle _C$ stands for the minimal subset of $\mc P(B)$
which contains $\mc P$ and which satisfies condition~(\ref{cond 2 in
  loc subsys}) of Definition~\ref{defi:loc subsys}.  In our case $\mc
P$ consists of chains of negative gradient trajectories of $f$ and
their inverses. We define 
$$
\Phi^q = \big(\Phi_\alpha : FH^q_{[\epsilon a,   \epsilon b]}(\widehat
E_{\alpha(0)},   \epsilon,   
\nu) \longrightarrow FH^q_{[\epsilon a,   \epsilon b]}(\widehat
E_{\alpha(1)},   \epsilon,   
\nu)\big)_{\alpha \in \mc P}
$$
on the generators of $\mc P$ by
$$
\Phi_\gamma = \Phi_\gamma^{\gamma(-\infty), 
  \gamma(+\infty)}(\epsilon,   \nu), \quad \gamma \tx{ negative
  gradient trajectory of } f.
$$

The next result is of a topological nature, although its statement
involves Floer homology groups. The transition to Floer homology is
realized through a repeated use of Lemma~\ref{lem:gluing property for
  parallel transport}. 
 
\begin{prop} \label{prop:Floer local system}  
  \renewcommand{\theenumi}{\alph{enumi}}
  \begin{enumerate}
  \item \label{first item local system} 
    For any $q\in \ZZ$ and $-\infty\le a < b \le \infty$ the triple
    $(C,   \mc  P,   \Phi^q)$ defines a local subsystem with fiber
    $FH^q_{[\epsilon a, \epsilon b]}(\widehat F,   \epsilon,   \nu)$,
    denoted by 
$$
\underline{\mc{FH}}^q_{[\epsilon a,  
      \epsilon b]}(\widehat F,   \epsilon,   \nu);
$$
  \item \label{second item local system} 
 If $B$ is connected, the support of the above local
    system is connected;
  \item \label{third item local system}
  For any  $p_0 \in \tx{Crit}(f)$ the canonical inclusion 
     $$\pi_1(\underline{\mc{FH}}^q_{[\epsilon a,  
      \epsilon b]}(\widehat F,   \epsilon,   \nu),   p_0)
    \hookrightarrow \pi_1(B,   p_0)
     $$
    is an isomorphism. 
  \end{enumerate}
\end{prop}

\begin{rmk} \label{rmk:truncFlocsys} {\bf (Truncated Floer local system)} \rm
Assertion~(\ref{third item local system}) implies together with
    Proposition~\ref{prop:extending_local_subsystems} that, for each
    $q\in\ZZ$, there is a unique local system on $B$ extending
    $\underline {\mc{FH}}^q_{[\epsilon a,\epsilon b]}(\widehat F,   \epsilon,
    \nu)$. We call it the \emph{truncated Floer local system} and denote it by
  $$\mc{FH}^q_{[\epsilon a, \epsilon b]}(\widehat F,   \epsilon,   \nu). 
  $$
\end{rmk} 

\noindent {\it \small Proof of Proposition~\ref{prop:Floer local
    system}.} 
Our proof crucially uses Lemma~\ref{lem:gluing property for
  parallel transport}, as well as the assumption that the unstable
manifolds of $-\nabla^{J_B}f$ provide a CW-decomposition of $B$
(cf. Section~\ref{computation of the E 1 term}). We denote the
$k$-skeleton by  
$$
B^k = \bigcup_{\tx{ind}_{\tx{Morse}}(x)\le k} W^u(x) .
$$ 

(\ref{first item local system}) We have to prove 
that parallel transport along a loop $\alpha\in \mc P$
which is null-homotopic {\it in $B$} is trivial.
Let $\alpha = \gamma_0\gamma_1 \ldots  \gamma_N$, where
$\gamma_i$ is a trajectory of  $\pm \nabla f$ running from $p_i$ to
$p_{i+1}$ and $p_0=p_{N+1}$. By
adding a chain $\eta_0\ldots \eta_\ell
\eta^{-1}_\ell \ldots \cdot \eta_0$ running from $p_0$ to a local
minimum of $f$ we can assume that
$\tx{ind}_{\tx{Morse}}(p_0)=0$. We are interested in the monodromy
along $\alpha$ as an automorphism of $FH^*(\widehat E_{p_0})$. By
successively applying Lemma~\ref{lem:gluing property for parallel
  transport} and deforming the $\gamma_i$'s to the boundary of $\mc M(p_i,
  p_j)$ we can assume that all trajectories $\gamma_i$ are of
index one. We can moreover cancel pairs $\gamma_i\gamma_{i+1}$, $i\neq
0,N$ which
are of the form $\eta\eta^{-1}$. 

{\it Claim 1.} There exists a null-homotopic 
chain $\alpha_1$ based at $p_0$ which consists of
trajectories connecting critical points of index $0$ and $1$ 
and which satisfies  $\Phi_{\alpha_1}=\Phi_\alpha$. 

\begin{figure}[h]   
        \begin{center}  
\input{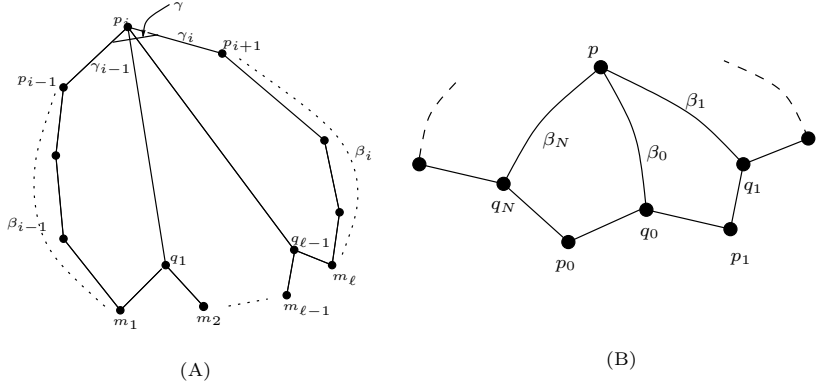} 
\caption{Parallel transport is trivial along null-homotopic loops. 
         \label{Picture for local system}
        }   
        \end{center}  
\end{figure}

The proof goes by induction over $m(\alpha)=\max_i   
\tx{ind}_{\tx{Morse}}(p_i)$. Assume $m(\alpha)\ge 2$. We prove the
existence of a null-homotopic chain $\alpha'$ based at $p_0$ such that
$m(\alpha')=m(\alpha)-1$ and $\Phi_{\alpha'}=\Phi_\alpha$. Let us
choose $p_i$ such that $\tx{ind}_{\tx{Morse}}(p_i)=m(\alpha)$. We
claim that the pair $\gamma_{i-1}\gamma_i$ connecting $p_{i-1}$, $p_i$
and $p_{i+1}$ can be replaced without affecting parallel transport 
by a chain connecting $p_{i-1}$ and
$p_{i+1}$ with intermediate critical points of index at most
$m(\alpha)-1$. Moreover, this chain is homotopic to
$\gamma_{i-1}\gamma_i$ with fixed endpoints, hence the resulting loop
is still null-homotopic. Let us choose a path 
$\gamma:[0,  1] \longrightarrow W^u(p_i)$ with $\gamma(0) \in \tx{im}
\,  \gamma_{i-1}$, $\gamma(1)\in \tx{im} \,  \gamma_i$ and which is
transverse in $W^u(p_i)$ to the manifolds $W^u(p_i) \cap W^s(q)$,
$q\in \tx{Crit}(f)$. In particular there exist points
$0=t_0<t_1<\ldots <t_{\ell-1}<t_\ell=1$ such that the trajectory
$[p_i,  \gamma(t_j)[$, $1\le j \le \ell-1$ 
lands on a critical point $q_j$ of index $1$ and, for every $t\in
]t_{j-1},  t_j[$, $1\le j\le \ell$ 
the trajectory $[p_i,   \gamma(t)[$ lands on a critical point $m_j$ of
index $0$ (see Figure~\ref{Picture for local system}(A)). The boundaries of the moduli spaces $\mc M(x, 
y)/\RR$ consist of broken trajectories and, together with the existence of
the curve $\gamma$, this ensures that the
trajectories $\gamma_{i-1}$, $[p_i,  q_j]$, $1\le j\le \ell -1$
and $\gamma_i$ can be completed to broken chains which
lie pairwise in the boundary of connected components of spaces 
$\mc
  M(p_i,   m)/\RR$, with $m$ a critical point of index $0$. 
More
precisely, there exist index decreasing 
chains of index $1$ trajectories from $p_{i-1}$
to $m_1$ and from $p_{i+1}$ to $m_\ell$, denoted $\beta_{i-1}$ and
$\beta_i$, as well as trajectories $[m_j,   q_j]$, $[q_j,  
m_{j+1}]$, $1\le j\le \ell-1$ such that the pairs 
$(q_{i-1}^{-1}\beta_{i-1},   [p_i,  q_1] \cdot [q_1,  m_1])$,
$([p_i,  q_j] \cdot [q_j,  m_{j+1}],   [p_i,  q_{j+1}] \cdot
[q_{j+1},   m_{j+1}])$, $1\le j \le \ell-2$ and $([p_i,  
q_{\ell-1}] \cdot [q_{\ell-1},   m_\ell],   \gamma_i \cdot
\beta_i)$ lie in the boundary of the same components of $\mc M(p_i,  m_j)/\RR$, $1\le j\le \ell$ respectively. The chain
$\gamma_{i-1}\gamma_i$ is now to be replaced by 
$$\beta_{i-1}\ \cdot 
\prod _{1\le j\le \ell-1} [m_j,  q_j] \cdot[q_j,  m_{j+1}] \quad
\cdot \  \beta_i^{-1} .$$
The construction can be repeated until all points of index
$m(\alpha)$ are eliminated.

{\it Claim 2.} We have $\Phi_{\alpha_1}=\tx{Id}$. 

For any CW-decomposition $B=\bigcup_{k\ge 0}B^k$ 
the group $\pi_1(B^1)$ is free, the group
$\pi_1(B^2)$ is the quotient of $\pi_1(B^1)$ by the normal subgroup
generated by the boundary cycles of the $2$-cells and the map 
$\pi_1(B^2)\longrightarrow \pi_1(B)$ induced by the inclusion
$B^2\hookrightarrow B$ is an isomorphism. In order to prove that
parallel transport is trivial along the null-homotopic loop $\alpha_1$
we can therefore assume without loss of generality that $\alpha_1$ is
the boundary cycle of a $2$-cell $W^u(p)$,
$\tx{ind}_{\tx{Morse}}(p)=2$ i.e. we can write 
$$\alpha_1 = \prod _{j=0}^N [p_j,   q_j] \cdot [q_j,   p_{j+1}],$$
with $\tx{ind}_{\tx{Morse}}(p_j)=0$, $\tx{ind}_{\tx{Morse}}(q_j)=1$,
$p_{N+1}=p_0$. Each $q_j$ is understood to be equipped with a trajectory
$\beta_j=[p,  q_j]$ such that 
$$\bigcup_j\beta_j =
\bigcup_{\tx{ind}_{\tx{Morse}}(q)=1} W^u(p)\cap W^s(q) .$$
The ordering of the $q_j$'s is such that $\beta_j$ and $\beta_{j+1}$,
$0\le j\le N$ lie
in the boundary of one same component of $\mc M(p,   m)$,
$\tx{ind}_{\tx{Morse}}(m)=0$ (see Figure~\ref{Picture for local system}(B)).

A deformation argument based on Lemma~\ref{lem:gluing property for parallel
  transport} and entirely similar to that of Case 1. shows that
parallel transport along $\alpha_1$ is the same as parallel transport
along the loop 
$$[p_0,  q_0] \cdot \prod_{j=0}^{N-1} \beta_j^{-1} \beta_{j+1} \cdot
[q_N,  p_0] = [p_0,  q_0] \cdot \beta_0^{-1}\beta_N \cdot [q_N, 
p_0] .$$
The last loop is the boundary of a component of $\mc M(p,  p_0)$ and,
again by Lemma~\ref{lem:gluing property for parallel
  transport}, the induced parallel transport is trivial.

b) Every critical point is connected by a trajectory to a point of
index $0$. On the other hand, because $B$ is connected and the
unstable manifolds form a CW-decomposition, the $1$-skeleton has to be
connected and any two index $0$ points are therefore connected by a chain of
index $1$ trajectories. This shows that the support of the local
system is connected.

c) We have to prove surjectivity of the map under study. We
use the fact that the map $\pi_1(B^1)\longrightarrow \pi_1(B)$ is
surjective, which means that every 
homotopy class in $B$ has a representative
which is supported in the $1$-skeleton $B^1=\bigcup
_{\tx{ind}_{\tx{Morse}}(p)=1}W^u(p)$. From this it is easy to 
find a representative
given by a chain of trajectories connecting points of index $0$ and
$1$. 
\hfill{$\square$}


\subsection{Proof of the main theorems} \label{sec:proofs}

We prove in this section Theorems A and B stated in the Introduction. 
We need four Lemmas which
describe the behaviour of 
the homology groups $FH^*_{[\epsilon a,\epsilon
  b]}(\epsilon K_\nu,\omega_\epsilon)$, of the complexes
$FC^*_{[\epsilon a,\epsilon b]} (Y^\epsilon_\nu,J_\nu)$, and of the
local systems $\mc{FH}^q_{[\epsilon a,\epsilon b]}(\widehat
F,\epsilon\nu)$ as one of the parameters $\epsilon$ and $\nu$ varies
and the other one is fixed. We choose $-\infty \le a < b \le \infty$
such that $\tau a,\tau b\notin \tx{Spec}(\partial E)$.

\begin{lem} \label{lem:on inverse limit with big epsilon}
Let $\nu \in \NN$ be fixed. For any $0< \epsilon <\epsilon' \le
 \epsilon_0$ we have natural isomorphisms
 $$FH^*_{[\epsilon a,  \epsilon b]}(\epsilon K_\nu, \omega_\epsilon)
 \stackrel {\psi_\nu^{\epsilon',  \epsilon}} 
 \longrightarrow FH^*_{[\epsilon' a,  \epsilon' b]}(\epsilon' K_\nu,
 \omega_{\epsilon'})  $$ 
 which satisfy $\psi_\nu^{\epsilon'',  \epsilon'}\circ
 \psi_\nu^{\epsilon',   \epsilon} = \psi_\nu^{\epsilon'', 
   \epsilon}$, $\epsilon < \epsilon' <\epsilon''$.
 The Floer complexes are based on $\tau$-periodic orbits and are
 understood to involve transverse almost complex structures. 
\end{lem} 

\demo The periodic orbits involved in the two complexes are the same.
Let us consider a homotopy of Hamiltonians and symplectic forms
$(\epsilon(s)K_\nu,\omega_{\epsilon(s)})$, $s\in \RR$ from $(\epsilon
K_\nu,\omega_\epsilon)$ to $(\epsilon' K_\nu,\omega_{\epsilon'})$. The
periodic orbits involved in the Floer complexes 
remain geometrically the same during the homotopy, 
but their action gets multiplied by $\epsilon(s)$. As a consequence,
the extremities of the action interval $[\epsilon(s)a,  
\epsilon(s)b]$ which interpolates between $[\epsilon a,\epsilon b]$
and $[\epsilon'a , \epsilon' b]$
are not crossed by any periodic orbit.
This ensures that the continuation morphism on Floer homology is
bijective~\cite[Prop.~1.1]{Viterbo99}.
\hfill{$\square$}

\begin{lem} \label{lem:filtered chain morphisms with small epsilon}
Let $\nu \in \NN$ be fixed. For any $0< \epsilon < \epsilon' \le
  \epsilon_\nu$ there are chain equivalences of filtered chain
  complexes 
$$
\xymatrix{
FC^*_{[\epsilon a,  \epsilon b]}(Y^\epsilon_\nu,J_\nu)
\ar[rr]^{\psi_\nu^{\epsilon',   \epsilon}} & & 
FC^*_{[\epsilon' a,  \epsilon' b]}(Y^{\epsilon'}_\nu,J_\nu),
}
$$
which induce in homology the isomorphisms $\psi_\nu^{\epsilon', 
  \epsilon}$ of Lemma~\ref{lem:on inverse limit with big epsilon}. 
\end{lem} 

\demo Let us consider a homotopy of pseudo-gradient vector fields 
and symplectic forms $(Y_\nu^{\epsilon(s)},\omega_{\epsilon(s)})$,
$s\in \RR$ between  
$(Y_\nu^\epsilon,\omega_\epsilon)$ and
$(Y_\nu^{\epsilon'},\omega_{\epsilon'})$. The induced chain map is
given by a count of solutions of the parametrized Floer equation $u_s+
J_s u_t =Y_\nu^{\epsilon(s)}(t,u(s,t))$ for some split almost complex
structure $J_s=\widetilde J_B \oplus J_s^{\tx{vert}}$. The projection
$v=\pi\circ u$ of such a solution solves an equation of the form
\begin{equation} \label{eq:epsilongrad} 
v_s+J_B v_t = \epsilon(s) J_BX_f.
\end{equation}  
Since the almost complex structure $J_B$ can be chosen to be regular
for~\eqref{eq:epsilongrad}, the induced 
chain map preserves the filtration. Note that solutions
of~\eqref{eq:epsilongrad} are actually reparametrized negative
gradient trajectories of $f$. That the induced morphism in homology
coincides with the map $\psi_\nu^{\epsilon',\epsilon}$ of
Lemma~\ref{lem:on inverse limit with big epsilon} is the usual
directed simple system property of Floer homology. 

The chain morphism $\psi_\nu^{\epsilon',   \epsilon}$ admits an
inverse up to chain homotopy obtained by considering the reversed
homotopy from $(Y_\nu^{\epsilon'},\omega_{\epsilon'})$ to
$(Y_\nu^\epsilon,\omega_\epsilon)$, and is therefore a chain
equivalence. 
\hfill{$\square$}

\begin{lem} \label{lem:filtered monotonicity morphisms}
Let $\nu\le \nu'$, so that $\epsilon_{\nu'}<\epsilon_\nu$. For any
  $0<\epsilon\le\epsilon_{\nu'}$ we have
  monotonicity morphisms 
\begin{equation} \label{eq:the special monotonicity morphism} 
\xymatrix{
FC^*_{[\epsilon a,\epsilon b]}(Y^\epsilon_{\nu'},J_{\nu'})
\ar[rr]^{\sigma^\epsilon_{\nu,\nu'}}  
& & 
FC^*_{[\epsilon a,\epsilon b]}(Y^\epsilon_\nu,J_\nu)
}
\end{equation}
which preserve the filtrations. 
\end{lem}

\demo We recall that
$$
Y_\nu^\epsilon = J_\nu(X^\epsilon_{\epsilon h_\nu} + \widetilde
{X_{\epsilon c_\nu f}} + (X^\epsilon_{\epsilon G_\nu})^{\tx{v}}),
\qquad 
Y_{\nu'}^\epsilon = J_{\nu'}(X^\epsilon_{\epsilon h_{\nu'}} + \widetilde
{X_{\epsilon c_{\nu'} f}} + (X^\epsilon_{\epsilon
  G_{\nu'}})^{\tx{v}}). 
$$
The key point is to deform  $Y_\nu^\epsilon$ to $Y_{\nu'}^\epsilon$
through a deformation which is ``good'' in the sense of
Definition~\ref{defi:good}. This will be achieved by catenating three
good deformations obtained from the following data: 
\begin{itemize} 
\item a linear interpolation between $Y_\nu^\epsilon$ and $J_\nu
  X^\epsilon_{\epsilon h_\nu}$; 
\item an increasing homotopy $h(s,  S)$ between $h_\nu$ and $h_{\nu'}$ through
  convex Hamiltonians which are linear of slope
  $\lambda_{\tx{max}}(s)$ for $S\ge 1$;
\item a linear interpolation between $J_{\nu'} X^\epsilon_{\epsilon
    h_{\nu'}}$ and $Y_{\nu'}^\epsilon$. 
\end{itemize} 

Firstly, each of the two linear interpolations defines a strong
pseudo-gradient for each value of the deforming parameter $s\in \RR$
and, moreover, the strong pseudo-gradient inequality is satisfied with
a uniform constant by Proposition~\ref{prop:general ps gr}. As a
consequence, the two interpolations define good deformations in the
sense of Definition~\ref{defi:good}. Secondly, any choice of
interpolation between $J_\nu$ and $J_{\nu'}$ defines a genuine
gradient deformation from $\nabla ^{J_\nu} h_\nu$ to
$\nabla^{J_{\nu'}} h_{\nu'}$, and in particular a good deformation in
the sense of Definition~\ref{defi:good}, with uniform constant equal
to $1$ in the pseudo-gradient inequality. By catenation we obtain
a good deformation from $Y_\nu^\epsilon$ to $Y_{\nu'}^\epsilon$, which
moreover satisfies the strong pseudo-gradient inequality with a
uniform constant.  

Arguing as in Section~\ref{sec:main section transv} one can show that 
transversality can be achieved by a generic choice of homotopies of vertical
almost complex structures and time-dependent horizontal distributions.
The trajectories of the Floer equation on $\widehat E$ project on
solutions of the equation 
\begin{equation} \label{eq:gradpar}
v_s +J_B v_t = \epsilon c(s) J_B X_f
\end{equation}
on the base, for some smooth function $c:\RR\to[0,\infty[$ which is
nonzero near $\pm\infty$. 
Since the almost complex structure $J_B$ can be chosen to be regular
for~\eqref{eq:gradpar}, 
the corresponding monotonicity morphisms preserve the filtrations.  
Note that solutions of~\eqref{eq:gradpar} are actually reparametrized
negative gradient lines of $f$.  
\hfill{$\square$}

\begin{lem} \label{lem:varying local system}
 (a) Let $\nu\in \NN$ and $q\in \ZZ$ be fixed. For any $0< \epsilon < \epsilon' \le
  \epsilon_\nu$ the filtered chain equivalences
  $\psi_\nu^{\epsilon',\epsilon}$ induce natural isomorphisms of local systems 
$$
\xymatrix{
\mc{FH}^q_{[\epsilon a,\epsilon b]}(\epsilon,\nu)
\ar[rr]^{\psi_\nu^{\epsilon',\epsilon}} & &
\mc{FH}^q_{[\epsilon' a,\epsilon' b]}(\epsilon',\nu).
}
$$

 (b) Let $q\in \ZZ$ and $\nu\le \nu'$, so that
  $\epsilon_{\nu'}<\epsilon_\nu$. For any
  $0<\epsilon\le\epsilon_{\nu'}$ the monotonicity morphisms
  $\sigma^\epsilon_{\nu,\nu'}$ induce morphisms of local systems  
  $$
\xymatrix{
\mc{FH}^q_{[\epsilon a,\epsilon b]}(\epsilon,\nu')
\ar[rr]^{\sigma^\epsilon_{\nu,\nu'}}  
& & 
\mc{FH}^q_{[\epsilon a,\epsilon b]}(\epsilon,\nu).
}
  $$
\end{lem}

\demo (a) The filtered chain equivalence $\psi_\nu^{\epsilon',\epsilon}$ of
Lemma~\ref{lem:filtered chain morphisms with small epsilon} induces an
isomorphism of spectral sequences
$$
\psi_r:E_r^{p,q}(\epsilon,\nu,a,b)\to
E_r^{p,q}(\epsilon',\nu,a,b), \qquad r\ge 1.
$$ 
The terms $E_1$ are direct sums of stalks of the corresponding local
systems, and we need to show that parallel transport commutes with
$\psi_1$. We recall from Lemma~\ref{lem:first properties of Morse
  parallel transport} that, in a suitable trivialization, the parallel
transport map is induced by a deformation of the Floer
equation. On the other hand, it follows from the proof of
Lemma~\ref{lem:filtered chain morphisms with small epsilon} that the
isomorphism $\psi_1$ is induced by deforming the Floer equation in
each fiber $\widehat E_p$, $p\in\mathrm{Crit}(f)$. The compositions
$\Phi(\epsilon',\nu)\circ \psi_1$ and $\psi_1\circ \Phi(\epsilon,\nu)$
are therefore induced by two deformations of the Floer equation having
the same endpoints, and as such coincide at the level of homology. 

(b) Let us denote by 
$$
\sigma_r:E_r^{p,q}(\epsilon,\nu',a,b)\to
E_r^{p,q}(\epsilon,\nu,a,b), \qquad r\ge 1
$$
the map of spectral sequences induced by
$\sigma^\epsilon_{\nu,\nu'}$. We need to show that
$$
\Phi(\epsilon,\nu)\circ \sigma_1 = \sigma_1\circ
\Phi(\epsilon',\nu),
$$ 
and the argument is entirely similar to the one given at (a). The only
difference is that $\sigma_1$ is not an isomorphism anymore, and this
is reflected in the weaker statement we need to prove, namely the
existence of a morphism of local systems, rather than the existence
of an isomorphism. 
\hfill{$\square$}

\begin{rmk} \label{rmk:Flocsys} {\bf (Floer local system)} \rm
As a consequence of Lemma~\ref{lem:varying local system} we define the
\emph{Floer local system} 
$$
\mc{FH}^q_{[a,b]}(\widehat F)=\lim_{\stackrel \longleftarrow \nu}
\mc{FH}^q_{[\epsilon_\nu a,\epsilon_\nu b]}(\epsilon_\nu,\nu),
$$
where the inverse limit is considered with respect to the maps  
$$
\xymatrix{
\mc{FH}^q_{[\epsilon_{\nu'} a,\epsilon_{\nu'}b]}(\epsilon_{\nu'},\nu')
\ar[rrr]^{\psi_\nu^{\epsilon_\nu,\epsilon_{\nu'}} \circ \,
  \sigma_{\nu',\nu}^{\epsilon_{\nu'}}} & & & 
\mc{FH}^q_{[\epsilon_\nu a,\epsilon_\nu b]}(\epsilon_\nu,\nu), \qquad \nu\le\nu'.
}
$$
\end{rmk}

\noindent {\it \small Proof of Theorem~A.}
We first show the existence of the requested spectral sequence. 
Let us consider the following diagram of chain complexes:

{\scriptsize
\begin{equation*}
\xymatrix
@C=17pt
@R=20pt@W=1pt@H=1pt
{FC^*_{[\epsilon_1 a,  \epsilon_1 b]}(Y^{\epsilon_1}_1, 
  J_1) 
& FC^*_{[\epsilon_1 a,  \epsilon_1 b]}(\epsilon_1 K_2, 
\omega_{\epsilon_1}) \ar@{.>}[l] 
& FC^*_{[\epsilon_1 a,  \epsilon_1 b]}(\epsilon_1 K_3,  
\omega_{\epsilon_1}) \ar@{.>}[l] 
&  \ar@{.>}[l] \\
FC^*_{[\epsilon_2 a,  \epsilon_2 b]}(Y^{\epsilon_2}_1, 
  J_1) \ar[u]
& FC^*_{[\epsilon_2 a,  \epsilon_2 b]}(Y^{\epsilon_2}_2, 
  J_2) \ar[l] \ar@{.>}[u] 
& FC^*_{[\epsilon_2 a,  \epsilon_2 b]}(\epsilon_2 K_3,  
\omega_{\epsilon_2}) \ar@{.>}[l] \ar@{.>}[u] 
&  \ar@{.>}[l] \\
FC^*_{[\epsilon_3 a,  \epsilon_3 b]}(Y^{\epsilon_3}_1, 
  J_1) \ar[u]
& FC^*_{[\epsilon_3 a,  \epsilon_3 b]}(Y^{\epsilon_3}_2, 
  J_2) \ar[l] \ar[u]
& FC^*_{[\epsilon_3 a,  \epsilon_3 b]}(Y^{\epsilon_3}_3, 
  J_3) \ar[l] \ar@{.>}[u] 
&  \ar@{.>}[l] \\
\vdots \ar[u] & \vdots \ar[u] & \vdots \ar[u] & 
}
\end{equation*}
}

The \emph{vertical}
arrows induce in cohomology the isomorphisms $\psi_\nu^{\epsilon', 
  \epsilon}$, $0<\epsilon < \epsilon' \le \epsilon_0$ of
Lemmas~\ref{lem:on inverse limit with big epsilon}
and~\ref{lem:filtered chain morphisms with small epsilon}. Moreover,
the ones drawn with continuous lines preserve the filtrations given
by~(\ref{eq:filtration in constr of sp seq}).   

The \emph{horizontal} arrows are given by monotonicity morphisms, and
Lemma~\ref{lem:filtered monotonicity morphisms} ensures that the ones
drawn with continuous lines preserve the filtrations.

By naturality of the continuation and
monotonicity morphisms in Floer homology, the above diagram 
commutes up to chain homotopy, and therefore commutes in homology.

Let us denote 
$$FC^*(\nu,a,b)= FC^*_{[\epsilon_\nu a,  \epsilon_\nu
  b]}(Y^{\epsilon_\nu}_\nu,  
  J_\nu) .$$
According to the diagram there are filtered chain morphisms 
$$FC^*(\nu,a,b) \longleftarrow FC^*(\nu',a,b), \ \nu \le
\nu'$$
which induce morphisms between the associated spectral
sequences~(\ref{eq:Floer spectral sequence}) 
$$E_r^{p,q}(\nu,a,b) \longleftarrow
E_r^{p,q}(\nu',a,b) .$$

We claim that
$$
E_r^{p,q}(a, b) = 
\lim_{\stackrel \longleftarrow \nu} E_r^{p,q}(\nu,a,b) 
$$
is a spectral sequence which converges to
$FH^*_{[\epsilon_0 a,  \epsilon_0 b]}(\omega_{\epsilon_0})$. 

In order to show that $E_r^{p,q}(a,b)$ is a spectral sequence we
need to show that every exact sequence 
$$\xymatrix{0 \ar[r] & \tx{im} \, d_r^{p-r,q+r-1}(\nu) \ar[r] & \ker \, 
  d_r^{p,q}(\nu) \ar[r] & E_{r+1}^{p,q}(\nu) \ar[r] & 0
}$$
remains exact in the inverse limit. This is the point where we need to
use field coefficients. In this situation the inverse limit is an exact
functor if all terms involved in the inverse systems are finite
dimensional vector spaces~\cite[Theorem~VIII,5.7]{ES52}, which is the
case in our setting.  

The limit of the spectral sequence $E_r^{p,q}(a,b)$ is
$\displaystyle \lim_{\stackrel \longleftarrow \nu} FH^*_{[\epsilon_\nu
  a,\epsilon_\nu b]}(Y_\nu^{\epsilon_\nu},J_\nu)$ by definition. Since
the vertical arrows in the diagram induce canonical
isomorphisms in homology we obtain canonical isomorphisms
$$
FH^*_{[\epsilon_\nu a,\epsilon_\nu b]}(Y_\nu^{\epsilon_\nu},J_\nu)
\simeq FH^*_{[\epsilon_1 a,\epsilon_1 b]}(\epsilon_1
K_\nu,\omega_{\epsilon_1})
$$ 
and therefore  
\begin{equation*}
  \lim_{\stackrel \longleftarrow \nu} FH^*_{[\epsilon_\nu a,\epsilon_\nu
  b]}(Y_\nu^{\epsilon_\nu},J_\nu)  = \lim_{\stackrel \longleftarrow \nu}
FH^*_{[\epsilon_1 a,\epsilon_1
  b]}(\epsilon_1 K_\nu,\omega_{\epsilon_1}) = FH^*_{[\epsilon_1 a,\epsilon_1
  b]}(\omega_{\epsilon_1}).
\end{equation*}
The last equality holds by definition. On the other hand we have
  $FH^*_{[\epsilon_1 a,\epsilon_1 b]}(\omega_{\epsilon_1}) \simeq
  FH^*_{[\epsilon_0 a,\epsilon_0 b]}(\omega_{\epsilon_0})$ by
  the deformation invariance Theorem~\ref{thm:invariance under def of
  symp form}, and the claim is proved.  

\medskip 

The spectral sequences $E_r^{p,q}(a,  b)$ are functorial with
respect to the truncation of the action morphisms since all the
intermediate complexes and maps involved in the construction have this
property. 

\medskip 

The definition of $E_2^{p,q}(a,b)$ and that of $\mc{FH}^q_{[a,b]}(\widehat
F)$ given in Remark~\ref{rmk:Flocsys} directly imply that 
$$
E_2^{p,q}(a,b)=H^p(B;\mc{FH}^q(\widehat F)). 
$$

\medskip 

The last assertion of the theorem states that the spectral
sequence $E_r^{p,q}(a,b)$, $r\ge 2$ and the local system
$\mc{FH}^q_{[a,b]}(\widehat F)$ are canonical, i.e. they do not
depend on the various choices involved (the almost complex structures,
the horizontal distribution used to lift $X_f$, the function $f$, the
perturbations $G$ etc.)  This is proved by the continuation method in
Floer homology, i.e. by considering suitable $s$-dependent
interpolating families. The key point is that all constructions can be
performed so that the resulting morphisms preserve the filtrations,
thus inducing isomorphisms between the corresponding spectral
sequences. A glimpse of this phenomenon has already appeared in
Lemmas~\ref{lem:filtered chain morphisms with small
  epsilon}\,--\,\ref{lem:filtered monotonicity morphisms} and we omit
further details.  
\hfill{$\square$}

\medskip 

In the proof that follows the notation $0^-$, $0^+$ stands for a small
enough negative, respectively positive number. 

\medskip 

\noindent {\it \small Proof of Theorem~B.}
Let us recall from Section~\ref{sec:the geometric Hamiltonians} that,
given $\nu\in\NN$, the Hamiltonian $K_\nu$ and the vector field
$Y^\epsilon_\nu$ are of 
the form 
$$
K_\nu=h_\nu+c_\nu\widetilde f +\varphi_\nu +G_\nu, \qquad 
Y_\nu^\epsilon = J_\nu(X_{h_\nu} + \epsilon c_\nu\widetilde{X_f} +
X_{\varphi_\nu}^{\tx{v}} + X_{G_\nu}^{\tx{v}}),
$$ 
where $\varphi_\nu$ is a
$C^2$-small time-independent perturbation of $h_\nu+c_\nu\widetilde f$
supported in $E\setminus \{1-\frac \delta 4 \le
S\le 1\}$ for some $\delta>0$, and $G_\nu$ is a small time-dependent
perturbation supported in an arbitrarily small neighbourhood of
$\partial E$. We use the notation $\widetilde
K_\nu=h_\nu+c_\nu\widetilde f+\varphi_\nu$ and $E_{1-\frac \delta 4} =
E\setminus \{1-\frac \delta 4 \le S\le 1\}$. 

Let us choose $\mu>0$ smaller than $\min(\tx{Spec}(\partial E))$.
A deformation argument shows that, as $\nu\in\NN$ varies, the
  spectral sequences $E_r^{p,q}(\nu,0^-,\mu)$ are canonically
  isomorphic at the page $r=1$, and therefore at all pages $r\ge
  1$. It is thus enough to prove that  
  $E_r^{p,q}(\nu,0^-,\mu)$ is canonically isomorphic to the
  Leray-Serre spectral sequence under the assumption that the maximal
  slope of $h_\nu$ is smaller than $\min(\tx{Spec}(\partial E))$. 
This implies both assertions of Theorem~B, namely the isomorphism of local
  systems $\mc{FH}^q_{[0^-, 0^+]}(\widehat F)\simeq \mc
  H^{k+q}(F,\partial F)$ and the isomorphism of spectral sequences 
  $E_r^{p,q}(0^-,0^+)\simeq {_{_{LS}}}E_r^{n+p,k+q}$, $r\ge 2$. 

We claim that, for a choice of $\nu$ as above, the complex
$FC^*_{[0^-,\epsilon\mu]}(\epsilon,\nu)$ coincides with the Morse 
complex of $Y_\nu^\epsilon$ for
$0<\epsilon\le\epsilon_\nu$ if $c_\nu$, $|\varphi_\nu|$ and $|G_\nu|$
are small enough. Note that $Y_\nu^\epsilon$ is a
pseudo-gradient for the Morse function $\widetilde K_\nu$ on $E$. To
prove the claim we must show that Floer trajectories connecting
critical points of $\widetilde K_\nu$ are independent of time. Since
the Hamiltonian $\widetilde K_\nu$ is $C^2$-small on $E_{1-\frac
  \delta 4}$ and $E$ is symplectically aspherical, it is enough to show
that these Floer trajectories are
contained in $E_{1-\frac \delta 4}$. Let us argue
by contradiction and assume that this is false. Then we find an
$\epsilon>0$, sequences $c_n\to 0$, $|\varphi_n|\to 0$, $|G_n|\to
0$ and $t_n\in \Ss^1$, $n\ge 1$, 
as well as a sequence $u_n$ of maps solving the equation $\partial _s
u_n+J_\nu \partial_t u_n = Y_n^\epsilon\circ u_n$ for 
$$
Y_n^\epsilon = J_\nu(X_{h_\nu} + \epsilon c_n\widetilde{X_f} +
X_{\varphi_n}^{\tx{v}}+ X_{G_n}^{\tx{v}}),
$$
satisfying $u_n(s,\cdot)\to x^{\pm}$, $s\to\pm\infty$ for some distinct
points $x^\pm\in\mathrm{Crit}(\widetilde K_\nu)$, as well as
$$
u_n(0,t_n)\notin E_{1-\frac \delta 4}.
$$ 
The uniform pseudo-gradient property for $\mc Y_n^\epsilon=J_\nu \frac
d {dt}- Y_n^\epsilon$ provides
a uniform bound on the energy $E(u_n)$ for all $n$, whereas a
uniform $C^0$-bound on the sequence $u_n$ follows from
Theorem~\ref{estimation C 0}. By Floer-Gromov
compactness we obtain in the 
limit a nonconstant map $u$ satisfying $u_t+J_\nu u_s=J_\nu
X_{h_\nu}$ and $u(s,\cdot)\to y^{\pm}$ for two distinct
points $y^{\pm} \in E_{1-\frac \delta 4}$, and such that $u(0,t)\notin
E_{1-\frac \delta 4}$ for some $t\in\Ss^1$. Since $h_\nu=h_\nu(S)$ we
can apply the maximum principle as in~\cite[Lemma~1.8]{Viterbo99} and
show that $u$ cannot have a local maximum in $\partial E\times
[1-\frac \delta 4,\infty[$ (see also Remark~\ref{rmk:split
  hamiltonians}). The map $u$ must therefore be entirely contained in
$E_{1-\frac \delta 4}$, a contradiction.  

 Using the notation $E_r^{p,q}(\nu,a,b)$ for the spectral sequence
 associated to the filtered complex $FC^*_{[\epsilon_\nu a,\epsilon_\nu
 b]}(Y_\nu^{\epsilon_\nu},J_\nu)$, we must prove that we have a
 canonical isomorphism of spectral sequences 
 \begin{equation} \label{eq:iso_LS}
 E_r^{p,q}(\nu,0^-,\mu) \simeq {_{_{LS}}}E_r^{n+p,k+q}. 
 \end{equation} 
  At this point we need to recall the construction of the Leray-Serre
 spectral sequence using cellular homology, as defined in
 Section~\ref{sec:MHlocsyst}. Assume one has CW-decompositions
 of $E$ and $B$ with the property that the projection of a cell
 $e^k_i$ in $E$ is a cell $\pi(e^k_i)$ in $B$. The cellular complex 
 $$
 \mathrm{Cell}^k(E)=\bigoplus_i \ZZ\langle e^k_i\rangle
 $$
 is then naturally filtered by
 \begin{equation} \label{eq:filteredcell} 
 F_p\mathrm{Cell}^k(E)=\bigoplus_{\dim\,\pi(e^k_i)\ge p} \ZZ\langle
 e^k_i\rangle, 
 \end{equation} 
 and the spectral sequence associated to this filtration is 
 ${_{_{LS}}}E_r^{p,q}$. This is precisely the definition given
 in~\cite{McC} tailored to the setup of cellular cohomology. 

 The isomorphism~\eqref{eq:iso_LS} follows now from the fact that the
 underlying filtered complexes are isomorphic. Indeed, the
 complex $FC^*_{[\epsilon_\nu a,\epsilon_\nu 
 b]}(Y_\nu^{\epsilon_\nu},J_\nu)$ was shown to coincide with the Morse
 complex, the latter is tautologically identified with the cellular
 complex on $E$ by our standing assumption on the behaviour of
 $K_\nu$ near its critical points (cf. Section~\ref{computation of the E
 1 term}), and the filtration~\eqref{eq:filtration in constr of sp seq}
 on the Morse complex tautologically coincides with the
 filtration~\eqref{eq:filteredcell} on the cellular complex due to our
 standing assumption on the behaviour of $f$ near its critical points.  
 The shift in the grading comes from the fact that the Floer
 complex is graded by minus the Conley-Zehnder index instead of the
 Morse index. 
\hfill{$\square$}

\medskip 

{\it Acknowledgements.} I am indebted to
Claude Viterbo for having suggested the topic to me, for his
confidence and for his constant enthusiastic support. I am equally
indebted to Dietmar Salamon, from whom I learnt that finding the right
degree of generality for a problem represents half the way to its
solution. Work on this paper has spread over a long period, during
which the author was successively supported by Universit\'e Paris Sud,
\'Ecole Polytechnique (Palaiseau), \'Ecole Normale Sup\'erieure de
Lyon, ETH Z\"urich, and Universit\'e Louis Pasteur, Strasbourg.

 
\appendix

\section{On symplectic forms and the taming property}
\label{sec:appendix symplectic forms}


\begin{prop} 

Let $\omega_1$ and $\omega_2$ be two symplectic forms on the closed
manifold $B$. We denote $\mc{J}(B)$ the set of (smooth) almost complex
structures on $B$ and $\mc J_\tau(B,   \omega_i) \subset \mc J (B)$
the contractible open subset of almost complex structures that are tamed by
$\omega_i$, $i=1,2$. 

\begin{enumerate} 
 \item Assume $\mc J_\tau(B,   \omega_1) \subset \mc J_\tau(B,  
 \omega_2)$. Then $\mc J_\tau(B,   \omega_1) = \mc J_\tau(B,  
 \omega_2)$ and $\omega_1 = f\omega_2$, with $f$ a strictly positive
 function on $B$. If $\dim \, B \ge 4$ the function $f$ is constant on
 every component of $B$. 
 \item \label{second assertion in prop relating a c str to symp forms}
   The previous statement holds with $\mc J_\tau (B,   \omega_i)$
 replaced by $\mc J(B,   \omega_i)$, the set of all
 $\omega_i$-compatible almost complex structures on $B$. 
\end{enumerate} 

\end{prop} 

\demo As $\mc J(B,  \omega_i) \subset \mc
J_\tau(B,  \omega_i)$ it is enough to prove~(\ref{second assertion in prop relating a
  c str to symp forms}). We can work fiberwise in each tangent space
and, without loss of generality, assume that $\omega_1$ and $\omega_2$
are linear symplectic forms in $\RR^{2n}$. We denote $\mc
J(\omega_i)=\mc J(\RR^{2n},\omega_i)$, $i=1,2$. The statement is obvious
for $n=1$ and we assume $n\ge 2$. Fix $v,w\in \RR^{2n}$ which are 
noncolinear. We
claim that $\omega_1(v,  w)=0$ if and only if $\omega_2(v,  w)=0$.
Assume first $\omega_1(v,  w)=0$. There exists
$J\in \mc J(\omega_1)$ such that $V=\tx{Sp}\langle v,  Jv
\rangle$ and $W=\tx{Sp}\langle w,  Jw\rangle$ are orthogonal with
respect to $\omega_1$. Let $T=\tx{Sp}\langle v,  Jv,  w,  Jw
\rangle$. By hypothesis the space $V$ is symplectic for
$\omega_2$ and its orthogonal $V^{\perp_{\omega_2}}$ in $T$, denoted $V'$, is a
complement of $V$ and therefore generated by $w+av+Jv$ and
$Jw+cv+dJv$, where $a,b,c,d\in\RR$ are suitable constants. We show that
$a=b=c=d=0$. Assume for example $d\neq 0$. We have $\omega_1(v, 
Jw+cv+dJv)=d\cdot \omega_1(v,  Jv) \neq 0$. There exists therefore $J'\in
\mc J (\omega_1)$ such that $Jw+cv+dJv=\tx{sign}(d) \cdot J'v$. Then
$\omega_2(v,  Jw+cv+dJv)=\tx{sign}(d)\cdot \omega_2(v,  J'v) \neq
0$, which contradicts the definition of $V'$. One proves in the same
way that $a=b=c=0$. Assume now that $\omega_2(v,  w)=0$. If
$\omega_1(v,  w)\neq 0$, there exists $J\in \mc J(\omega_1)$ such
that $w=\tx{sign}(\omega_1(v,  w))\cdot Jv$, hence $\omega_2(v, 
w)=\pm \omega_2(v,  Jv)\neq 0$, which is again a contradiction.  
One can easily see now that we also have
$\omega_1(v,  w)>0$ if and only if $\omega_2(v,  w)>0$. 

Let now $e_1,f_1,\ldots,e_n,f_n$ be a symplectic basis for $\omega_1$
i.e. $\omega_1=\sum_{i=1}^n e_i^*\wedge f_i^*$. The above implies that
$\omega_2 = \sum_{i=1}^n \lambda_ie_i^*\wedge f_i^*$, $\lambda_i>0$.
We show that all the $\lambda_i$ are equal: for $1\le i,j\le n$ we have
$\omega_1(e_i+e_j,  f_i-f_j)=0$, hence $0=\omega_2(e_i+e_j, 
f_i-f_j) = \lambda_i-\lambda_j$. Therefore $\omega_1=\lambda\omega_2$,
$\lambda>0$ and $\mc J(\omega_1)=\mc J(\omega_2)$, $\mc
J_\tau(\omega_1)=\mc J_\tau(\omega_2)$. 

On the manifold $B$ we can write $\omega_1 = f\omega_2$, with
$f:B\longrightarrow \RR^+$ smooth. Assume $\dim \, B\ge 4$. By closedness
we get $df\wedge \omega_2=0$. Assume $df$ is not identically zero i.e.
there exist $x\in B$ and $X\in T_xB$ such that $df_x(X)>0$. The germ
of $\Sigma = f^{-1}(f(x))$ at $x$ is a smooth hypersurface of
dimension at least $3$, on which $\lambda=\iota_X\omega_2$ is nonzero.
The form $\omega_2$ is nondegenerate on $\ker \,  \lambda$ and this
implies $(df\wedge \omega_2)_x\neq 0$, which is a contradiction. 
\hfill{$\square$}

\bibliographystyle{abbrv}
\bibliography{LS}

\end{document}